  \documentclass[runningheads]{llncs}
  \usepackage[T1]{fontenc}
  \usepackage{graphicx}
  \usepackage{cite}
  \usepackage{amsmath,amssymb,amsfonts}
  \usepackage{algorithmic}
  \usepackage{textcomp}
  \usepackage{xcolor}
  \usepackage{dsfont}
  \usepackage[ruled,linesnumbered]{algorithm2e}
  \usepackage{float}
  \usepackage{enumerate}
  \usepackage{epstopdf}
  \usepackage{bbding}
  \usepackage{pifont}
  \usepackage{subfigure}
  \usepackage{epsfig}

\begin{document}	
  	\title{Accelerating Stochastic Recursive and Semi-stochastic Gradient Methods with
   Adaptive Barzilai-Borwein Step Sizes}
  	\titlerunning{AR-VR-SGD-RHBB}
  	\author{JiangShan Wang, YiMing Yang,  Zheng  Peng}
  	% \authorrunning{J. Long, C. Wang, C. Ou et al.}
  	\institute{{School of Mathematics and Computational Science, Xiangtan University}\\ {Xiangtan, Hunan 411105, China}}%\\
  	%\email{longjiangshan@whu.edu.cn \\wchenxu@whu.edu.cn \\ouchanghai@whu.edu.cn \\tangming@whu.edu.cn}
  	\maketitle              % typeset the header of the contribution
  	\begin{abstract}
  		The mini-batch versions of StochAstic Recursive grAdient algoritHm and Semi-Stochastic Gradient Descent method, employed the random Barzilai-Borwein step sizes (shorted as MB-SARAH-RBB and mS2GD-RBB), have surged into prominence through timely step size sequence. Inspired by modern adaptors and variance reduction techniques, we propose two new variant rules in the paper, referred to as RHBB and RHBB+, thereby leading to four algorithms MB-SARAH-RHBB, MB-SARAH-RHBB+, mS2GD-RHBB and mS2GD-RHBB+ respectively. RHBB+ is an enhanced version that additionally incorporates the importance sampling technique. They are aggressive in updates, robust in performance and self-adaptive along iterative periods. We analyze the flexible convergence structures and the corresponding complexity bounds in strongly convex cases. Comprehensive tuning guidance is theoretically provided for reference in practical implementations. Experiments show that the proposed methods consistently outperform the original and various state-of-the-art methods on frequently tested data sets. In particular, tests on the RHBB+ verify the efficacy of applying the importance sampling technique to the step size level. Numerous explorations display the promising scalability of our iterative adaptors.
  		\keywords{Variance reduction \and stochastic optimization \and random hedge Barzilai-Borwein method \and importance sampling \and iterative adaptors}
  	\end{abstract}
  	
  	\section{Introduction}
  	We focus on the following stochastic optimization problem
  	$$
  	\min _{w \in \mathbb{R}^d} P(w)=\mathbb{E}_{\mathbb{P}}[f(w ; \xi)_{\xi\in \Omega}],
  	$$
  	where $\xi$ is a random instance of an input-output pair $(x_i, z_i)$, with input representation vector $x_i$ and target output $z_i$. Hence, $f(\cdot)$ takes the form 
  	$$f(w ; \xi)=f\left(w ; (x_i, z_i)\right).$$
  	Define $P(\cdot)$ by empirical expectation on probability space $(\Omega, \mathbb{P})$, where  $\Omega=\left\{\xi_1, \cdots, \xi_{n}\right\}$ is a finite support set and $\mathbb{P}$ is the probability measure over $\Omega$. In this case, it transforms into an unconstrained finite-sum model, i.e.,
  	\begin{equation}
  	\label{1}
  	\min _{w \in \mathbb{R}^d} P(w)=\int_{\Omega} f\left(w ; (x_i, z_i)\right) d \mathbb{P}(x_i,z_i)\approx  \frac{1}{n} \sum_{0\le i\le |\Omega|} f(w; (x_i, z_i)).
  	\end{equation}
  	
  	Problem (\ref{1}) covers a broad range of applications. Assume that a lower bound $P(w_*)$ of $P(\cdot)$ exists, the goal is to iteratively update $w$ to reduce $P(\cdot)$ steadily and swiftly. Given a sequence of $n$ labeled pairs $\{\left(x_1, z_1\right), \ldots,\left(x_{n}, z_{n}\right)\}$ into finite dimensional spaces $\{(\mathbb{R}^d, \mathbb{R})\}$, the linear least squares regression is of $f(w) \stackrel{\text { def }}{=}$ $\left(x_i^T w-z_i\right)^2$. In terms of the logistic regression, we exploit $f(w) \stackrel{\text { def }}{=} \log \left(1+\exp \left(-z_i x_i^T w\right)\right)$. Indeed, efficient regularization may be taken into account for specific purposes, then it develops into a composite model, i.e.,
  	$$\min _{w \in \mathbb{R}^d} \frac{1}{n} \sum_{0\le i\le |\Omega|} f(w; (x_i, z_i))+\psi(w),$$
  	where $\psi(\cdot)$ is a proper, closed and convex penalty on parameters. Its subdifferential at $w$ defines $\partial \psi(w)=\left\{v \in \mathbb{R}^d \mid \psi(d) \geq \psi(w)+ v^T \left(d-w\right), \forall d \in \operatorname{dom}(\psi)\right\}$. Throughout the paper, we mainly utilize a smoothing regularization process, thus $\psi(\cdot)$ is differentiable and $\partial \psi(w)=\{\nabla \psi(w)\}$. 
  	
  	Due to the productive and scalable frameworks, stochastic optimization is prevalent among large-scale problems or complex relationships. It offers a cost-effective alternative to deterministic schemes, stated in a universal form as
  	$$w^{(t)}=w^{(t-1)}-\eta_t g_t\left(w^{(t-1)}, \nu_t\right).$$ 
   Here, $g(\cdot)$ represents the gradient estimator, $\nu$ the randomness pointer. Vanilla stochastic gradient descent (SGD) \cite{Robbins}\cite{shalev} specifies $g(\cdot)$ into the basic moment estimate, where it enjoys cheap computational cost per update and the independence with $n$ in term of complexity. However, this straightforward estimator inevitably introduces variance and noise to the steps. Diminishing step size $\eta_t=O(1 / t)$ is then forced to employed in \cite{luo}\cite{solodov}\cite{nemirovski} for a sublinear convergence rate of $O(1 / t)$ (Moulines et al.\cite{moulines}), which should satisfy
  	$$
  	\sum_{t=1}^{\infty} \eta_t=\infty \text {  and  } \sum_{t=1}^{\infty} \eta_t^2<\infty.
  	$$
  	Such updates come with a side effect of halt in convergence near the eventual limit\cite{eon}. Batch methods of \cite{dekel}\cite{cotter}\cite{shalev} decrease the intrinsic variance through a bunch of samples, but at the cost of further computational workload. Therefore, parallel processing power becomes indispensable for explosion-scale data.
  	
  	Shifts on $g(\cdot)$ have been intensively studied with a vast number of research papers. Canonical examples include but not limit to SAG/SAGA \cite{Schmidt}\cite{defazio}, SVRG \cite{johnson}, MSVRG\cite{reddi}, S2GD \cite{konevcny2}, mS2GD \cite{konevcny}, MISO \cite{mairal}, SARAH\cite{nguyen2}, iSARAH \cite{nguyen3}, MB-SARAH \cite{nguyen2}, SPIDER \cite{fang}, etc. Konecny et al. \cite{konevcny} proposed mS2GD and showed it reaches a pre-defined accuracy with less overall work than a method without batching. They established a threshold for the batch size, at which more than linear speedup can be achieved, it's worthy to further explore. MB-SARAH is presented by Nguyen et al. \cite{nguyen2} for solving non-convex problems. Recursive updates free from the storage of past gradients and avoid oscillation of the Euclidean norm of $g(\cdot)$ in inner loops, which are well-suited for modern complex scenarios. Researches have also extended to the dual space of (\ref{1}) by updating random dual variables or variable blocks, such as RCDM\cite{nesterov2}, AsySPDC\cite{liu2}, SDCA\cite{shalev2}, SPDC\cite{zhang}, mSDCA\cite{takac}, ASDCA\cite{shalev3} and QUARTZ\cite{zheng}, Prox-SDCA \cite{zhao}. The subsequent drawback lies in that these algorithms rely on a tuning step size by hand, which can be time consuming in practice.
  
  	Several methods of auto-tuning prevail among the stochastic algorithms. Barzilai-Borwein method (BB) of the second order tuning is outstanding in the trend, due to its simplicity and numerical efficiency. Sopyła et al.  \cite{sopyla} employed the BB in vanilla SGD to solve linear SVM in dual space. Tan et al. \cite{tan} incorporated the BB into SVRG (SVRG-BB) and established the linear convergence in strongly convex cases. To further accelerate the rates, Yang et al.  \cite{yang3} introduced the BB to mS2GD (mS2GD-BB) for nonsmooth and strongly convex functions. On the basis of Hessian and its eigenvalues, Ma et al.  \cite{ma} proposed the stabilized Barzilai-Borwein (SBB) method to match SVRG (SVRG-SBB) for the ordinal embedding problems, which avoids the denominator tending to zero. Yang et al.  \cite{yangy} considered the inexact SARAH (iSARAH-BB) in order to reduce the cost in deterministic steps, then showed its robustness in implementations. Byrd et al. \cite{byrd} utilized batch methods to approximate the quasi-Newton property. Recently, the timely random Barzilai-Borwein method (RBB) emerged and was primarily applied in MB-SARAH (MB-SARAH-RBB\cite{yang}) and mS2GD (mS2GD-RBB\cite{yang2}) algorithmic settings. The promising performance outperformed and matched state-of-the-art algorithms. However, they still have flaws and aspects that can be explored and improved.
  	
  	In the context of the RBB rule, when the batch size is insufficient, it can result in a high level of noise, primarily caused by the random step sizes. This noise leads to an increasing or oscillating trend towards divergence. 
   %On the other hand, increasing the batch size reduces the variance of stochastic curvature but unfortunately also slows down convergence, which is discouraging.
   %In the RBB rule, an insufficient batch may lead to steps with a high level of noise (most noise from step sizes), which incurs an increasing trend or oscillating trend towards divergence. 
   As batch size increases, the variance of stochastic curvature decreases, but using a larger batch discouragingly slows down the convergence. In addition, it's insensitive to the iterative periods, and the well-worn uniform sampling deserves to be generalized. Therefore, we introduce the random hedge Barzilai-Borwein method (RHBB) in pursuit of improvement. We further incorporate our RHBB with the importance sampling technique and develop another enhanced version, RHBB+.

   \vspace{3pt}  
  	The key contributions in the paper are summarized as follows:
  	\begin{enumerate}
  		\item[1)] We propose the random hedge Barzilai-Borwein step size rule (RHBB) for MB-SARAH and mS2GD algorithmic settings, obtaining MB-SARAH-RHBB and mS2GD-RHBB algorithms. The adaptive acceleration mechanism is analyzed, trade-off rules are studied and the tuning guidance is provided.
  		\item[2)] We incorporate the importance sampling technique into the RHBB to make use of the distributed efficiency in data sets (e.g. the elementwise sparsity), which yields the enhanced version RHBB+ and the corresponding MB-SARAH-RHBB+ and mS2GD-RHBB+ algorithms.
  		\item[3)] We establish the global convergence of all four algorithms in strongly convex cases. Further, we analyze the ill-conditioned scenario to theoretically validate the robustness of the new algorithms. In MB-SARAH-RHBB/RHBB+, the square of the full gradient converges linearly in expectation. And the expected distance of iterates to the global optimum has linear convergence in mS2GD-RHBB/RHBB+. 
  		\item[4)] We conduct extensive experiments to demonstrate the exceptional performance of our algorithms. Next, we explore a tentative, incremental scheme for the iterative adaptor and view immense potential for its scalability.
  	\end{enumerate}
  	
  	\section{Common Assumptions and Inequalities}
  	We add subscripts to distinguish element functions, e.g. $f_i$ denotes the $i$-th component. Unless otherwise specified, $\|\cdot\|$ denotes Euclidean norm in this paper. Subsequently, we provide following common assumptions.\\
  	\textbf{Assumption 1} ({\bf Smoothness}). Each convex $f_i: \mathbb{R}^d \rightarrow \mathbb{R}$ is $L$-smooth over any compact set of its domain, i.e., there exists an $L>0$, for all $w, w^{\prime} \in \text{dom}(f_i)$
  	$$
  	\left\|\nabla f_i(w)-\nabla f_i\left(w^{\prime}\right)\right\| \leq L\left\|w-w^{\prime}\right\|.
  	$$
  	
  	Due to $\nabla f_i(w)-\nabla f_i\left(w^{\prime}\right)=H\left(w-w^{\prime}\right)$ where $H=\nabla^2 f_i\left(\hat{w}\right)$ is the Hessian at midpoint $\hat{w}$, the largest eigenvalue of $H$ is no more than $L$. In view of it, we obtain a boundary for the square distance between the gradients, i.e.,
  	\begin{equation}
  	\label{L}
  	\begin{aligned}
  	\left\|\nabla f_i(w)-\nabla f_i\left(w^{\prime}\right)\right\|^2
  	& =\left(w-w^{\prime}\right)^{\top} H^2\left(w-w^{\prime}\right) \\
  	&\leq  L\left(w-w^{\prime}\right)^{\top} H\left(w-w^{\prime}\right) \\
  	&=  L\left(\nabla f_i(w)-\nabla f_i\left(w^{\prime}\right)\right)^{\top}\left(w-w^{\prime}\right).
  	\end{aligned}
  	\end{equation}
  	
  	The individual $L$-smoothness implies the overall $P(w)=\frac{1}{n} \sum_{0\le i\le |\Omega|} f_i(w)$ is $L$-smooth as well. Equivalently, we derive the following bound 
  	\begin{equation}
  	\label{inequation1}
  	P(w) \leq P(w^{\prime})+\nabla P(w^{\prime})^T(w-w^{\prime})+\frac{L}{2}\|w-w^{\prime}\|^2.
  	\end{equation}
  	For the $L$-smoothness of the overall $P(\cdot)$, we have another equivalent claim, i.e., 
  	\begin{equation}
  	\label{L2}
  	P(w) \geq P(w^{\prime})+\nabla P(w)^T(w-w^{\prime})+\frac{1}{2 L}\|\nabla P(w)-\nabla P(w^{\prime})\|^2.
  	\end{equation}
  	\textbf{Assumption 2a} ({\bf Strong convexity I}). $P(w)$ is $\mu$-strongly convex, i.e., there exists $\mu>0$ such that, for all $w, w^{\prime} \in \text{dom}(P)$
  	\begin{equation}
  	\label{C}
  	(\nabla P(w)-\nabla P(w^{\prime}))^T(w-w^{\prime}) \geq \mu\|w-w^{\prime}\|^2.
  	\end{equation}
  	\textbf{Assumption 2b} ({\bf Strong convexity II}). Each component $f_i$ is $\mu$-strongly convex, i.e., there exists $\mu>0$ for each $f_i$ such that for all $w, w^{\prime} \in \text{dom}(f_i)$\\
  	\begin{equation}
  	\label{C2}
  	(\nabla f_i(w)-\nabla f_i(w^{\prime}))^T(w-w^{\prime}) \geq \mu\|w-w^{\prime}\|^2.
  	\end{equation}
  	
  	Assumption 2b can imply Assumption 2a, but not vice versa. Assumption 2b is a stricter premise which requires the strong convexity on each $f_i$.
  	
  	For the $\mu$-strongly convex $P(\cdot)$, we have another equivalent claim as
  	$$
  	P(w) \geq P(w^{\prime})+\nabla P(w^{\prime})^T(w-w^{\prime})+\frac{\mu}{2}\|w-w^{\prime}\|^2 .
  	$$
  	
  	Define the global optimum $w_*=\operatorname{argmin}_w P(w)$, it further indicates that
  	\begin{equation}
  	\label{convex1}
  	2 \mu\left[P(w)-P\left(w_*\right)\right] \leq\|\nabla P(w)\|^2, \forall w \in \mathbb{R}^d.
  	\end{equation}
  	To see this, we have by strong convexity
  	$$P\left(w_*\right) \geq P(w)+\nabla P(w)^T\left(w_*-w\right)+\frac{\mu}{2}\left\|w-w_*\right\|^2.$$
  	Then, through some basic derivation, we have 
  	$$
  	\begin{aligned}
  	& 2 \mu\left[P(w_*)-P\left(w\right)\right]+\|\nabla P(w)\|^2 \geq 2 \mu \nabla P(w)^T\left(w_*-w\right)+\mu^2 \left\|w-w_*\right\|^2+\|\nabla P(w)\|^2,\\
  	& 2 \mu\left[P(w_*)-P\left(w\right)\right]+\|\nabla P(w)\|^2 \geq \left\|\nabla P(w)+\mu\left(w_*-w\right)\right\|^2 \geq 0,\\
  	& 2 \mu\left[P(w)-P\left(w_*\right)\right] \leq\|\nabla P(w)\|^2.
  	\end{aligned}
  	$$
  	
  	\section{Motivation}
  	\subsection{Barzilai-Borwein Method and the Random Versions}
  	Barzilai-Borwein method, originally developed in the pioneer work of Barzilai and Borwein [1], shows great preeminence in solving nonlinear optimization problems and has been widely improved up to now. 
  	
  	We automatically hope $\theta_k I$ approximates $H_k$ in the $k$-th epoch, where $\theta_k$ denotes the $k$-th step size, $I$ the identity matrix and $H_k$ the inverse of Hessian. To minimize the residual of the secant equations, i.e., $\left\|\left(1 / \theta\right) s_k-y_k\right\|_2^2$ and $\left\|\theta y_k-s_k\right\|_2^2$, we have the following step size solutions respectively
  	$$\theta_k^{\mathrm{BB} 1}=\frac{s_{k}^T s_{k}}{s_{k}^T y_{k}}, \quad
  	\theta_k^{\mathrm{BB} 2}=\frac{s_{k}^T y_{k}}{y_{k}^T y_{k}},$$
  	where $s_{k}=w_{k+1}-w_{k}$, $y_{k}=\nabla f\left(w_{k+1}\right)-\nabla f\left(w_{k}\right)$.
  	
  	Among the well-known Barzilai-Borwein methods, either BB1 or BB2 is expected to be computed at the start of each epoch and employed uniformly during the entire consecutive stochastic stages (see \cite{tan}\cite{yang3}\cite{ma}\cite{yangy} for a brief reference). Yang et al. \cite{yang}\cite{yang2} advocated to calculate the Barzilai-Borwein step size timely by stochastic curvature (batch curvature, similar to Castera et al. \cite{castera}) within each stochastic stage. In mS2GD algorithmic setting, they proposed a random version of BB1, we instinctively deduce the random BB2 by analogy
  	$$
  	\tilde{\eta}_k^{\mathrm{RBB} 1}=\frac{1}{|S_1|}\cdot \frac{\left\|w_k-w_{k-1}\right\|^2}{\left(\left(w_k-w_{k-1}\right)^T\left(\nabla P_{S_1}\left(w_k\right)-\nabla P_{S_1}\left(w_{k-1}\right)\right)\right)},
  	$$
  	$$
  	\tilde{\eta}_k^{\mathrm{RBB} 2}=\frac{1}{|S_2|}\cdot \frac{\left(\left(w_k-w_{k-1}\right)^T\left(\nabla P_{S_2}\left(w_k\right)-\nabla P_{S_2}\left(w_{k-1}\right)\right)\right)}{\left\|\nabla P_{S_2}\left(w_k\right)-\nabla P_{S_2}\left(w_{k-1}\right)\right\|^2},
  	$$
  	where $\nabla P_{S_1}\left(w_k\right)=\frac{1}{|S_1|} \sum_{i \in S_1} \nabla f_i\left(w_k\right), \nabla P_{S_2}\left(w_{k}\right)=\frac{1}{|S_2|} \sum_{j \in S_2} \nabla f_j\left(w_{k}\right)$. The $S_1, S_2\subset\{1, \ldots, n\}$ are randomly selected subsets with size $|S_1|$ and $ |S_2|$. As in MB-SARAH algorithmic setting, another parameter $\gamma$ is multiplied to adjust the current RBB for a better Hessian approximation, i.e.,
  	$$
  	\eta_k^{\mathrm{RBB} 1}=\frac{\gamma}{|S_1|}\cdot \frac{\left\|w_k-w_{k-1}\right\|^2}{\left(\left(w_k-w_{k-1}\right)^T\left(\nabla P_{S_1}\left(w_k\right)-\nabla P_{S_1}\left(w_{k-1}\right)\right)\right)},
  	$$
  	$$
  	\eta_k^{\mathrm{RBB} 2}=\frac{\gamma}{|S_2|}\cdot \frac{\left(\left(w_k-w_{k-1}\right)^T\left(\nabla P_{S_2}\left(w_k\right)-\nabla P_{S_2}\left(w_{k-1}\right)\right)\right)}{\left\|\nabla P_{S_2}\left(w_k\right)-\nabla P_{S_2}\left(w_{k-1}\right)\right\|^2}.
  	$$
  	
  	\begin{figure*}[htbp]
  	\centering
  	\subfigure[a8a]
  	{\includegraphics[width=0.327\textwidth]{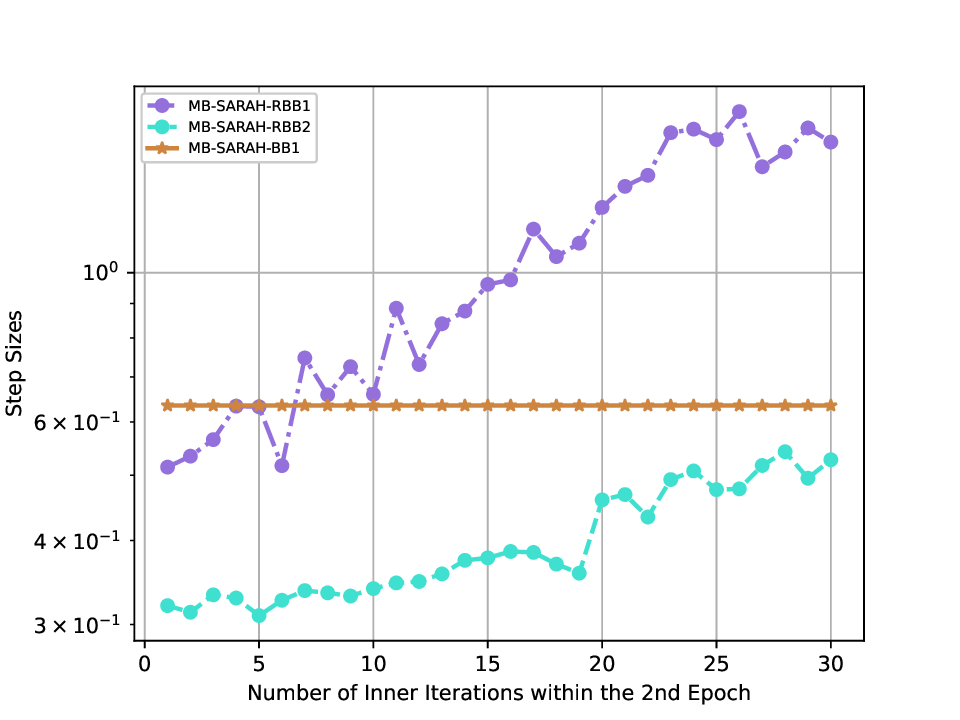}}
  	\subfigure[covtype]
  	{\includegraphics[width=0.327\textwidth]{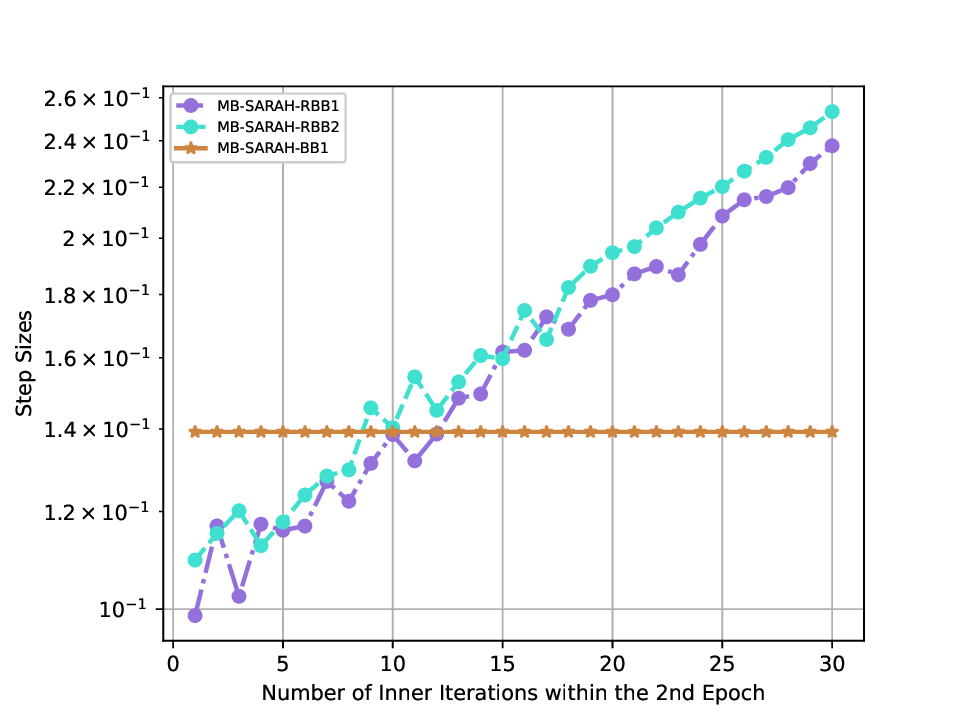}}
  	\subfigure[phishing]
  	{\includegraphics[width=0.327\textwidth]{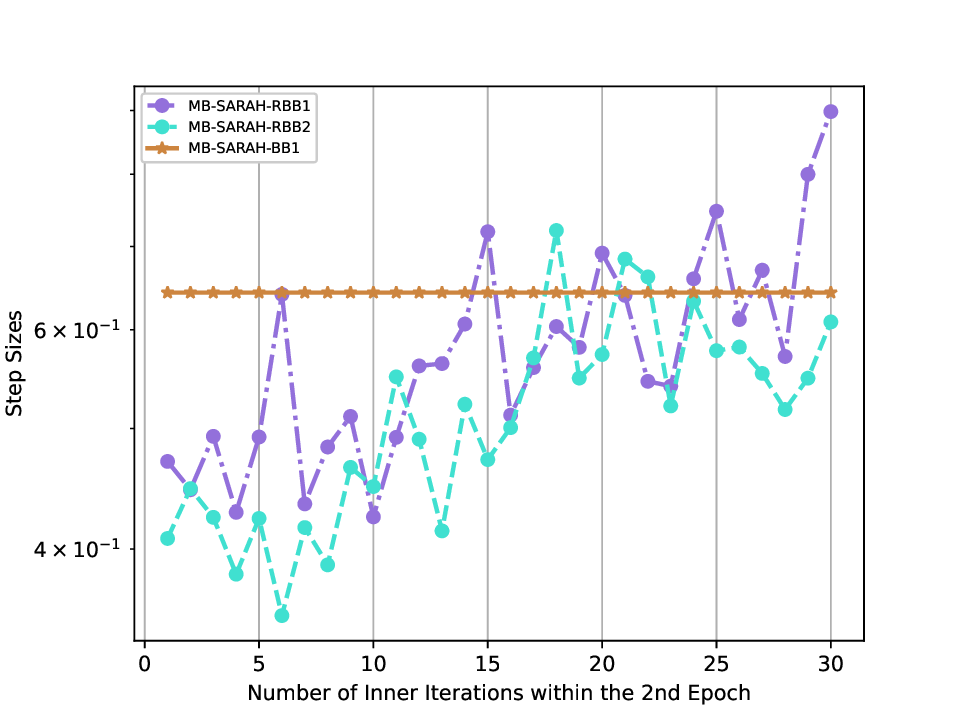}}
  	\subfigure[a8a]
  	{\includegraphics[width=0.327\textwidth]{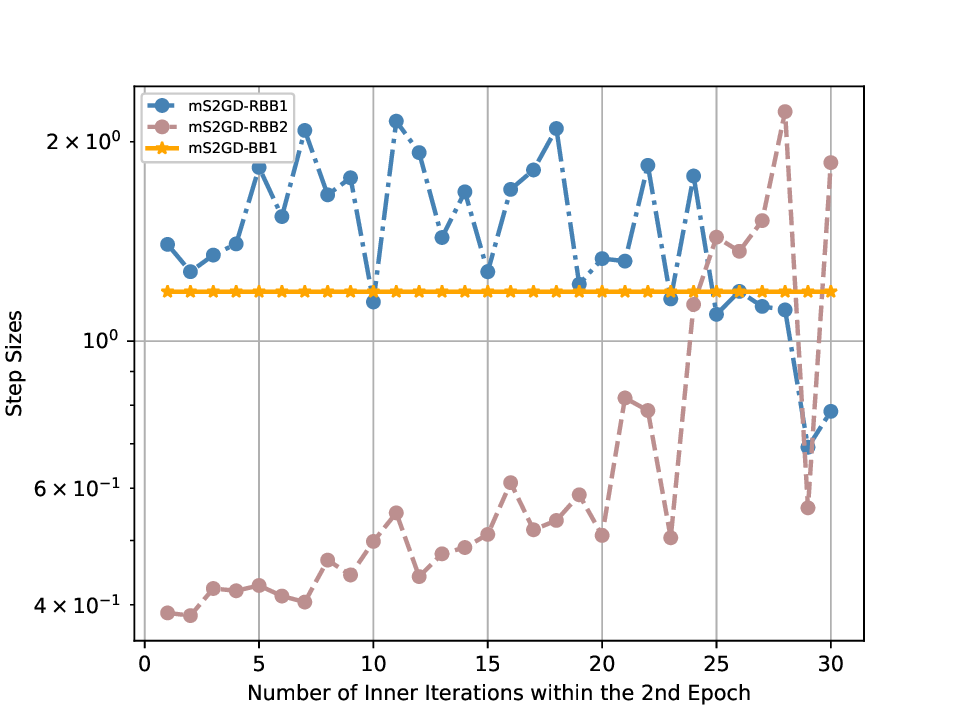}}
  	\subfigure[covtype]
  	{\includegraphics[width=0.327\textwidth]{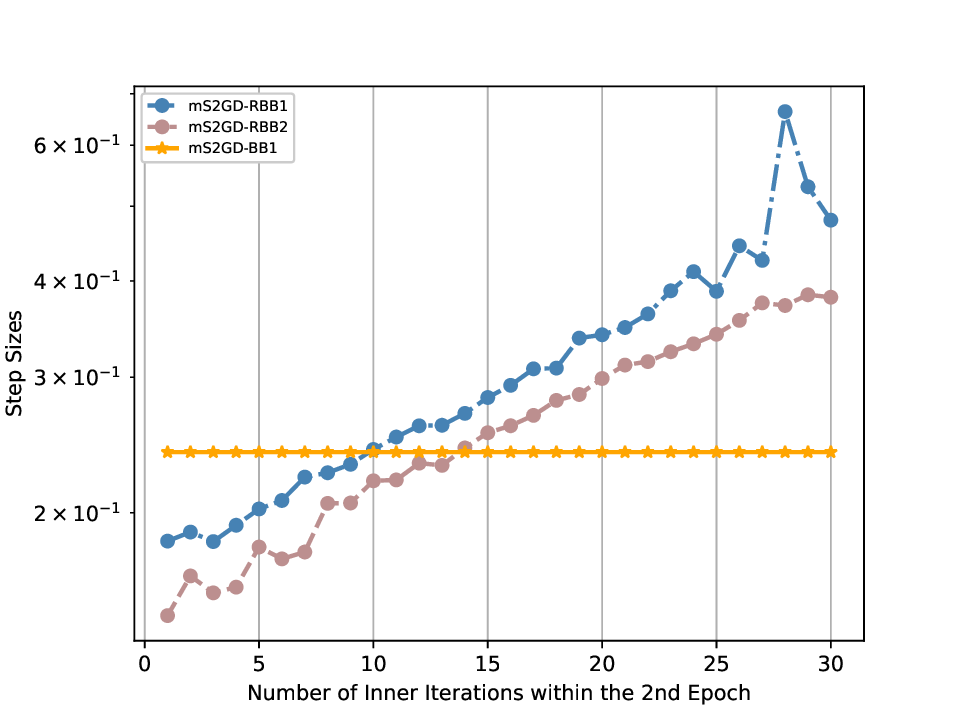}}
  	\subfigure[phishing]
  	{\includegraphics[width=0.327\textwidth]{figures/Fig1_a.eps}}
  	\caption{\footnotesize (a)(b)(c): Trajectories of BB1 step sizes and RBB step sizes in MB-SARAH. (d)(e)(f): Trajectories of BB1 step sizes and RBB step sizes in mS2GD.}
  	\label{fig1}
  	\end{figure*}
  	
  	Let us observe a series of trajectories in Fig. \ref{fig1}, where the $x$-axis denotes the number of inner iterations within the $2$-nd epoch and the $y$-axis represents the corresponding step sizes. All relevant parameters are set as suggested in \cite{yang}\cite{yang2}. For illustration, we discard the first epoch that is insufficient to compute the BB1 step sizes (instead a specified constant is applied during the first epoch). It can be observed that BB1 remains unvaried in a single value throughout the entire epoch, while RBB1 and RBB2 update in a timely manner, evolving with RBB1 taking the precedence or both intertwining mutually. 
  	
  	In fact, the BB1-type is more preferable due to its aggressive finesse, however in a wide range of cases, it still may not reach the extreme of efficiency as well. Note that hedging is an innocuous way to mitigate risks in the financial sector. Inspired by this, we ensure or deliberately expand the effective magnitude of RBB1, while offsetting any over-utility from the opposite direction via its `twin' RBB2 to devise our RHBB step size rule. It can enlarge the adaptive step sizes smoothly and controllably. Further, to improve the adaptivity along iterative periods, we expect to include an adaptor to adjust the step sizes. Specifically, RHBB is based on an affine combination of RBB1, RBB2 via an adaptive parameter $\alpha^{h(\sigma_1s+\sigma_2k)}$, where $\alpha^{h(\sigma_1s+\sigma_2k)}>1$. The adaptor $h$ is an exponential decay rate of the affine magnitude. In early epochs, the $h$ boosts the step sizes at a low cost to accelerate training. As approaching the global optimum, it then enforces them conservative to ensure the final convergence. In fact, the adaptor $h$ is monotone decreasing with respect to the linear indicators $\left(\sigma_1s+\sigma_2k\right)$, and iteratively satisfies
  	$$|\alpha- h(\sigma_1s+\sigma_2k)|>\epsilon(s), \forall k\in m, \forall s.$$
  	In practice, we make $\sigma_1, \sigma_2 \in \{0, 1\}$ and $\epsilon(s)>\frac{s}{m}$, the gap should be expanded in latter periods. It's of great distinction in the structural sense from the composite Barzilai-Borwein method (CBB in \cite{li}) and the composite adaptive Barzilai-Borwein method (CABB in \cite{li}), which utilize a convex combination with parameters within range $(0, 1)$. The CBB and CABB comprise two components, each of which extracts partial resources from either BB1 or BB2. Essentially, all Barzilai-Borwein methods enjoy the `calculation' adaptivity, we attach another adaptor, $h(\cdot)$, to enhance the adaptivity along iterations.
  	
  	The employment of RBB2 introduces another type of stochastic curvature. Notice that we capture the stochastic curvature from two probabilistic subsets $S_1$ and $S_2$, suggesting that we use quite another subset to do hedging. Due to the quasi-Newton property, RHBB indeed reduces the noise in the second order level. For the sake of convergence, we pick the larger batch size as the batch correction. Hence, we have the RHBB for MB-SARAH algorithmic setting as (with total $\alpha^{h(\sigma_1s+\sigma_2k)}>1$):
  	\begin{equation}
  	\label{SARAH-RHBB}
  	\begin{aligned}
  	(\eta_k^s)^{\mathrm{RHBB}}
  	& = \frac{\gamma}{\max\{|S_1|, |S_2|\}}\cdot\left( \frac{\alpha^{h(\sigma_1s+\sigma_2k)}\cdot\left\|w_k^s-w_{k-1}^s\right\|^2}{\left(\left(w_k^s-w_{k-1}^s\right)^T\left(\nabla P_{S_1}\left(w_k^s\right)-\nabla P_{S_1}\left(w_{k-1}^s\right)\right)\right)}\right.\\
  	&+ \left.\frac{\left(1-\alpha^{h(\sigma_1s+\sigma_2k)}\right)\cdot\left(\left(w_k^s-w_{k-1}^s\right)^T\left(\nabla P_{S_2}\left(w_k^s\right)-\nabla P_{S_2}\left(w_{k-1}^s\right)\right)\right)}{\left\|\nabla P_{S_2}\left(w_k^s\right)-\nabla P_{S_2}\left(w_{k-1}^s\right)\right\|^2}\right).
  	\end{aligned}
  	\end{equation}
  	
  	Next, we multiply an extra balance parameter, $\gamma_2$, to the RHBB step sizes in mS2GD algorithmic setting, i.e.,
  	\begin{equation}
  	\label{mS2GD-RHBB}
  	\begin{aligned}
  	(\tilde{\eta}_k^s)^{\mathrm{RHBB}}
  	& = \frac{\gamma_2}{\max\{|S_1|, |S_2|\}}\cdot\left( \frac{\alpha^{h(\sigma_1s+\sigma_2k)}\cdot\left\|w_k^s-w_{k-1}^s\right\|^2}{\left(\left(w_k^s-w_{k-1}^s\right)^T\left(\nabla P_{S_1}\left(w_k^s\right)-\nabla P_{S_1}\left(w_{k-1}^s\right)\right)\right)}\right.\\
  	&+ \left.\frac{\left(1-\alpha^{h(\sigma_1s+\sigma_2k)}\right)\cdot\left(\left(w_k^s-w_{k-1}^s\right)^T\left(\nabla P_{S_2}\left(w_k^s\right)-\nabla P_{S_2}\left(w_{k-1}^s\right)\right)\right)}{\left\|\nabla P_{S_2}\left(w_k^s\right)-\nabla P_{S_2}\left(w_{k-1}^s\right)\right\|^2}\right).
  	\end{aligned}
  	\end{equation}
   In the latter section, we will elaborate the role of $\gamma_2$ in relaxing the stochastic hedge effect. And related trade-off rules will as well be studied.
  
  	We choose the batch scheme for possible opportunities of parallel processing. Notice from Tan et al. \cite{tan} that they applied a convex combination to approximate the full gradient at the snapshots, hence the absolute operation has been taken upon the denominators in the step size calculations. Here in MB-SARAH and mS2GD algorithmic settings, we bear no worry on possible negative step sizes. According to Castera et al.\cite{castera}, when the curvature condition $(g_k^s)^T\nabla^2P(w_k^s)g_k^s$ does not keep positive ($g_k^s$ is an update direction), it's advisable to set $\eta_k^s=c$ ($c>0$). Different from AS in Liu et al. \cite{liu3}, we conduct and use real-time estimations instead of the accumulated moving average.

  	\subsection{Importance Sampling technique}
  	In terms of theory, this technique leads to the improvement of leading constants in the complexity estimates (Richtarik et al \cite{richtarik}, Needell et al \cite{needell}). The overhead associated with configuring distributions and withdrawing samples is negligible, and hence the net effect \cite{csiba} is speedup.
  	
  	Uniform sampling enables unbiased estimators but sacrifices potential opportunities of variance reduction, algorithms nowadays have strived for the opposite. Most sampling techniques have been applied to the gradient estimates, which include but not limit to Prox-SVRG \cite{xiao}, Prox-SDCA \cite{zhao} and SARAH-I \cite{liu}. Hence, we inventively exploit the sampling schemes in the step size level, customizing probability distributions to filter stochastic variations. In fact, we configure $Q \sim$ $\left\{q_1, q_2, \ldots q_n\right\}$ according to Zhao et al \cite{zhao}.
  	
  	The basic moment estimate of $\nabla P(\cdot)$, over the uniform distribution, on the subset $S\subset \Omega$, is in the form of
  	\begin{equation}
  	\label{stochastic estimate}
  	\nabla P_{S}\left(w_k\right)=\frac{1}{|S|} \sum_{0\le i \le |S|} \nabla f_i\left(w_k\right).
  	\end{equation}
  	In light of (\ref{stochastic estimate}), consider that from a general distribution we have
  	\begin{equation}
  	\label{P+}
  	\nabla P_{S}^+\left(w_k\right)=\frac{1}{|S|} \sum_{0\le i \le |S|} \frac{\nabla f_i\left(w_k\right)}{n\cdot q_i} = \frac{1}{|S|} \sum_{0\le i \le |S|} \nabla f^+_i\left(w_k\right).
  	\end{equation}
  	Uniformity delivers $q_1=...=q_{n}=\frac{1}{n}$, suggesting a special case (\ref{stochastic estimate}) of $(\ref{P+})$.
  	
  	By substituting (\ref{stochastic estimate}) with (\ref{P+}), we extend the RHBB (\ref{SARAH-RHBB}) (\ref{mS2GD-RHBB}) to support general distributions, resulting in the enhanced rule RHBB+, i.e.,
  	\begin{equation}
  	\label{SARAH-RHBB+}
  	\begin{aligned}
  	(\eta_k^s)^{\mathrm{RHBB+}}
  	& = \frac{\gamma}{\max\{|S_1|, |S_2|\}}\cdot\left( \frac{\alpha^{h(\sigma_1s+\sigma_2k)}\cdot\left\|w_k^s-w_{k-1}^s\right\|^2}{\left(\left(w_k^s-w_{k-1}^s\right)^T\left(\nabla P_{S_1}^+\left(w_k^s\right)-\nabla P_{S_1}^+\left(w_{k-1}^s\right)\right)\right)}\right.\\
  	&+ \left.\frac{\left(1-\alpha^{h(\sigma_1s+\sigma_2k)}\right)\cdot\left(\left(w_k^s-w_{k-1}^s\right)^T\left(\nabla P_{S_2}^+\left(w_k^s\right)-\nabla P_{S_2}^+\left(w_{k-1}^s\right)\right)\right)}{\left\|\nabla P_{S_2}^+\left(w_k^s\right)-\nabla P_{S_2}^+\left(w_{k-1}^s\right)\right\|^2}\right),
  	\end{aligned}
  	\end{equation}
  	\begin{equation}
  	\label{mS2GD-RHBB+}
  	\begin{aligned}
  	(\tilde{\eta}_k^s)^{\mathrm{RHBB+}}
  	& = \frac{\gamma_2}{\max\{|S_1|, |S_2|\}}\cdot\left( \frac{\alpha^{h(\sigma_1s+\sigma_2k)}\cdot\left\|w_k^s-w_{k-1}^s\right\|^2}{\left(\left(w_k^s-w_{k-1}^s\right)^T\left(\nabla P_{S_1}^+\left(w_k^s\right)-\nabla P_{S_1}^+\left(w_{k-1}^s\right)\right)\right)}\right.\\
  	&+ \left.\frac{\left(1-\alpha^{h(\sigma_1s+\sigma_2k)}\right)\cdot\left(\left(w_k^s-w_{k-1}^s\right)^T\left(\nabla P_{S_2}\left(w_k^s\right)-\nabla P_{S_2}^+\left(w_{k-1}^s\right)\right)\right)}{\left\|\nabla P_{S_2}^+\left(w_k^s\right)-\nabla P_{S_2}^+\left(w_{k-1}^s\right)\right\|^2}\right).
  	\end{aligned}
  	\end{equation}
  	
  	\section{Algorithms}
  	We shall first clarify the notations: $\{v_k^s\}$ denotes the estimate sequence of $\nabla P(\cdot)$ in MB-SARAH-RHBB/RHBB+, each estimate is in a recursive form with
  	\begin{equation}
  	\label{v}
  	v_k^s=\nabla P_S\left(w_k^s\right)-\nabla P_S\left(w_{k-1}^s\right)+v_{k-1}^s.
  	\end{equation}
  	In mS2GD-RHBB/RHBB+, we signify the estimate array as $\{\tilde{v}_k^s\}$ and have 
  	\begin{equation}
  	\label{v2}
  	\tilde{v}_k^s=\nabla P_S\left(w_k^s\right)-\nabla P_S(\widetilde{w})+\nabla P(\widetilde{w}).
  	\end{equation}	
  	Notably, $\{w_k^s\}$ represents the inner iterative sequence within the $s$-th outer epoch, and we use $\{\widetilde{w}_{s}\}$ for the outer series and $\widetilde{w}$ for the snapshots.

  	\begin{algorithm}
  		\label{alg1}
  		\caption{MB-SARAH-RHBB/RHBB+}
  		\KwIn{$\widetilde{w}_0$, update frequency $m$, batch sizes $b, b_1, b_2$, constant sequence $\{\eta_0^s\}$, modification parameter $\gamma>0$, hedge base $\alpha$ and monotone decreasing function $h$, probability distribution $Q$.}
  		\vspace{2pt}
  		\KwOut{approximate solution $\widetilde{w}_s$.}
  		
  		\vspace{4pt}
  		\textbf{Outer Loop:} \For{$s = 1,2,...,$}
  		{$w_0^s=\widetilde{w}_{s-1}$\\
  			\vspace{4pt}
  			$v_0^s=\frac{1}{n} \sum_{i \in \Omega} \nabla f_i\left(w_0^s\right)=\nabla P\left(w_0^s\right)$\\
  			\vspace{4pt}
  			$w_1^s=w_0^s-\eta_0^s v_0^s$\\
  			\vspace{4pt}
  			\textbf{Inner Loop:} \For{$k = 1,2,...,m-1$}
  			{   \vspace{4pt}
  				Pick subset $S \subset\{1, \ldots, n\}$ of size $b$ uniformly at random\\
  				\vspace{4pt}
  				Update the stochastic recursive gradient $v_k^s$ by\
  				$$v_k^s=\nabla P_S\left(w_k^s\right)-\nabla P_S\left(w_{k-1}^s\right)+v_{k-1}^s$$
  				
  				Compute the step size $\eta_k^s$ by {\bf (1) Option I} or {\bf (2) Option II} \\
  				\vspace{8pt}
  				\textbf{(1) Option I : RHBB rule}\\
  				\vspace{4pt}
  				Configure $Q$ as uniform probability distribution \\
  				\vspace{4pt}
  				Pick subset $S_1 \subset\{1, \ldots, n\}$ of size $b_1$ randomly according to $Q$\\
  				\vspace{4pt}
  				Pick subset $S_2 \subset\{1, \ldots, n\}$ of size $b_2$ randomly according to $Q$\\
  				\vspace{4pt}
  				Calculate $\eta_k^s$ according to (\ref{SARAH-RHBB})\\ 
  				
  				\vspace{16pt}
  				\textbf{(2) Option II : RHBB+ rule}\\
  				\vspace{4pt}
  				Configure $Q$ to our needs \\
  				\vspace{4pt}
  				Pick subset $S_1 \subset\{1, \ldots, n\}$ of size $b_1$ randomly according to $Q$\\
  				\vspace{4pt}
  				Pick subset $S_2 \subset\{1, \ldots, n\}$ of size $b_2$ randomly according to $Q$\\
  				\vspace{4pt}
  				Compute $\nabla P_{S_1}^+, \nabla P_{S_2}^+$ according to (\ref{P+})\\
  				\vspace{4pt}
  				Calculate $\eta_k^s$ according to (\ref{SARAH-RHBB+})\\
  				\vspace{16pt}
  				
  				Update the iterate by		
  				$$w_{k+1}^s=w_k^s-\eta_k^s v_k^s$$
  				
  			}
  			$\widetilde{w}_s=w_m^s$
  		}
  	\end{algorithm}

  	\begin{algorithm}
  		\label{alg2}
  		\caption{mS2GD-RHBB/RHBB+}
  		\KwIn{$\widetilde{w}_0$, update frequency $m$, batch sizes $b, b_1, b_2$, constant sequence $\{\tilde{\eta}_0^s\}$, balance parameter $\gamma_2\geq1$, hedge base $\alpha$ and monotone decreasing function $h$, probability distribution $Q$.}
  		\vspace{2pt}
  		\KwOut{approximate solution $\widetilde{w}_s$.}
  		
  		\vspace{4pt}
  		\textbf{Outer Loop:} \For{$s = 1,2,...,$}
  		{$\widetilde{w}=\widetilde{w}_{s-1}$\\
  			\vspace{4pt}	
  			$w_0^s=\widetilde{w}$\\
  			\vspace{4pt}
  			$\varphi=\frac{1}{n} \sum_{i \in \Omega} \nabla f_i\left(\widetilde{w}\right)=\nabla P\left(\widetilde{w}\right)$\\
  			\vspace{4pt}
  			$\tilde{v}_0^s=\varphi$\\
  			\vspace{4pt}
  			$w_1^s=w_0^s-\tilde{\eta}_0^s \tilde{v}_0^s$\\
  			\vspace{4pt}
  			\textbf{Inner Loop:} \For{$k = 1,2,...,m-1$}
  			{   \vspace{4pt}
  				Pick subset $S \subset\{1, \ldots, n\}$ of size $b$ uniformly at random\\
  				\vspace{4pt}
  				Update the semi-stochastic gradient by\
  				$$\tilde{v}_k^s=\nabla P_S\left(w_k^s\right)-\nabla P_S(\widetilde{w})+\varphi$$
  				
  				Compute the step size $\tilde{\eta}_k^s$ by {\bf (1) Option I} or {\bf (2) Option II}: \\
  				\vspace{8pt}
  				\textbf{(1) Option I : RHBB rule}\\
  				\vspace{4pt}
  				Configure $Q$ as uniform probability distribution \\
  				\vspace{4pt}
  				Pick subset $S_1 \subset\{1, \ldots, n\}$ of size $b_1$ randomly according to $Q$\\
  				\vspace{4pt}
  				Pick subset $S_2 \subset\{1, \ldots, n\}$ of size $b_2$ randomly according to $Q$\\
  				\vspace{4pt}
  				Calculate $\tilde{\eta}_k^s$ according to (\ref{mS2GD-RHBB})\\
  				\vspace{16pt}
  				
  				\textbf{(2) Option II : RHBB+ rule}\\
  				\vspace{4pt}
  				Configure $Q$ to our needs \\
  				\vspace{4pt}
  				Pick subset $S_1 \subset\{1, \ldots, n\}$ of size $b_1$ randomly according to $Q$\\
  				\vspace{4pt}
  				Pick subset $S_2 \subset\{1, \ldots, n\}$ of size $b_2$ randomly according to $Q$\\
  				\vspace{4pt}
  				Compute $\nabla P_{S_1}^+, \nabla P_{S_2}^+$ according to (\ref{P+})\\
  				\vspace{4pt}
  				Calculate $\tilde{\eta}_k^s$ according to (\ref{mS2GD-RHBB+})\\
  				\vspace{16pt}
  				
  				Update the iterate:		
  				$$w_{k+1}^s=w_k^s-\tilde{\eta}_k^s \tilde{v_k}^s$$
  			}
  			$\widetilde{w}_s=w_m^s$
  		}
  	\end{algorithm}
  	
  	\newpage
  	
  	\begin{remark}
  		At beginning of each epoch $s$, the constant step sizes $\eta_0^s$ and $\tilde{\eta}_0^s$ are used in the first deterministic step of full pass. Our RHBB or RHBB+ are placed in following stochastic stages to match stochastic recursive or semi-stochastic estimators and create smooth paths for the convergence. Distribution $Q$ can be tailored to the needs of particular data sets.
  	\end{remark}
  	
  	\section{Convergence Analysis}
  	%	\begin{definition}
  	Hereafter, we use following notations:
  	let the batch correction $\overline{b}=\max\{|S_1|, |S_2|\}=\max\{b_1, b_2\}$. 
  	%	\end{definition}
  	%\begin{definition}
  	Under the probability distribution $Q \sim\left\{q_1, q_2, \ldots q_{n}\right\}$, let%we define the $L_q$, $L_r$, $\mu_q$ and $\mu_{r}$ as
  	\begin{equation}
  	\label{L_u_Q}
  	L_q=\max _i \frac{L}{n\cdot q_i}, \quad L_r=\frac{L}{L_q}, \quad \mu_q=\min _i \frac{\mu}{n\cdot q_i}, \quad \mu_r = \frac{\mu_q}{\mu}.
  	\end{equation}
  	%	\end{definition}
  	
  	Then, we have straightforward results: $L_q\geq L$, $\mu_q\leq\mu$ and $L_r\leq1$, $\mu_r\leq1$. Besides, we obtain an approximate condition number $\kappa^+ = \frac{L_q}{\mu_q}=\frac{\kappa}{L_r\mu_r}\geq\kappa$.
   
   Since the adaptor $h(\cdot)$ is monotone decreasing, we need the following uniform boundness assumption in the convergence analysis.
   
  	\noindent
  	\textbf{Assumption 3} ({\bf Uniform boundness}). The iterative adaptor $h(\cdot)$ is continuous over the bounded closed domain, monotone decreasing with respect to the epoch count $s$ and the inner count $k$, i.e., there exist constants $\hat{\alpha}$, $\tilde{\alpha}$ such that
  	$$
  	\label{assumption3}
  	1<\tilde{\alpha}<\alpha^{h(\sigma_1s+\sigma_2k)}<\hat{\alpha}.
  	$$

  	\subsection{MB-SARAH-RHBB and MB-SARAH-RHBB+}
  	To begin with, we provide subsequent Lemma \ref{lemma1} to show the summative boundary of $\mathbb{E}\left[\left\|\nabla P\left(w\right)\right\|^2\right]$ within the $s$-th epoch (inner loop).
  	\begin{lemma}
  		\label{lemma1}
  		Suppose that Assumption 1, 2a and 3 hold. The subsets $S$, $S_1$, $S_2$ are selected uniformly at random of size $b, b_1, b_2$ respectively. Then, for any $s\geq1$ in MB-SARAH-RHBB, we have
  		$$
  		\begin{aligned}
  		& \sum_{k=0}^m \mathbb{E}\left[\left\|\nabla P\left(w_k^s\right)\right\|^2\right] \leq \frac{2 \mu \overline{b} L}{\hat{\alpha} \gamma L+(1-\tilde{\alpha})\gamma\mu} \mathbb{E}\left[P\left(w_0^s\right)-P\left(w_*\right)\right] \\
  		& +\sum_{k=0}^m \mathbb{E}\left[\left\|\nabla P\left(w_k^s\right)-v_k^s\right\|^2\right]-\left(1-\frac{\hat{\alpha} \gamma L^2+(1-\tilde{\alpha}) \gamma L \mu}{\mu \overline{b} L}\right) \sum_{k=0}^m \mathbb{E}\left[\left\|v_k^s\right\|^2\right].
  		\end{aligned}
  		$$
  		Furthermore, if Assumption 2b holds, and the subsets $S_1$ and $S_2$ are sampled according to the probability distribution $Q$ of size $b_1$ and $b_2$. For any $s\geq1$ in MB-SARAH-RHBB+, we have
  		$$
  		\begin{aligned}
  		& \sum_{k=0}^m \mathbb{E}\left[\left\|\nabla P\left(w_k^s\right)\right\|^2\right] \leq \frac{2 \mu_q \overline{b} L_q}{\hat{\alpha} \gamma L_q+(1-\tilde{\alpha})\gamma\mu_q} \mathbb{E}\left[P\left(w_0^s\right)-P\left(w_*\right)\right] \\
  		& +\sum_{k=0}^m \mathbb{E}\left[\left\|\nabla P\left(w_k^s\right)-v_k^s\right\|^2\right] - \left(1-\frac{\hat{\alpha} \gamma L L_q+(1-\tilde{\alpha}) \gamma L \mu_q}{\mu_q \overline{b} L_q}\right) \sum_{k=0}^m \mathbb{E}\left[\left\|v_k^s\right\|^2\right].
  		\end{aligned}
  		$$
  	\end{lemma}
  	\proof  Deferred to the Appendix A.

       \vspace{8pt}
  	Next, we prove that the deviation (expected distance) of the full gradient to the recursive estimates is upper bounded within the $s$-th epoch (inner loop).
  	\begin{lemma}
  		\label{lemma2}
  		Suppose that Assumption 1, 2a hold. The subsets $S$, $S_1$, $S_2$ are selected uniformly at random of size $b, b_1,  b_2$, respectively. Within the $s$-th epoch of MB-SARAH-RHBB, for any $1\leq k \leq m$, we have
  		\begin{equation}
  		\label{lemma21}
  		\mathbb{E}\left[\left\|\nabla P\left(w_k^s\right)-v_k^s\right\|^2\right] \leq \frac{n-b}{b\left(n-1\right)}\left(\frac{\hat{\alpha}\gamma L^2+\left(1-\tilde{\alpha}\right)\gamma\mu L}{\overline{b} \mu L}\right)^2\sum_{j=1}^k \mathbb{E}\left[\left\|v_{j-1}^s\right\|^2\right].
  		\end{equation}
  		If Assumption 2b holds further and the subsets $S_1$, $S_2$ are sampled according to the probability distribution $Q$. Within the $s$-th epoch of MB-SARAH-RHBB+, for any $1\leq k \leq m$, we then have
  		$$
  		\mathbb{E}\left[\left\|\nabla P\left(w_k^s\right)-v_k^s\right\|^2\right] \leq \frac{L_r^2\left(n-b\right)}{b\left(n-1\right)}\left(\frac{\hat{\alpha}\gamma L_q + (1-\tilde{\alpha})\gamma  \mu_q }{\overline{b} \mu_q }\right)^2\sum_{j=1}^k \mathbb{E}\left[\left\|v_{j-1}^s\right\|^2\right].
  		$$
  	\end{lemma} 	
  	\proof  Deferred to the Appendix B.

\vspace{8pt}
  	By employing Lemma \ref{lemma1} and Lemma \ref{lemma2}, we are adequate to provide the theoretical analysis of inner loops in Theorem \ref{theorem1}.
  	\begin{theorem}
  		\label{theorem1}
  		Suppose that Assumption 1, 2a hold. Pick the subsets $S, S_1, S_2 \subset \{1, \ldots, n\}$ of size $b, b_1, b_2$ uniformly at random. Parameters $b$, $\gamma$ are chosen under a simple and suitable condition \footnote{This condition will be specified as \eqref{theorem1i} in Appendix C.}. Within the $s$-th epoch of MB-SARAH-RHBB, for any finite $m>1$, we have
  		\begin{equation}
  		\label{rate1}
  		\mathbb{E}\left[\left\|\nabla P\left(w_m^s\right)\right\|^2\right] \leq \frac{2 \overline{b} \mu L}{\gamma(m+1)(\hat{\alpha} L+(1-\tilde{\alpha})\mu)}\left[P\left(w_0^s\right)-P\left(w_*\right)\right].
  		\end{equation}
  		If Assumption 2b holds further and the subsets $S_1$, $S_2$ are sampled according to the probability distribution $Q$. Parameters $b$, $\gamma$ are chosen under another simple and suitable condition \footnote{The condition will be specified as \eqref{theorem1ii} in Appendix C.}. Within the $s$-th epoch of MB-SARAH-RHBB+, for any finite $m>1$, we thus have
  		$$
  		\mathbb{E}\left[\left\|\nabla P\left(w_m^s\right)\right\|^2\right] \leq \frac{2 \overline{b} \mu_q L_q}{\gamma(m+1)(\hat{\alpha} L_q+(1-\tilde{\alpha})\mu_q)}\left[P\left(w_0^s\right)-P\left(w_*\right)\right].
  		$$
  	\end{theorem}
  	\proof  Deferred to the Appendix C.

\vspace{8pt}
  	Theorem \ref{theorem1} shows sublinear convergence rates of the inner loops, i.e., the inner $\{\left\|\nabla P\left(w_k^s\right)\right\|^2\}$ converges sublinearly in expectation with increasing $m$. Indeed, we're sufficient to fix $s=1$ to dispose of the outer epoch, Algorithm \ref{alg1} degenerates to MB-SARAH-IN-RHBB/RHBB+ (see \cite{nguyen2} for reference). To obtain an $\varepsilon$-accurate solution in MB-SARAH-IN-RHBB, the number of iterations, $m$, is put up so that $\mathbb{E}\left[\|\nabla P(w_m)\|^2\right] \leq \varepsilon$, which suggests that 
  	\begin{equation}
  	\label{con1}
  	\frac{2 \overline{b} \mu L}{\gamma\left(m+1\right)\left(\hat{\alpha} L+\left(1-\tilde{\alpha}\right)\mu\right)}\left[P\left(w_0\right)-P\left(w_*\right)\right] \leq \varepsilon.
  	\end{equation}
  	Assume that $P\left(w_0\right)-P\left(w_*\right)=\sigma$, (\ref{con1}) implies $m_{RH} = \lceil \frac{2\overline{b}\mu \sigma\kappa}{\varepsilon \gamma\left(\hat{\alpha} \kappa+1-\tilde{\alpha}\right)}-1\rceil$. Compare with $m_{R} = \lceil \frac{2\overline{b}\mu \sigma}{\varepsilon \gamma}-1\rceil$ in \cite{yang}, we have $\left(m_{RH}<m_{R}\right)$ due to $L > \mu$. In MB-SARAH-RHBB+, that's $m_{RH+} = \lceil \frac{2\overline{b}\mu_q \sigma\kappa^+}{\epsilon \gamma(\hat{\alpha} \kappa^+ +1-\tilde{\alpha})}-1\rceil$ to achieve the same $\varepsilon$-accuracy. Regardless of rounding errors, it's very likely that $m_{RH+}\leq m_{RH}$.
  	
  	Ineq. (\ref{con1}) indicates that we can locally manipulate $\hat{\alpha}$, $\tilde{\alpha}$ in the early epochs (e.g., by temporarily using a different $h$) to address the issue of a poor initial $w_0^0$ (or $w_0^s$) with an unexpectedly large $\sigma$. To our best know, this issue can not be effectively resolved in many existing methods, e.g., \cite{yang} \cite{yang3} \cite{nguyen} \cite{nguyen2} \cite{yangy} \cite{liu3}.  
  	
  	For a class of ill-conditioned objective functions $P(\cdot)$ under $L \gg \mu$, (\ref{con1}) implies $m_{RH+}+1\leq m_{RH}+1\approx \frac{1}{\hat{\alpha}}\left(m_{R}+1\right)$, suggesting the inner speedup is approximately proportional to $\mathcal{O}\left(\frac{1}{\hat{\alpha}}\right)$. In ill-conditioning, we tolerate towards the functional form and the decay rate of $h$, but remain focus on $\hat{\alpha}$.
  	
  	Next, we evaluate the workload in terms of incremental first order oracle (IFO) complexity model in Agarwal et al.\cite{agarwal}. In \cite{ghadimi}, it's SFO under stochastic settings. MB-SARAH-RHBB/RHBB+ are IFO algorithms that are specified through calls to an IFO, regardless of $P(\cdot)$. Each epoch invokes SFO at most $2bm$ times in the recursive gradient evaluations (\ref{v}), corresponding to an overall cost of $\mathcal{O}\left(n+2bm\right)$ SFOs. Since the adaptor $h$ is mentor-specified, we can force $\hat{\alpha}$ and $\tilde{\alpha}$ to be arbitrarily large and small. By setting $L+\frac{1-\tilde{\alpha}}{\hat{\alpha}}\mu=\mathcal{O}\left(L\right)$ and $L_q+\frac{1-\tilde{\alpha}}{\hat{\alpha}}\mu_q=\mathcal{O}\left(L_q\right)$, it's sufficient to have $m =\mathcal{O}\left(\frac{\overline{b}\mu}{\varepsilon \gamma\hat{\alpha}}\right)$ and $m =\mathcal{O}\left(\frac{\overline{b}\mu_q}{\varepsilon \gamma\hat{\alpha}}\right)$. Therefore, we obtain the following conclusions for the complexity bounds.
  	
  	\begin{corollary}
  		Suppose Assumption 1 and 2a hold. MB-SARAH-IN-RHBB converges sublinearly in expectation with a rate of $\mathcal{O}\left(\mu\overline{b} / \gamma m \hat{\alpha}\right)$, and the complexity to achieve an $\varepsilon$-accurate solution is in the order of $n+2bm=\mathcal{O}\left(n+\frac{b\overline{b}\mu}{\varepsilon \gamma\hat{\alpha}}\right)$. Suppose Assumption 2b holds further. MB-SARAH-IN-RHBB+ owns sublinear convergence rate of  $\mathcal{O}\left(\overline{b}\mu_q/\gamma m \hat{\alpha}\right)$, and the complexity for the  $\varepsilon$-accuracy corresponds to $n+2bm=\mathcal{O}\left(n+\frac{b\overline{b}\mu_q}{\varepsilon \gamma\hat{\alpha}}\right)$ units of work.
  	\end{corollary}
  	
  	On the basis of Theorem \ref{theorem1}, we can proceed to the analysis of multiple outer steps, and we establish the convergence of MB-SARAH-RHBB and MB-SARAH-RHBB+ in subsequent Theorem \ref{theorem2}.
  	\begin{theorem}
  		\label{theorem2}
  		Suppose that Assumption 1, 2a hold. Pick the subsets $S, S_1, S_2 \subset\{1, \ldots, n\}$ of size $b, b_1, b_2$ uniformly at random, and we choose $b, \gamma$ that satisfy condition (\ref{theorem1i}). In MB-SARAH-RHBB, for any $s>1$, we have
  		\begin{equation}
  		\label{rate2}
  		\mathbb{E}\left[\left\|\nabla P\left(\widetilde{w}_s\right)\right\|^2\right] \leq\left(\frac{\kappa\overline{b}}{\gamma\left(m+1\right)\left(\hat{\alpha}\kappa+1-\tilde{\alpha}\right)}\right)^s\left\|\nabla P\left(\widetilde{w}_0\right)\right\|^2.
  		\end{equation}
  		If Assumption 2b holds further and the subsets $S_1$, $S_2$ are sampled according to the probability distribution $Q$, and we choose $b, \gamma$ that satisfy condition (\ref{theorem1ii}). In MB-SARAH-RHBB+, for any $s>1$, we have
  		$$
  		\mathbb{E}\left[\left\|\nabla P\left(\widetilde{w}_s\right)\right\|^2\right] \leq\left(\frac{\mu_r\kappa^+\overline{b}}{\gamma\left(m+1\right)\left(\hat{\alpha}\kappa^++1-\tilde{\alpha}\right)}\right)^s\left\|\nabla P\left(\widetilde{w}_0\right)\right\|^2
  		.$$
  	\end{theorem}
  	\proof  Deferred to the Appendix D.

\vspace{8pt}
  	Theorem \ref{theorem2} indicates that the outer $\{\left\|\nabla P\left(\widetilde{w}_s\right)\right\|^2\}$ converges linearly in expectation. Assume that $\left\|\nabla P\left(\widetilde{w}_0\right)\right\|^2=\zeta$, to obtain $\mathbb{E}\left[\left\|\nabla P\left(\widetilde{w}_s\right)\right\|^2\right] <\varepsilon$ in MB-SARAH-RHBB, the number of outer epoch $s$ must satisfy
  	$$
  	\left(\frac{\kappa\overline{b}}{\gamma\left(m+1\right)\left(\hat{\alpha}\kappa+1-\tilde{\alpha}\right)}\right)^s\cdot\zeta\leq \varepsilon.
  	$$
  	It infers $s_{RH}=\lceil \frac{log(\zeta)-log(\varepsilon)}{log(\hat{\alpha}\kappa+1-\tilde{\alpha})-log(\kappa)+log(\gamma(m+1))-log(\overline{b})}\rceil$. Compared to \cite{yang} with $s_{R}=\lceil\frac{log(\zeta)-log(\epsilon)}{log(\gamma(m+1))-log(\overline{b})}\rceil$, the overhead of the outer epoch decreases. In MB-SARAH-RHBB+, both $\mu_r\leq1$, $\kappa^+\geq\kappa$ implies $s_{RH+}\leq s_{RH}$, indicating the iterative cost can be further reduced via the effective sampling. In ill-conditioning, we then have $s_{RH^+}\leq s_{RH}\approx\lceil\frac{log(\zeta)-log(\varepsilon)}{log(\hat{\alpha})+log(\gamma(m+1))-log(\overline{b})}\rceil$.
  	
    	Furthermore, our analysis can be refined to obtain smaller rate constants in some gradient dominated scenarios (see in Polyak et al.\cite{polyak}, Reddi et al. \cite{reddi}). If $P(\cdot)$ is $\delta$-gradient dominated with $\delta<\frac{1}{2\mu}$, we derive the rate constants
  	$$
  	\rho'_{RH}=\frac{2\overline{b}\mu L\delta}{\gamma(m+1)\left(\hat{\alpha}L+(1-\tilde{\alpha})\mu\right)}, \quad
  	\rho'_{RH+}=\frac{2 \overline{b} \mu_q L_q\delta}{\gamma(m+1)(\hat{\alpha} L_q+(1-\tilde{\alpha})\mu_q)},
  	$$
  	for MB-SARAH-RHBB and MB-SARAH-RHBB+, respectively. The theoretical convergence speed further increases in virtue of $2\mu\delta<1$.
  	
  	\begin{corollary}
  		\label{corollary2}
  		Suppose that Assumption 1 and 2a hold. MB-SARAH-RHBB converges linearly with the total complexity to achieve an $\varepsilon$-accurate solution as $\mathcal{O}\left(\left(n+\frac{b\overline{b}\mu}{\varepsilon \gamma\hat{\alpha}}\right)\log(1/\varepsilon)\right)$. Suppose that Assumption 2b holds further. MB-SARAH-RHBB+ obtain linear convergence rate, and the overall complexity for the same $\varepsilon$-accuracy is of order $\mathcal{O}\left(\left(n+\frac{b\overline{b}\mu_q}{\varepsilon \gamma\hat{\alpha}}\right)\log\left(1/\varepsilon\right)\right)$.
  	\end{corollary}
  	
  	Compared with MB-SARAH \cite{nguyen2}, MB-SARAH-RBB \cite{yang}, MB-SARAH-HD \cite{yang5}, iSARAH-BB \cite{yangy}, Corollary \ref{corollary2} indicates that MB-SARAH-RHBB/RHBB+ have lower complexity when using an appropriate adaptor $h(\cdot)$ and a proper $\overline{b}$.

  	\subsection{mS2GD-RHBB and mS2GD-RHBB+}
  	We exhibit the following lemma, based on Lemma $2$ from \cite{yang2}, to start the convergence analysis for mS2GD-RHBB/RHBB+.  
  	\begin{lemma}
  		\label{lemma3}
  		Suppose that Assumption 1, 2a hold. The subset $S$ is selected uniformly at random with size $b$. Then, we have an upper bound for the semi-stochastic estimate $\tilde{v}$ (\ref{v2}) as follows
  		\begin{equation}
  		\mathbb{E}\left[\left\|\tilde{v}_k^s\right\|^2\right] \leq \frac{4 L}{b}\left[P\left(w_{k-1}^s\right)\!-P\left(w_*\right)\!+P(\widetilde{w}_{s-1})\!-P\left(w_*\right)\right] \\
  		+\frac{2}{b}\left\|\nabla P\left(w_{k-1}^s\right)\right\|^2 .
  		\end{equation}
  	\end{lemma}
  	\proof  Deferred to the Appendix E.
   
       \vspace{8pt}
  	Based on Lemma \ref{lemma3}, we present subsequent Theorem 3 to demonstrate the linear convergence of mS2GD-RHBB and mS2GD-RHBB+.
  	\begin{theorem}
  		\label{theorem3}
  		Suppose that Assumption 1, 2a hold. Let $\kappa_r=\hat{\alpha}\kappa+1-\tilde{\alpha}$, and pick the subsets $S, S_1, S_2 \subset\{1, \ldots, n\}$ of size $b, b_1, b_2$ uniformly at random. Assume that $b\overline{b}>4\kappa_r\gamma_2$, and $h(\cdot)$ is chosen such that
  		\begin{equation}
  		\label{rate3}
  		\tilde{\rho}_1=\frac{\kappa b \overline{b}^2}{m\gamma_2\kappa_r\left(b \overline{b}-4\gamma_2\kappa_r\right)}+\frac{2\gamma_2 \kappa_r}{b \overline{b}-4\gamma_2\kappa_r}<1.
  		\end{equation}
  		Then, mS2GD-RHBB converges linearly in expectation with rate $\tilde{\rho}_1$, that's
  		$$
  		\mathbb{E}\left[P\left(\widetilde{w}_s\right)\right]-P\left(w_*\right) \leq (\tilde{\rho}_{1})^s\left[P\left(\widetilde{w}_0\right)-P\left(w_*\right)\right].
  		$$
  		If Assumption 2b holds further, let $\kappa_r^+=\hat{\alpha} \kappa^+ +1-\tilde{\alpha}$, and sample the subsets $S_1, S_2 \subset\{1, \ldots, n\}$ according to the probability distribution $Q$. Assume that $b\overline{b}>4\kappa_r^+\gamma_2L_r$, and $h(\cdot)$ is chosen such that
  		\begin{equation}
  		\label{rate4}
  		\tilde{\rho}_{2}=\frac{\mu_r\kappa^+b\overline{b}^2}{m\gamma_2\kappa_r^+\left(b \overline{b}-4\gamma_2\kappa_r^+ L_r\right)}+\frac{2 \gamma_2\kappa_r^+L_r}{b \overline{b}-4\gamma_2\kappa_r^+L_r}<1.
  		\end{equation}
  		Then, mS2GD-RHBB+ converges linearly in expectation with rate $\tilde{\rho}_2$, that's 
  		$$
  		\mathbb{E}\left[P\left(\widetilde{w}_s\right)\right]-P\left(w_*\right) \leq (\tilde{\rho}_{2})^s\left[P\left(\widetilde{w}_0\right)-P\left(w_*\right)\right].
  		$$
  	\end{theorem}
  	\proof  Deferred to the Appendix F.

        \vspace{8pt}
  	\paragraph{Further Discussion on $h(\cdot)$: }
	Here, we show how to find an $h(\cdot)$ to ensure the improvement in terms of theory. For clarity, we use $\tilde{\rho}_{R}$, $\tilde{\rho}_{RH}$, $\tilde{\rho}_{RH+}$ to denote the convergence rates of mS2GD-RBB, mS2GD-RHBB, mS2GD-RHBB+, respectively. In mS2GD-RBB, the update frequency $m$ and the batch sizes $b, b_1$ are chosen (here, $\overline{b}=b_1$) such that
  	\begin{equation}
  	\label{RBB}
  	\tilde{\rho}_{R}=\frac{\mu b\overline{b}^2+2mL}{\mu b \overline{b} m-4mL}<1.
  	\end{equation}
  	With the identical parameter set $\{m, b, \overline{b}\}$, mS2GD-RHBB possesses $$\tilde{\rho}_{RH}=\frac{\left(\frac{\kappa}{\kappa_r\gamma_2}\right)^2\cdot\mu b\overline{b}^2+2mL}{\left(\frac{\kappa}{\kappa_r\gamma_2}\right)\cdot\mu b \overline{b} m-4mL}.$$
  	Assume that $\tilde{\rho}_{R}=\frac{\mu b\overline{b}^2+2mL}{\mu b \overline{b} m-4mL}=c<1$, we then obtain
  	\begin{equation}
  	\label{RHBB good}
  	mL=\frac{c}{2+4c}m\mu b\overline{b} -  \frac{1}{2+4c}\mu b \overline{b}^2.
  	\end{equation}
  	According to (\ref{RHBB good}), the rate of $\tilde{\rho}_{RH}$ can be re-expressed into
  	$$
  	\tilde{\rho}_{RH}=\frac{1}{2}\cdot\frac{\left((1+2c)(\frac{\kappa}{\kappa_r\gamma_2})^2-1-c\right)\mu b\overline{b}^2+c\mu b\overline{b}^2 + c\mu b\overline{b}m}{\left(\frac{\left(1+2c\right)\kappa}{2\kappa_r\gamma_2}-1-c\right)\mu b\overline{b}m+\mu b\overline{b}^2 + \mu b\overline{b}m}.
  	$$
  	Let's mark the term $\left((1+2c)(\frac{\kappa}{\kappa_r\gamma_2})^2-1-c\right)$ by $A$ and the term $\left(\frac{\left(1+2c\right)\kappa}{2\kappa_r\gamma_2}-1-c\right)$ by $B$. If $A\overline{b}<c Bm$ satisfied, it follows $\tilde{\rho}_{RH}<\frac{1}{2}\tilde{\rho}_{R}=\frac{1}{2}c$. Note that $A<0$ if $\frac{\kappa_r}{\kappa}>\frac{1}{\gamma_2}\sqrt{\frac{1+2c}{c+1}}$ and $B<0$ if $\frac{\kappa_r}{\kappa}>\frac{2c+1}{\gamma_2\left(2c+2\right)}$. Therefore, the balance parameter $\gamma_2$ ($\gamma_2\geq1$) relaxes the essential boundaries of terms A and B. To meet $A\overline{b}<c Bm$, we should choose a $\gamma_2$ such that 
   \begin{equation}
   \label{condition1}
  	\overline{b}>\frac{\frac{\left(1+2c\right)c\kappa}{2\kappa_r\gamma_2}-c-c^2}{(1+2c)\cdot(\frac{\kappa}{\kappa_r\gamma_2})^2-1-c}\cdot m.
   \end{equation}
  	It suggests that configuring an $h(\cdot)$ that satisfies (\ref{condition1}) will realize a significant speed of $\tilde{\rho}_{RH}<\frac{1}{2}\tilde{\rho}_{R}$ in theory. Our rules allows to flexibly trade-off between $h(\cdot)$ and $\gamma_2$, however, it's not necessary to strictly tune out a rate constant that prompts more than twice improvement (shrinking to less than the half). In practice, we commonly set $\gamma_2$ and configure $\hat{\alpha}$, $\tilde{\alpha}$ slightly larger than $1$.
  	
  	For comprehensive analysis, we explicate it from the aspect of the effective range. According to (\ref{RBB}), we derive the effective range of $\tilde{\rho}_{RH}$ as 
  	\begin{equation}
  	\label{RHBB good2}
  	\tilde{\rho}_{RH} < \frac{\left(\left(\frac{\kappa}{\gamma_2\kappa_r}\right)^2-\frac{1}{3}\right)\mu b \overline{b}^2+\frac{1}{3}m\mu b\overline{b}} {\frac{2}{3}\mu b\overline{b}^2+\left(\frac{\kappa}{\gamma_2\kappa_r}-\frac{2}{3}\right)m\mu b\overline{b}}< 1+ \frac{\left(\left(\frac{\kappa}{\gamma_2\kappa_r}\right)^2-1\right)\overline{b}-\left(\frac{\kappa}{\gamma_2\kappa_r}-1\right)m}{\frac{2}{3}\left(\overline{b}-m\right)+\frac{\kappa}{\gamma_2\kappa_r}m}.
  	\end{equation}
  	Hence, restricting the batch correction $\overline{b}$ to a broad interval of
  	$$
  	\left(-\infty, \left(1-\frac{3\kappa}{2\gamma_2\kappa_r}\right)m\right) \cup \left(\left(1-\frac{\kappa}{\kappa+\gamma_2\kappa_r}\right)m, +\infty\right)
  	$$
  	will enforce the second term in (\ref{RHBB good2}) negative, which urges $\tilde{\rho}_{RH}$ to fall into a narrower interval (compared to the original $\tilde{\rho}_{R}<1$). By selecting the set $\{m, b, b_1\}$ identical to mS2GD-RBB, tuning $b_2$ can affect the lower bound of convergence speed, suggesting that the worst convergence result is also a fast one. A straightforward trade-off of $\overline{b}>m$ will enforce the second term in (\ref{RHBB good2}) negative, habitually applied in practice. 
  	
  	Furthermore, let's try to minimize the `ineffective' range
  	$$
  	\min _{h} \quad \left|\frac{3\kappa}{2\gamma_2\kappa_r}-\frac{\kappa}{\kappa+\gamma_2\kappa_r}\right|m.
  	$$
  	Due to the monotonicity, the `ineffective' interval shortens as $\hat{\alpha}$ increases or $\tilde{\alpha}$ decreases. It means that one can enlarge the value span of the adaptor $h(\cdot)$ or increase the decay speed (while satisfying $b\overline{b}>4\kappa_r\gamma_2$) to obtain a sufficient speedup in convergence.
  	
  	In mS2GD-RHBB+, the analysis follows a similar line of reasoning. The rate $\tilde{\rho}_{RH+}$ can be built smaller even than the $\tilde{\rho}_{RH}$, due to the facts $\mu_r\leq1$ and $\kappa^+\geq\kappa$. To achieve $\tilde{\rho}_{RH+}<\frac{1}{2}\tilde{\rho}_{R}=\frac{1}{2}c$, one should ensure $\frac{\kappa_r^+}{\kappa^+}>\frac{\mu_{r}}{\gamma_2}\sqrt{\frac{1+2c}{c+1}}$, hence, RHBB+ allows a freer selection of the exponential adaptor $h(\cdot)$.
  	
  	For a class of ill-conditioned functions $P(\cdot)$ with $L \gg \mu$, we have $\frac{\kappa}{\kappa_r} \approx \frac{\kappa^+}{\kappa_r^+} \approx \frac{1}{\hat{\alpha}}$. The associated limitations become dependent only on the upper bound $\hat{\alpha}$ of $h(\cdot)$, which is almost equivalent to the original $\hat{\alpha}>1$. This can be solved with ease at the initial inputs, saving plenty of tuning effort.

  	\vspace{8pt}
  	Note from Theorem \ref{theorem3}, it's feasible to discard the outer epoch and set up mS2GD-IN-RHBB/RHBB+ algorithms (similar to \cite{nguyen2}). Theorem \ref{theorem3} suggests that whenever the set $\{m, b_1, b_2\}$ are chosen, the second terms of (\ref{rate3}) (\ref{rate4}) can be regulated sufficiently small through $h(\cdot)$.
  	
  	By setting $\frac{\hat{\alpha}}{\tilde{\alpha}}L+\frac{1-\tilde{\alpha}}{\tilde{\alpha}}\mu=\mathcal{O}\left(\mu\right)$ and $\frac{\hat{\alpha}}{\tilde{\alpha}}L_q+\frac{1-\tilde{\alpha}}{\tilde{\alpha}}\mu_q=\mathcal{O}\left(\mu_q\right)$, from (\ref{rate3}) (\ref{rate4}) we obtain $m=\mathcal{O}\left(\frac{\overline{b}\kappa}{\gamma_2\tilde{\alpha}}\right)$ and $m=\mathcal{O}\left(\frac{\mu_{r}\overline{b}\kappa^+}{\gamma_2\tilde{\alpha}}\right)=\mathcal{O}\left(\frac{\overline{b}\kappa}{L_{r}\gamma_2\tilde{\alpha}}\right)$, correspondingly. 
  	
  	To satisfy $\mathbb{E}\left[P\left(\widetilde{w}_s\right)\right]-P\left(w_*\right) \leq \left(\tilde{\rho}_{RH}\right)^s\cdot\left[P\left(\widetilde{w}_0\right)-P\left(w_*\right)\right] \leq \varepsilon$ in mS2GD-RHBB, the number of outer epoch $s$ must satisfy
  	$$s\geq \frac{log\left(P(\widetilde{w}_0)-P(w_*)\right)-log\left(\varepsilon\right)}{-log\left(\tilde{\rho}_{RH}\right)}.$$
  	By the same token, we demand the $s$ in mS2GD-RHBB+ such that
  	$$s\geq \frac{log\left(P(\widetilde{w}_0)-P(w_*)\right)-log\left(\varepsilon\right)}{-log\left(\tilde{\rho}_{RH+}\right)}.$$
  	Therefore, to bound the number of oracles of IFO model, the following result for the total complexity is obtained.
  	\begin{corollary}
  		\label{corollary3}
  		Suppose that Assumption 1 and 2a hold. The complexity of mS2GD-RHBB to achieve an $\varepsilon$-accurate solution is $\mathcal{O}\left(\left(n+\frac{b\overline{b}}{\gamma_2\tilde{\alpha}}\kappa\right)log\left(\frac{1}{\varepsilon}\right)\right)$. Supposethat  Assumption 2b holds further. To obtain an $\varepsilon$-accurate solution, the overall complexity of mS2GD-RHBB+ is of order
  		$\mathcal{O}\left(\left(n+\frac{b\overline{b}}{L_{r}\gamma_2\tilde{\alpha}}\kappa\right)log\left(\frac{1}{\varepsilon}\right)\right)$.
  	\end{corollary}
  	
  	Compared with mS2GD \cite{konevcny}, mS2GD-BB \cite{yang3}, mS2GD-RBB \cite{yang2}, Corollary \ref{corollary3} indicates that, to achieve an $\varepsilon$-accurate solution, mS2GD-RHBB/RHBB+ have lower total complexity when choosing $h(\cdot)$ and $\overline{b}$ properly.

  	\section{Experiments}
  	\subsection{Experimental Settings}
  	To be specific, our experiments are performed on the well-worn problems of training $\ell_2$ regularized ridge regression, i.e.,
  	$$
  	\min _{w \in \mathbb{R}^d} P(w) = \frac{1}{n} \sum_{i=1}^n \log \left(1+\exp \left(-z_i x_i^T w\right)\right) +\frac{\lambda}{2}\|w\|^2.
  	$$
  	
  	For RHBB+, we specify the probability distribution $Q$ with two options. Since $f_i(w)=\log \left(1+\exp \left(-z_i x_i^T w\right)\right)$ with $z_i\in\{-1, 1\}$, then $\|\nabla f_i(w)\|\leq\|x_i\|\leq\sqrt{d}\|x_i\|_{\infty}$, hence for option I we set $q_i=\frac{\left\|x_i\right\|_{\infty}^{\tau}}{\sum_{j=1}^{n}\left\|x_j\right\|_{\infty}^{\tau}}$. For option II, we consider sparsity and set: $q_i=\frac{\left\|x_i\right\|_{0}^{\tau}}{\sum_{j=1}^{n}\left\|x_j\right\|_{0}^{\tau}}$. Here, coefficient $\tau$ is equipped to mitigate the batch influence on the importance sampling.
  	
  	We verify MB-SARAH-RHBB and mS2GD-RHBB on data sets $a8a$, $w8a$, $ijcnn1$, $covtype$, $phishing$ and $mushrooms$. Due to the statistical characteristics of distribution $Q$, we explore MB-SARAH-RHBB+ and mS2GD-RHBB+ on another three $australian$, $madelon$ and $german.numer$. All data sets are publicly available in LIBSVM \footnote{ \url{https://www.csie.ntu.edu.tw/~cjlin/libsvmtools/datasets/}.}. More details are referred to \textbf{Table \ref{table1}}.
  	
  	\begin{table}
  		\setlength{\abovecaptionskip}{0pt}
  		\setlength{\belowcaptionskip}{1pt}
  		\caption{DATA INFORMATION OF EXPERIMENTS}
  		\centering
  		\vspace{4pt}
  		\begin{tabular}{lccc}
  			\hline
  			\multicolumn{1}{l}{Datasets}&\multicolumn{1}{l}{\qquad Instances ($n$)} &\multicolumn{1}{l}{\qquad Features ($d$)} &\multicolumn{1}{c}{\qquad $\lambda$}
  			\\ \hline\noalign{\smallskip}
  			a8a & \qquad 22,696 & \qquad 123 & \qquad$10^{-2}$\\
  			w8a & \qquad 49,749 & \qquad 300 & \qquad$10^{-2}$ \\
  			ijcnn1 & \qquad 49,990 & \qquad 22 & \qquad$10^{-2}$ \\
  			covtype & \qquad 581,012 & \qquad 54 & \qquad$10^{-2}$ \\
  			phishing & \qquad 11,055 & \qquad 68 & \qquad$10^{-2}$ \\
  			mushrooms & \qquad 8,124 & \qquad 112 & \qquad$10^{-2}$ \\
  			australian & \qquad 690 & \qquad 14 & \qquad$10^{-2}$\\
  			madelon & \qquad 2,000 & \qquad 600 & \qquad$10^{-2}$\\
  			german.numer & \qquad 1,000 & \qquad 24 & \qquad$10^{-2}$\\
  			\hline
  		\end{tabular}
  		\label{table1}
  	\end{table}
  	
  	\subsection{Experiments investigating for Non-adaptive Hedge Effect}      
  	Our first aim is to investigate whether the hedge ideology is helpful to improve the numerical efficiency. By fixing $h(\sigma_1s+\sigma_2k)=1$, we separate adaptive technique from the hedge operation and verify the unvaried hedge effect. 
  	
  	We conduct experiments under $b_1=b_2=b_H$ first, where $b_H$ conveys the unified batch size and is used in the legends. For clarity, notations of this subsection are summarized in \textbf{Table 2}.
  	\begin{table}
  		
  		\setlength{\abovecaptionskip}{0pt}
  		\setlength{\belowcaptionskip}{1pt}
  		\caption{NOTATIONS DESCRIPTIONS}
  		\centering
  		\vspace{4pt}
  		\begin{tabular}{lclc}
  			\hline
  			\multicolumn{1}{l}{Notations}&\multicolumn{1}{l}{\qquad Hedge Bases} &\multicolumn{1}{l}{\qquad Step Sizes} &\multicolumn{1}{c}{\qquad Adaptivity}
  			\\ \hline\noalign{\smallskip}
  			MB-SARAH-RBB  & \qquad \ding{56} & \qquad RBB & \qquad \ding{56} \\
  			MB-SARAH-RBB+ & \qquad \ding{56} & \qquad RBB+ & \qquad \ding{56} \\
  			MB-SARAH-RHBB($\alpha$) & \qquad $\alpha$ & \qquad RHBB & \qquad \ding{56} \\
  			MB-SARAH-RHBB($\alpha$)+ & \qquad $\alpha$ & \qquad RHBB+ & \qquad \ding{56} \\
  			mS2GD-RBB  & \qquad \ding{56} & \qquad RBB & \qquad \ding{56} \\
  			mS2GD-RBB+ & \qquad \ding{56} & \qquad RBB+ & \qquad \ding{56} \\
  			mS2GD-RHBB($\alpha$) & \qquad $\alpha$ & \qquad RHBB & \qquad \ding{56} \\
  			mS2GD-RHBB($\alpha$)+ & \qquad $\alpha$ & \qquad RHBB+ & \qquad \ding{56} \\
  			\hline
  		\end{tabular}
  		\label{table2}
  	\end{table}

  	\subsubsection{Non-adaptive MB-SARAH-RHBB/mS2GD-RHBB}
  	
  	\ 
  	
  	\vskip 10pt
  	
  	\noindent\textbf{Parametric Settings: } We set $b=4$ and sample the subsets $S$, $S_1$, $S_2$ according to uniform distribution. Under $b=4$, we follow the guidelines in \cite{yang} and set $\gamma=1$ in MB-SARAH-RHBB. In mS2GD-RHBB, we conduct a conservative trade-off with a moderate $\gamma_2=1$. We fixed $b_H=40$ in general experiments and we varied $b_H=20, 30, 40, 50, 60$ in the last. We sequentially select $\alpha$ from the set $\{2, 3, 4, 5\}$ and the set $\{10, 11, 12, 13\}$. 
  	
\begin{figure*}[htbp]
  		\centering
  		\subfigure[a8a]
  		{\includegraphics[width=0.327\textwidth]{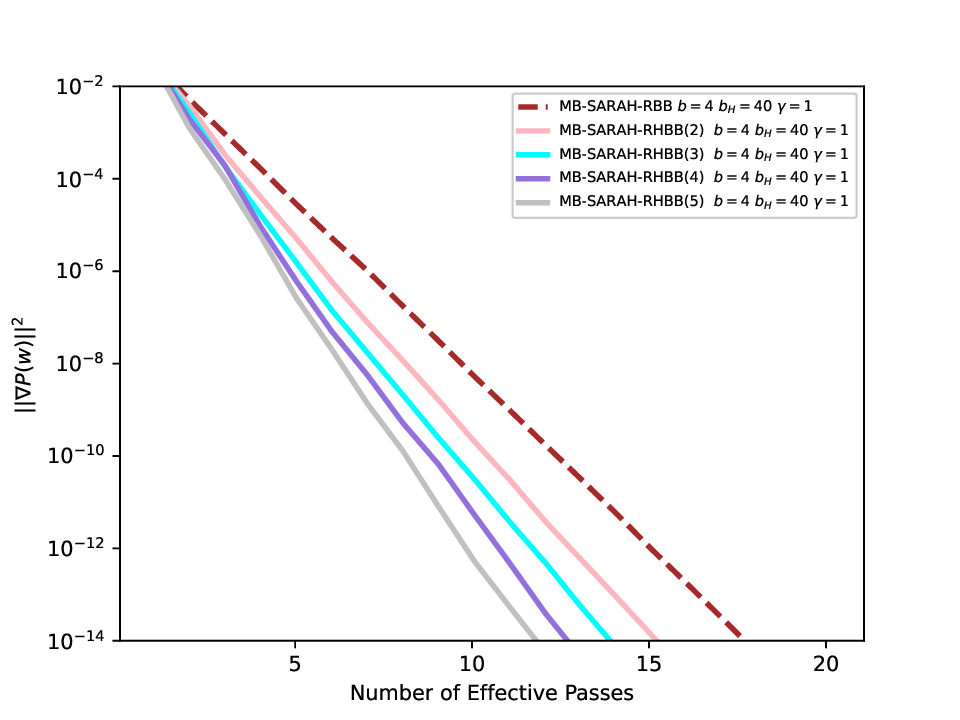}}
  		\subfigure[w8a]
  		{\includegraphics[width=0.327\textwidth]{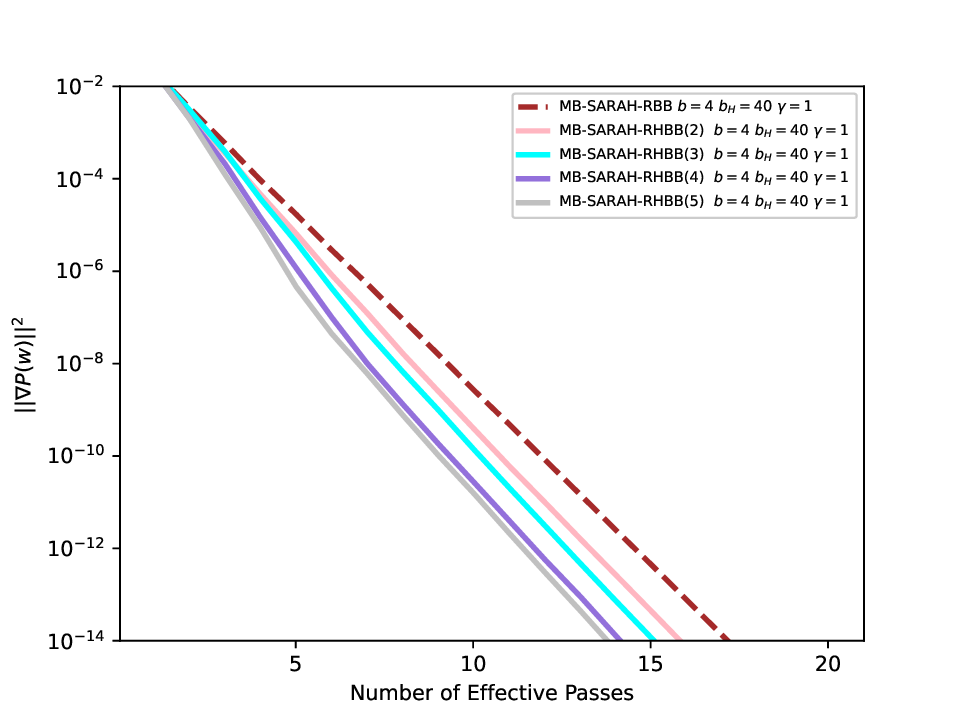}}
  		\subfigure[ijcnn1]
  		{\includegraphics[width=0.327\textwidth]{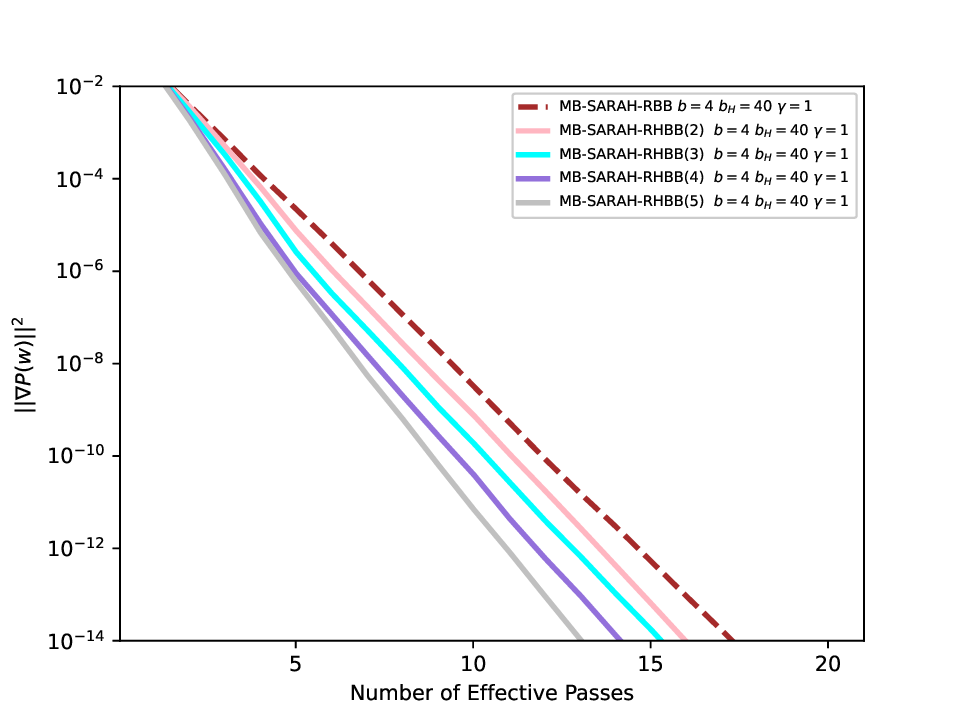}}
  		\subfigure[covtype]
  		{\includegraphics[width=0.327\textwidth]{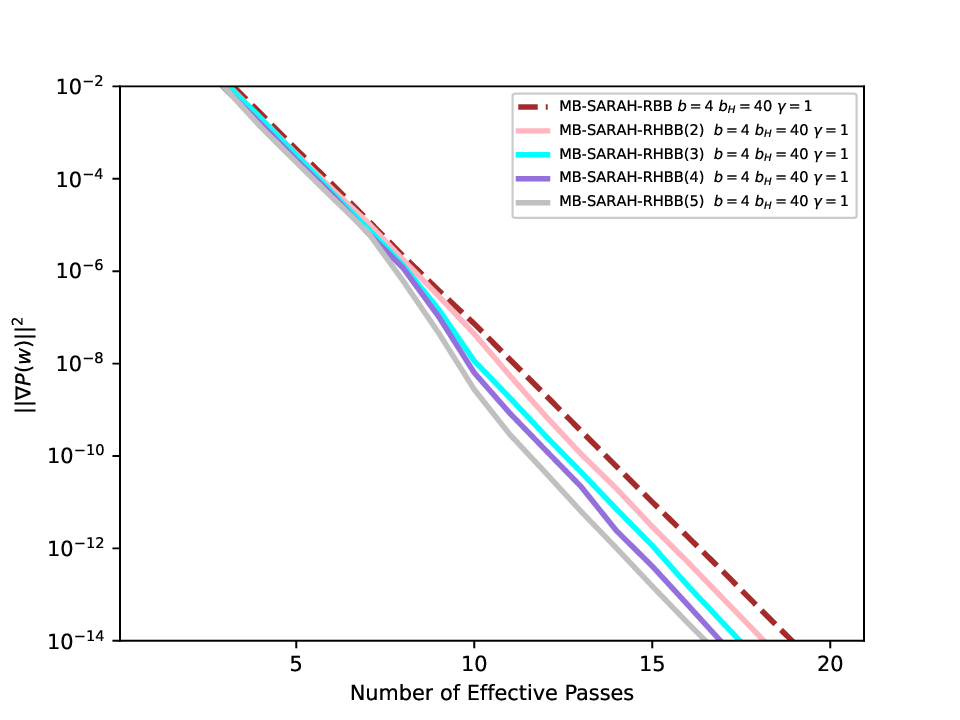}}
  		\subfigure[phishing]
  		{\includegraphics[width=0.327\textwidth]{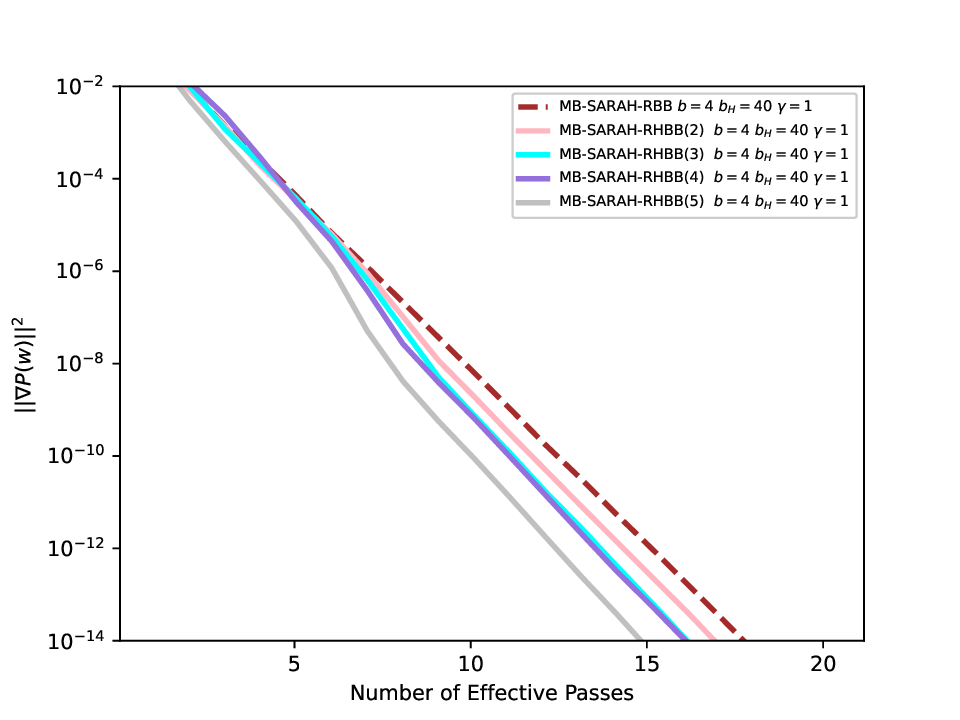}}
  		\subfigure[mushrooms]
  		{\includegraphics[width=0.327\textwidth]{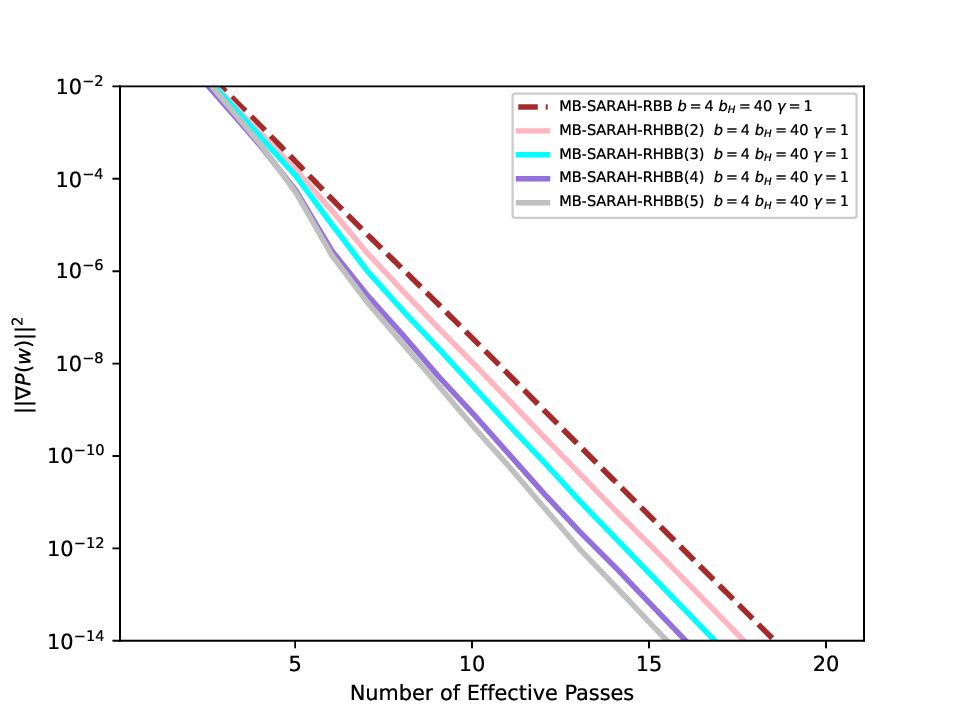}}
  		\caption{\footnotesize Comparisons of MB-SARAH-RHBB (dash line) and MB-SARAH-RBB (solid lines) with an $\alpha$ from $\{2,3,4,5\}$.}
  		\label{fig2}
  	\end{figure*}

\begin{figure*}[htbp]
  		\centering
  		\subfigure[a8a]
  		{\includegraphics[width=0.327\textwidth]{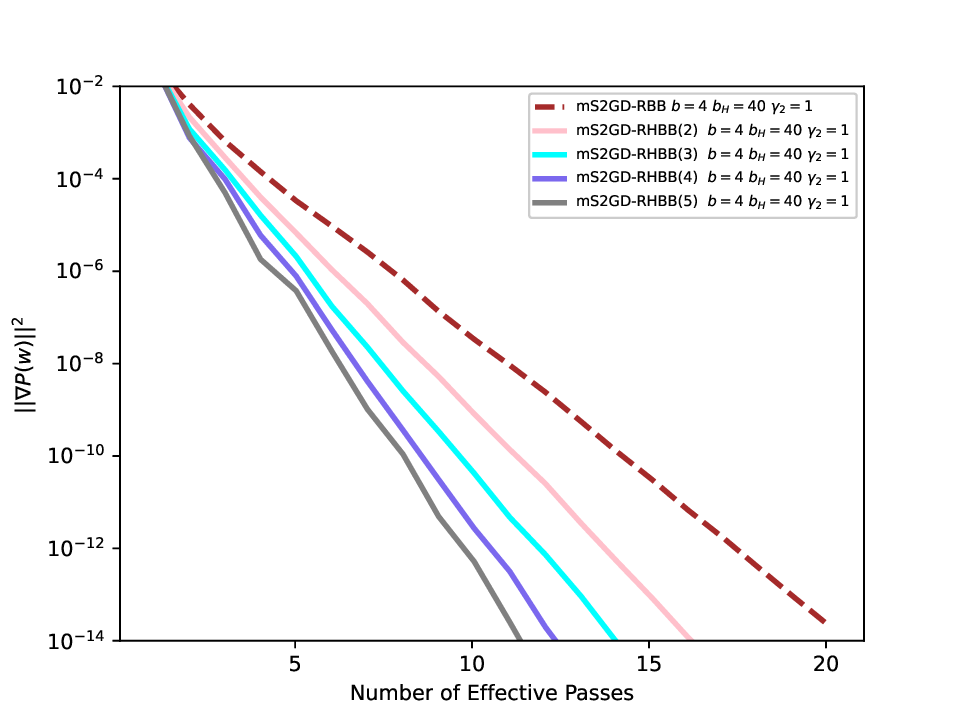}}
  		\subfigure[w8a]
  		{\includegraphics[width=0.327\textwidth]{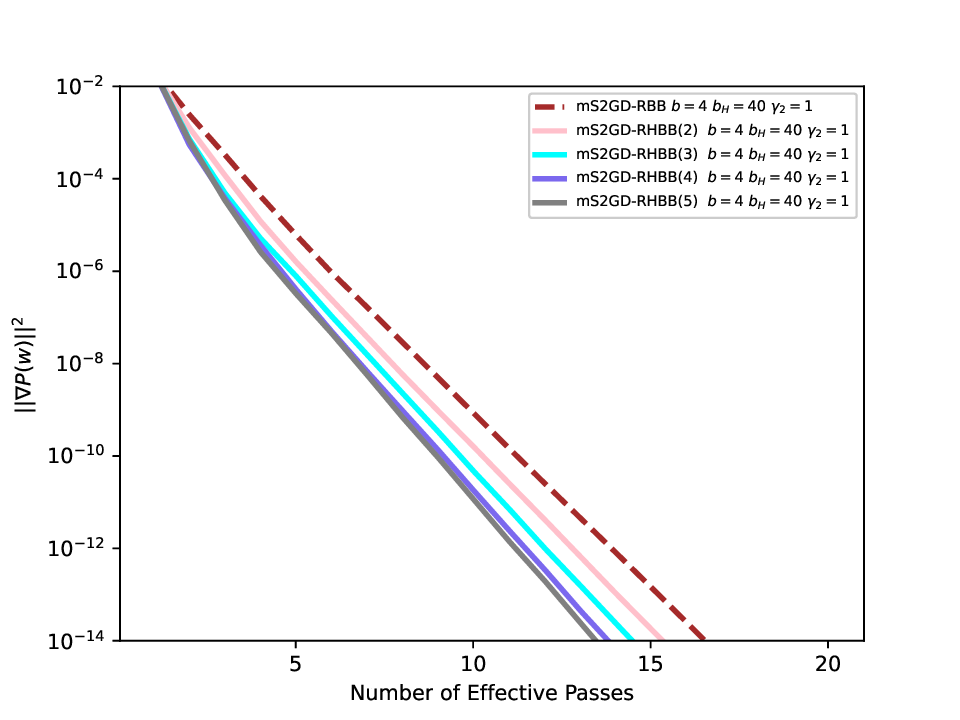}}
  		\subfigure[ijcnn1]
  		{\includegraphics[width=0.327\textwidth]{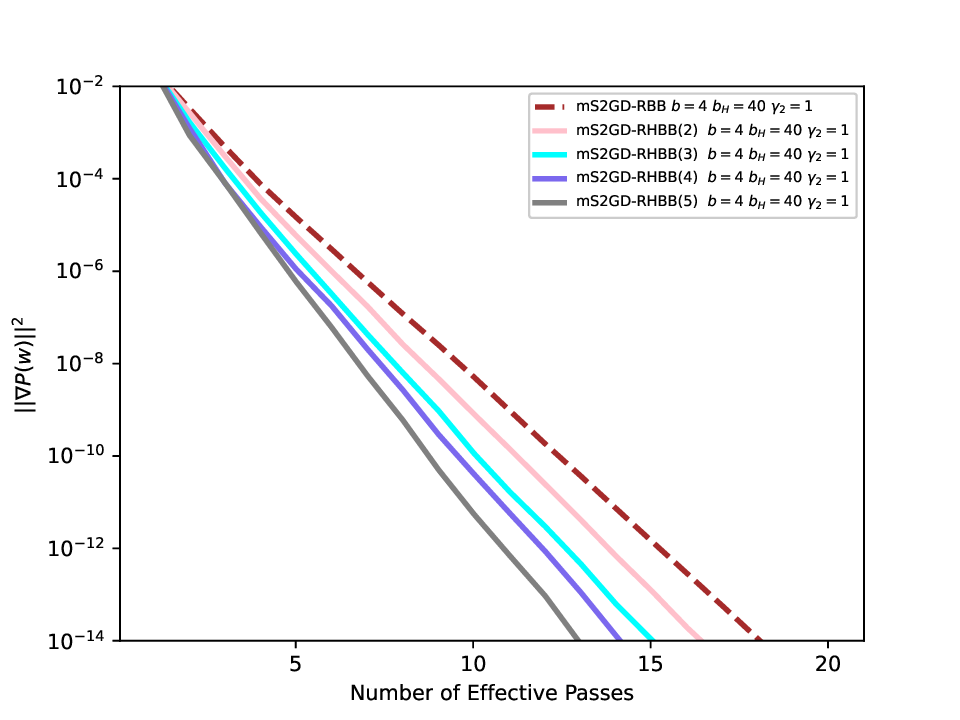}}
  		\subfigure[covtype]
  		{\includegraphics[width=0.327\textwidth]{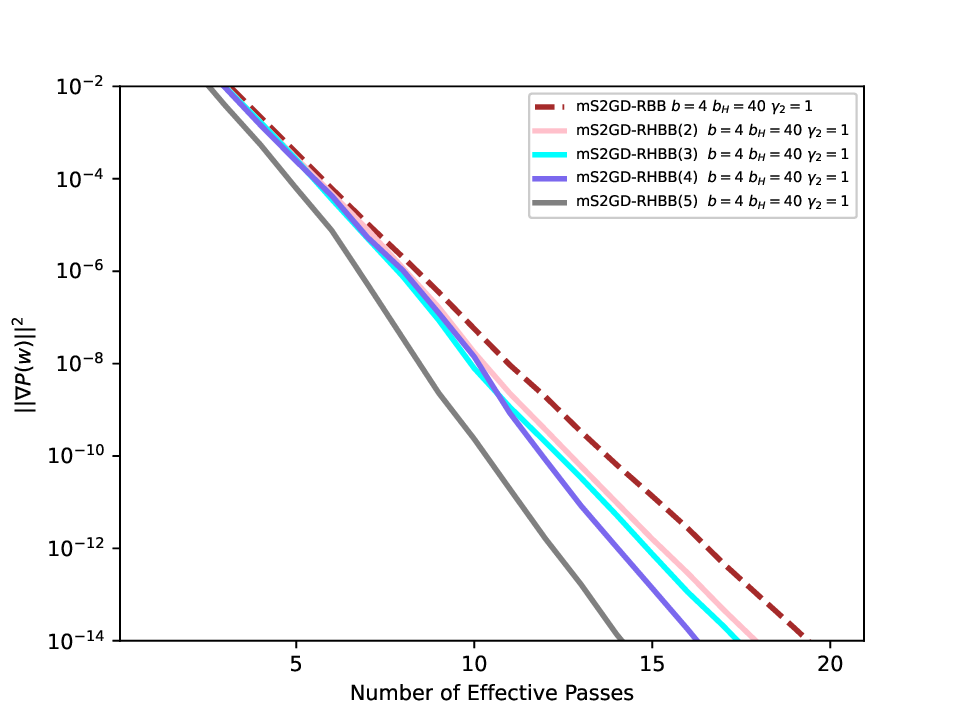}}
  		\subfigure[phishing]
  		{\includegraphics[width=0.327\textwidth]{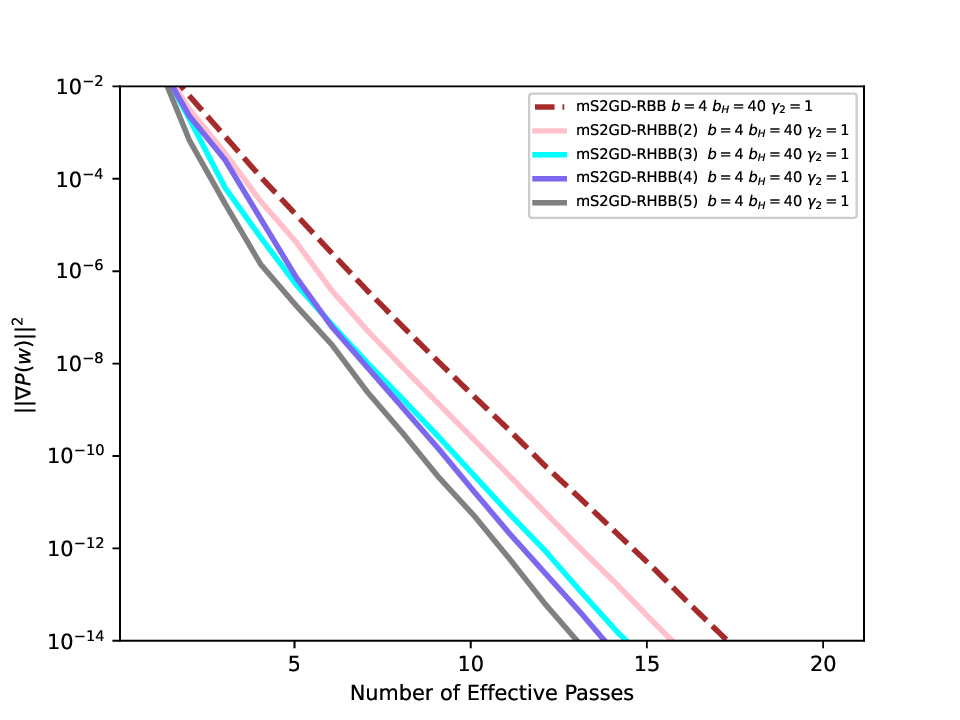}}
  		\subfigure[mushrooms]
  		{\includegraphics[width=0.327\textwidth]{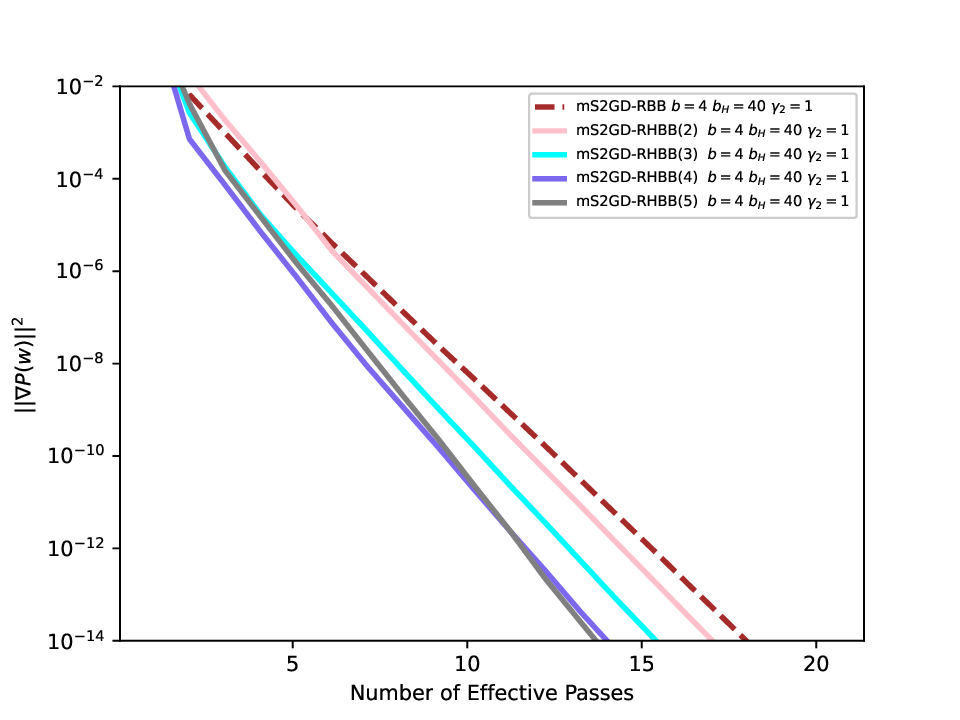}}
  		\caption{\footnotesize Comparisons of mS2GD-RHBB (dash line) and mS2GD-RBB (solid lines) with an $\alpha$ from $\{2,3,4,5\}$.}
  		\label{fig3}
  	\end{figure*}

\begin{figure*}[htbp]
  		\centering
  		\subfigure[a8a]
  		{\includegraphics[width=0.327\textwidth]{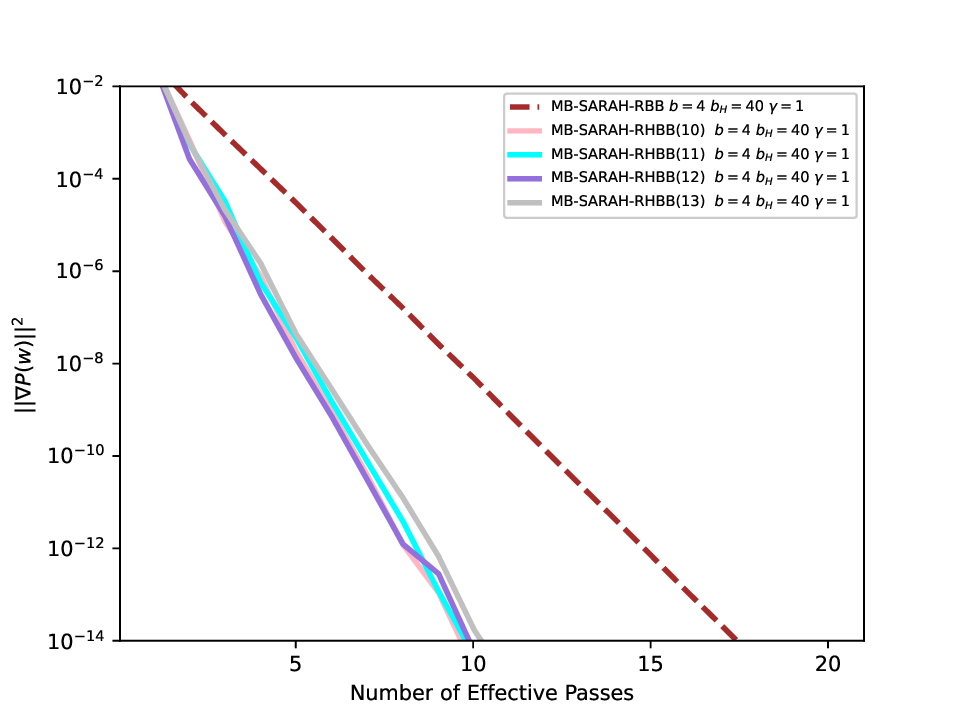}}
  		\subfigure[w8a]
  		{\includegraphics[width=0.327\textwidth]{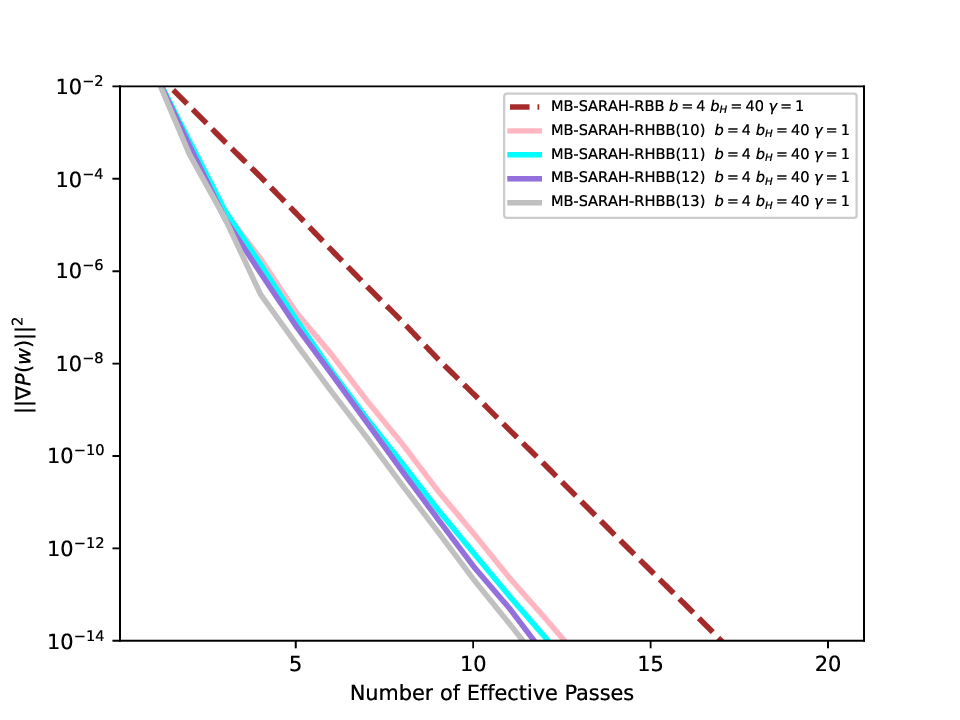}}
  		\subfigure[ijcnn1]
  		{\includegraphics[width=0.327\textwidth]{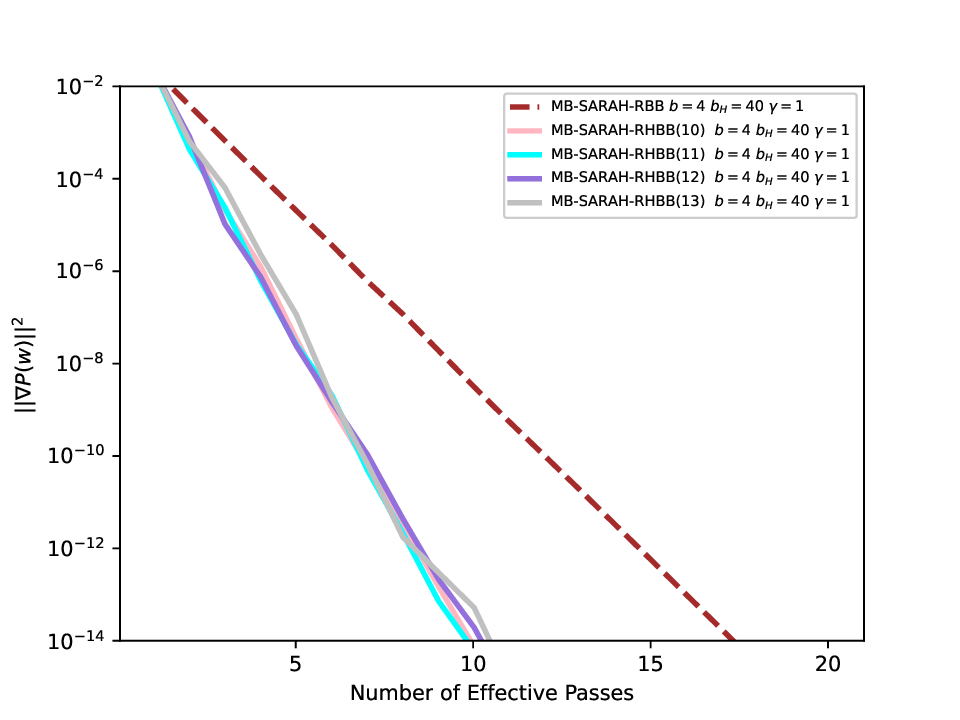}}
  		\subfigure[covtype]
  		{\includegraphics[width=0.327\textwidth]{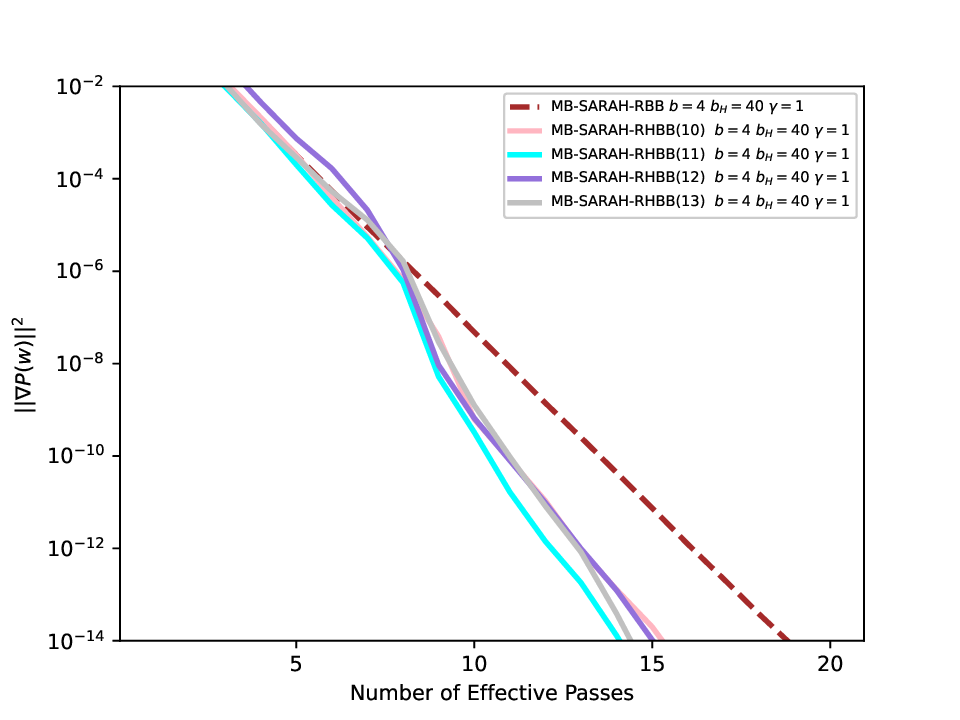}}
  		\subfigure[phishing]
  		{\includegraphics[width=0.327\textwidth]{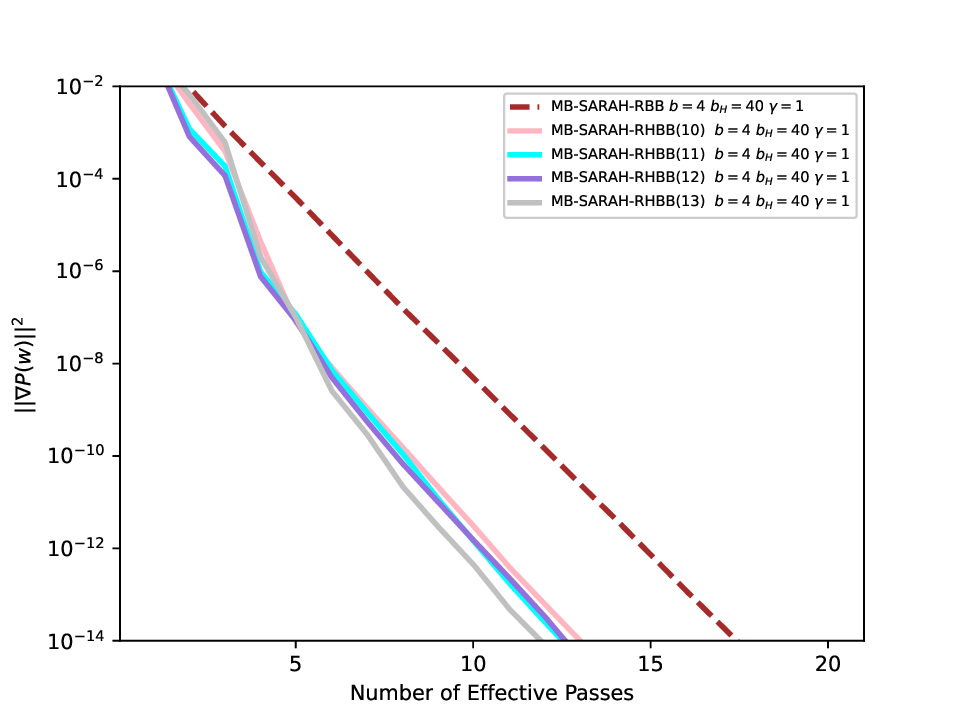}}
  		\subfigure[mushrooms]
  		{\includegraphics[width=0.327\textwidth]{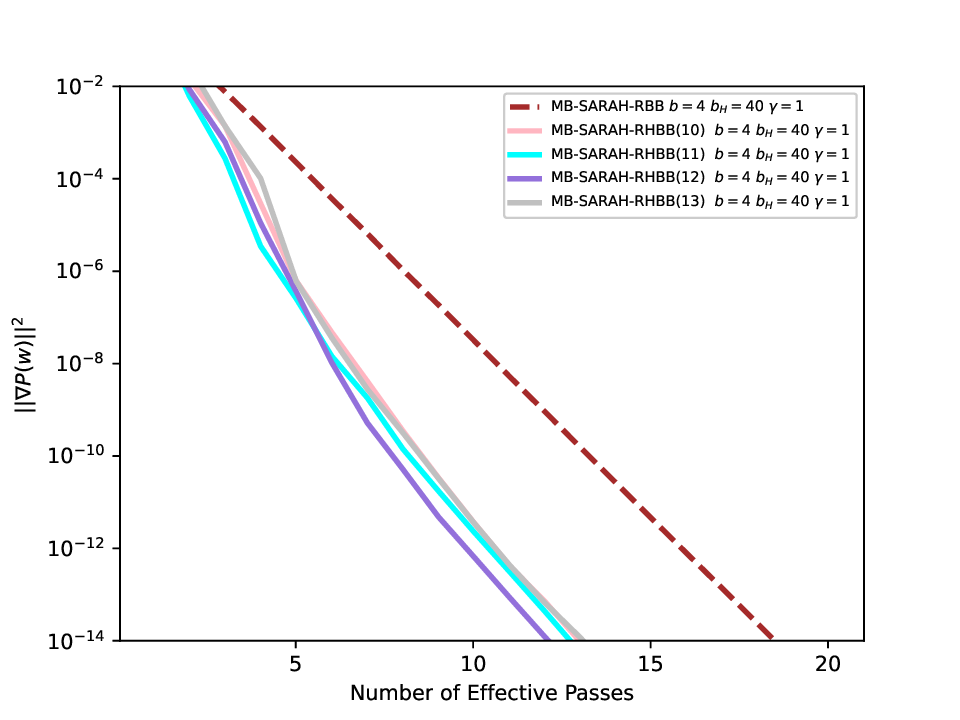}}
  		\caption{\footnotesize Comparisons of MB-SARAH-RHBB (dash line) and MB-SARAH-RBB (solid lines) with an aggressive $\alpha$ from $\{10,11,12,13\}$.}
  		\label{fig4}
  	\end{figure*}

\begin{figure*}[htbp]
  		\centering
  		\subfigure[a8a]
  		{\includegraphics[width=0.327\textwidth]{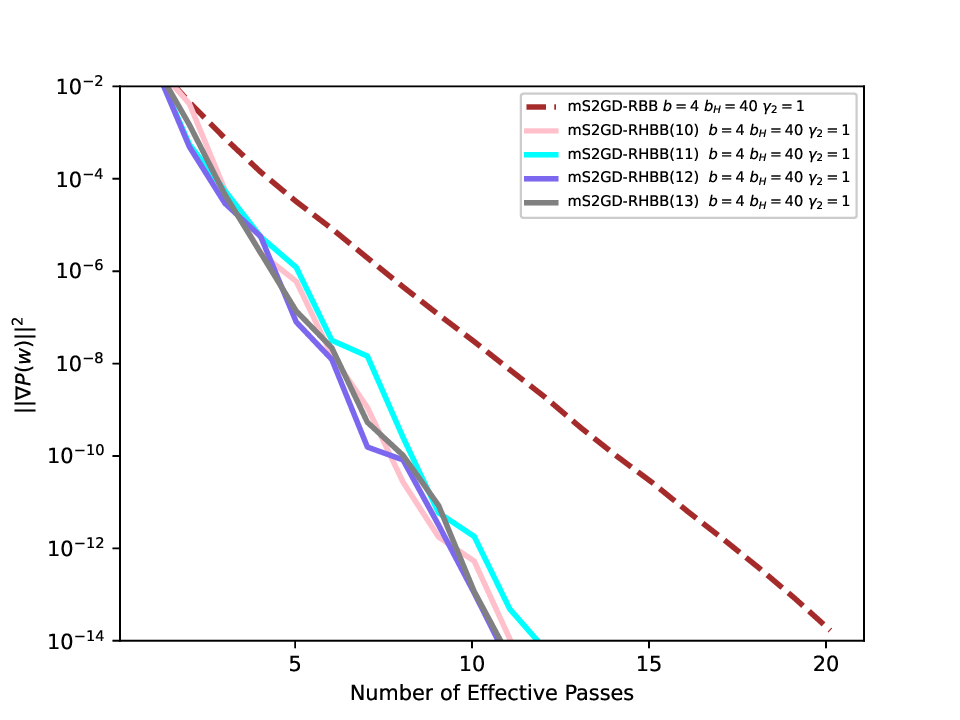}}
  		\subfigure[w8a]
  		{\includegraphics[width=0.327\textwidth]{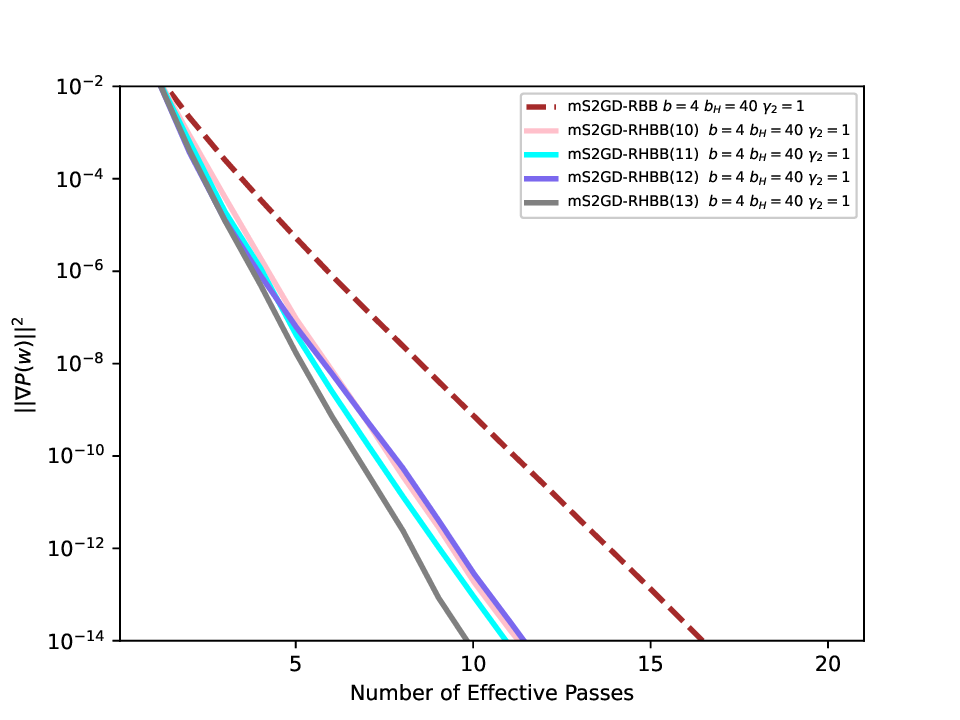}}
  		\subfigure[ijcnn1]
  		{\includegraphics[width=0.327\textwidth]{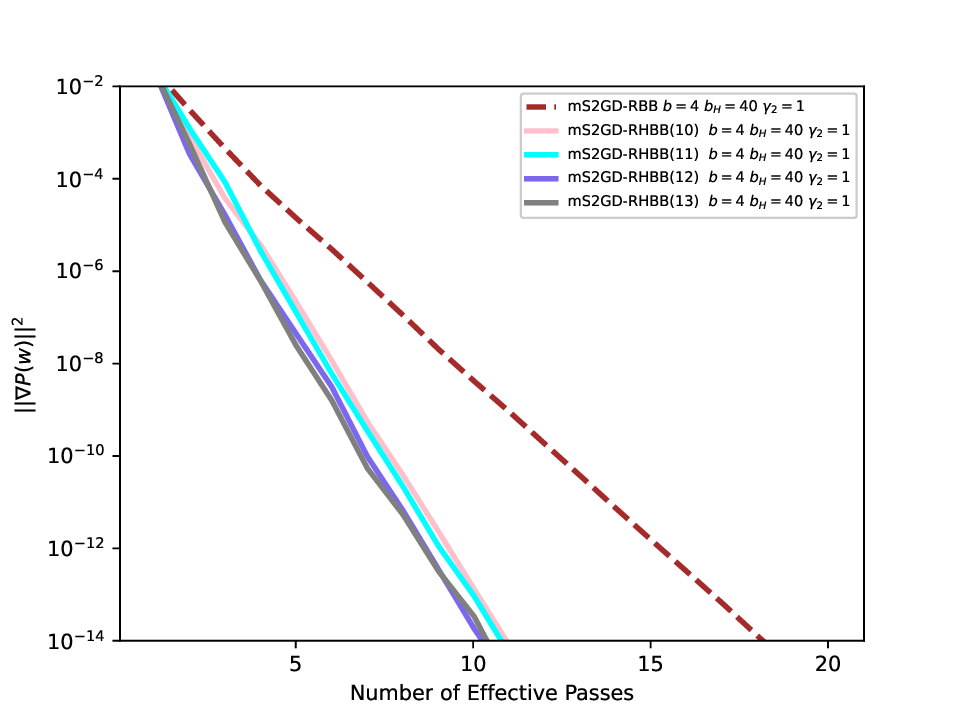}}
  		\subfigure[covtype]
  		{\includegraphics[width=0.327\textwidth]{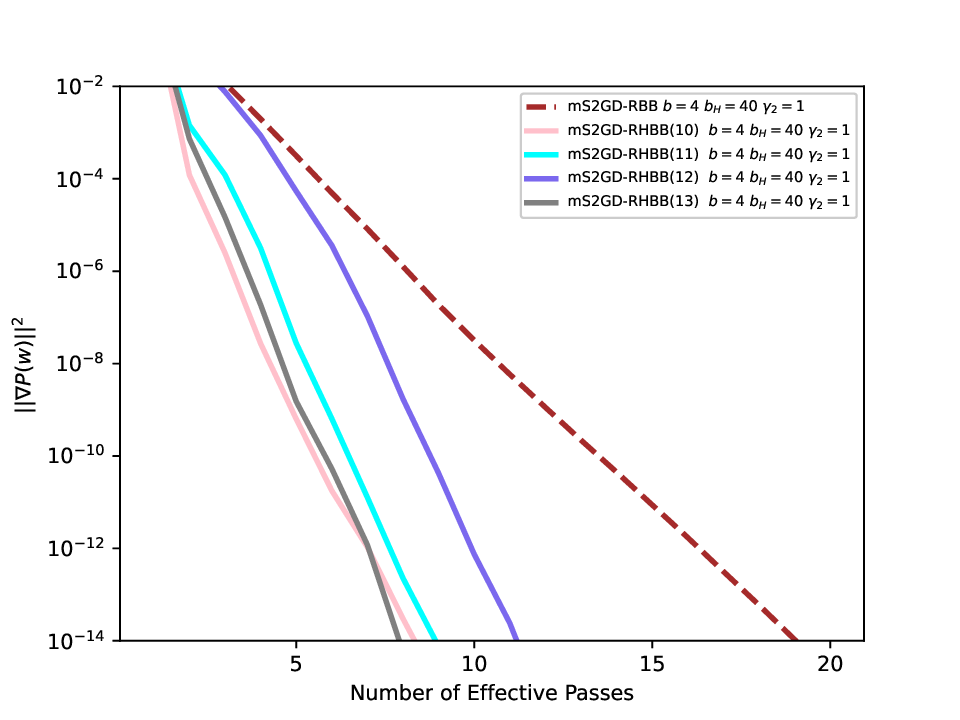}}
  		\subfigure[phishing]
  		{\includegraphics[width=0.327\textwidth]{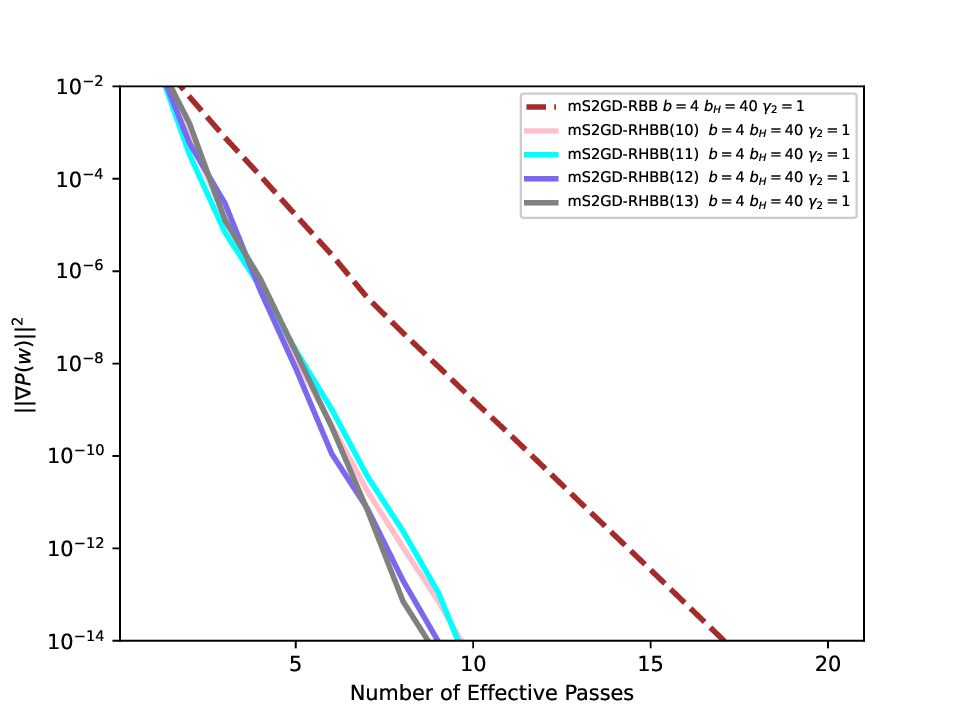}}
  		\subfigure[mushrooms]
  		{\includegraphics[width=0.327\textwidth]{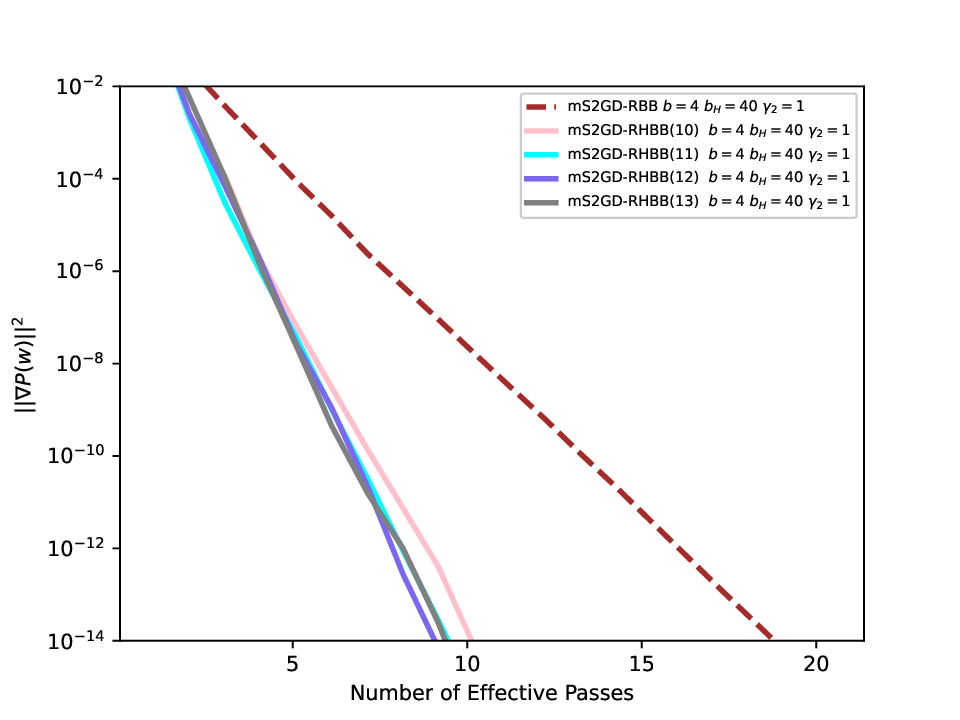}}
  		\caption{\footnotesize Comparisons of mS2GD-RHBB (dash line) and mS2GD-RBB (solid lines) with an aggressive $\alpha$ from $\{10,11,12,13\}$.}
  		\label{fig5}
  	\end{figure*}

\begin{figure*}[htbp]
  		\centering
  		\subfigure[unvaried]
  		{\includegraphics[width=0.4\textwidth]{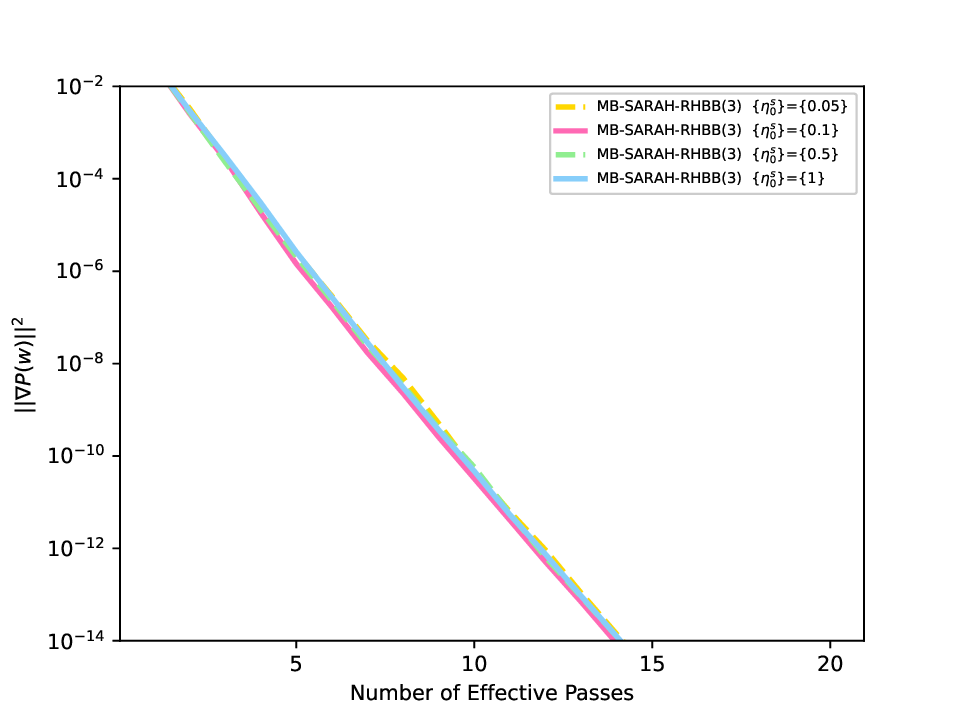}}
  		\subfigure[unvaried]
  		{\includegraphics[width=0.4\textwidth]{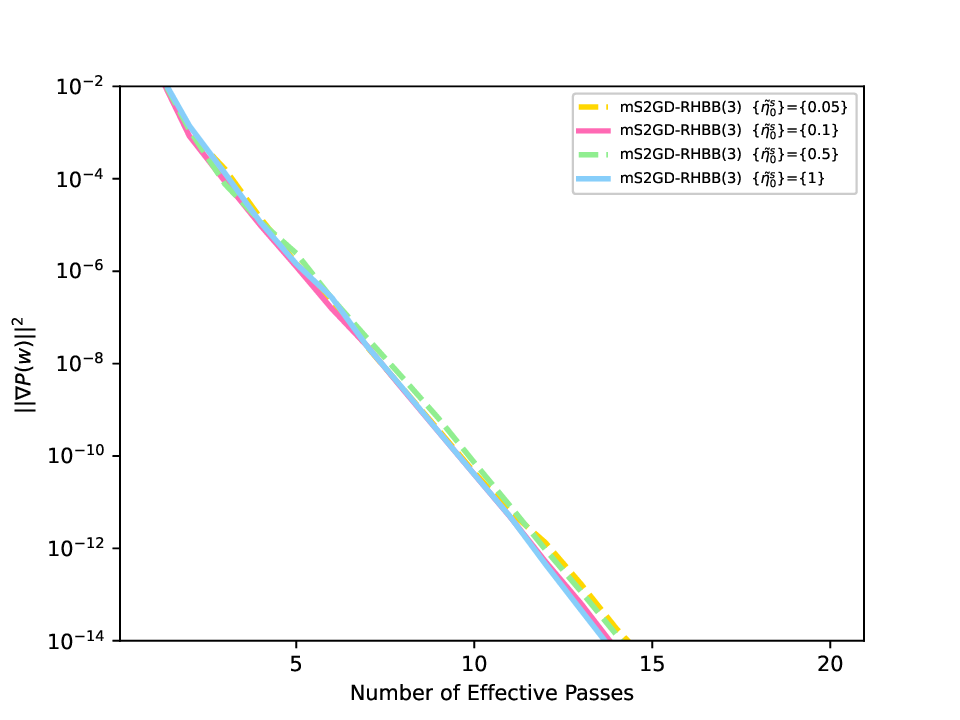}}
  		
  		\subfigure[mix]
  		{\includegraphics[width=0.4\textwidth]{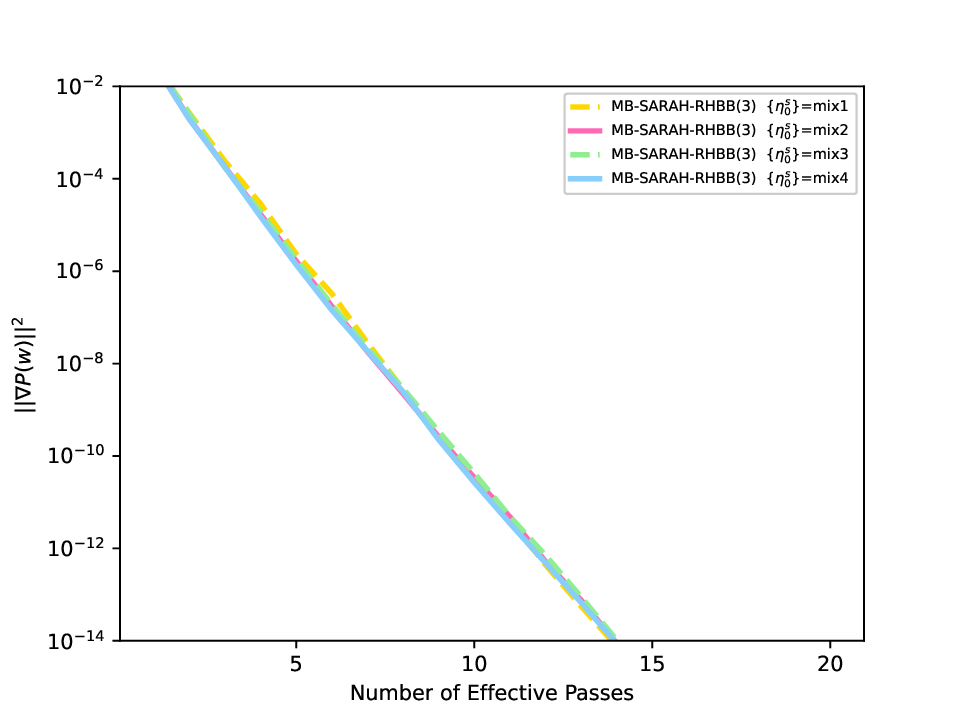}}
  		\subfigure[mix]
  		{\includegraphics[width=0.4\textwidth]{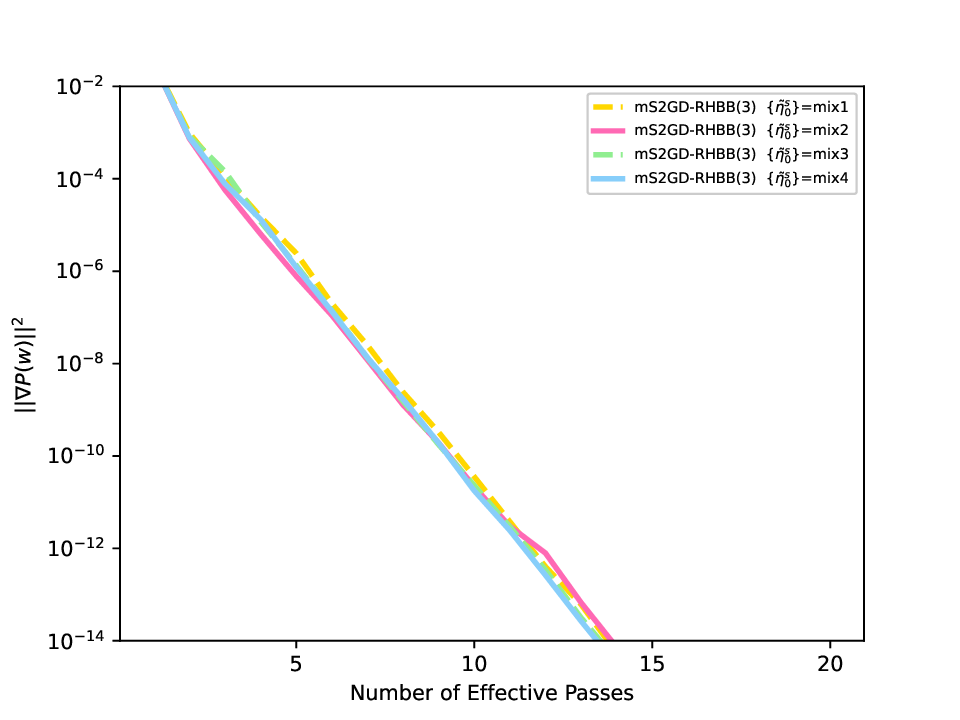}}
  		\caption{\footnotesize The performance of MB-SARAH-RHBB and mS2GD-RHBB under different initial step sizes on $a8a$. (a) MB-SARAH-RHBB under the unvaried initial step sizes; (b) mS2GD-RHBB under the unvaried initial step sizes; (c) MB-SARAH-RHBB under the mix initial step sizes; (d) mS2GD-RHBB under the mix initial step sizes.}
  		\label{fig6}
  	\end{figure*}

\begin{figure*}[htbp]
  		\centering
  		\subfigure[a8a]
  		{\includegraphics[width=0.327\textwidth]{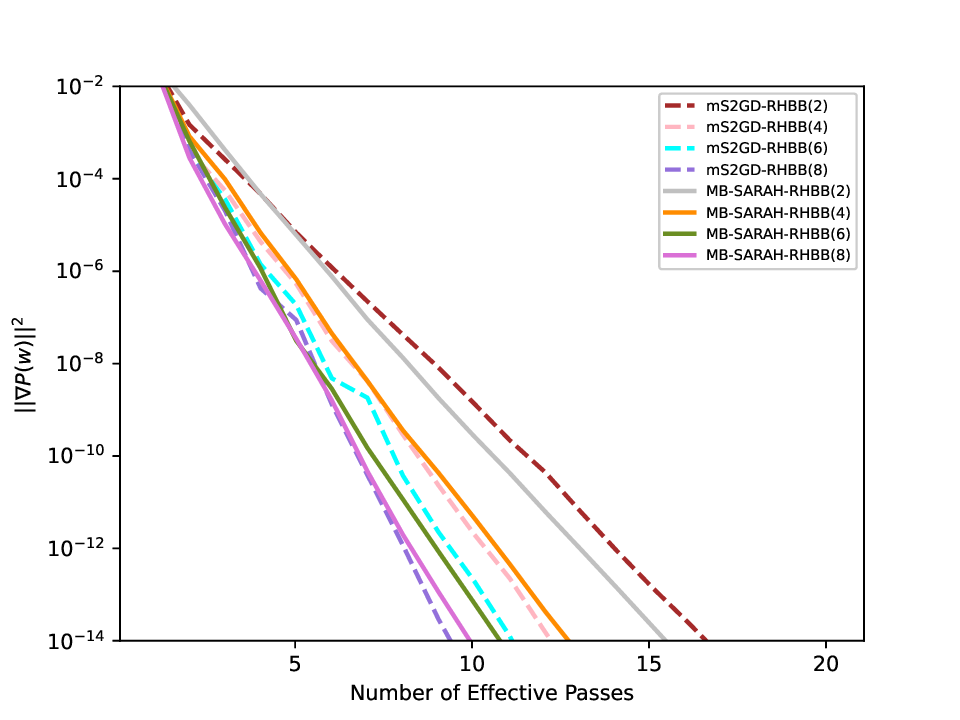}}
  		\subfigure[w8a]
  		{\includegraphics[width=0.327\textwidth]{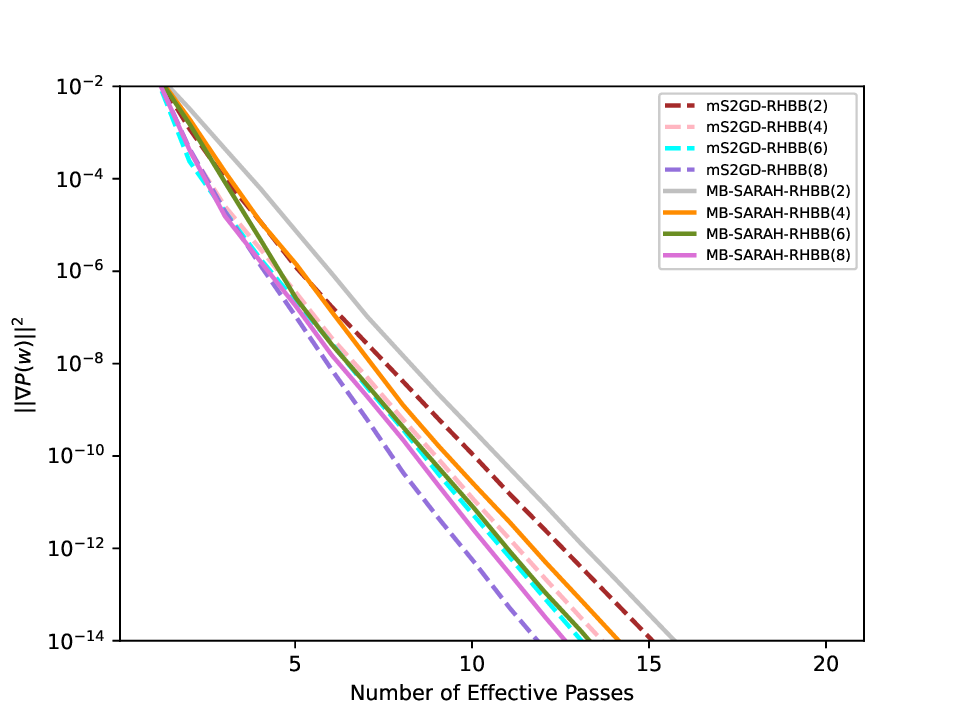}}
  		\subfigure[ijcnn1]
  		{\includegraphics[width=0.327\textwidth]{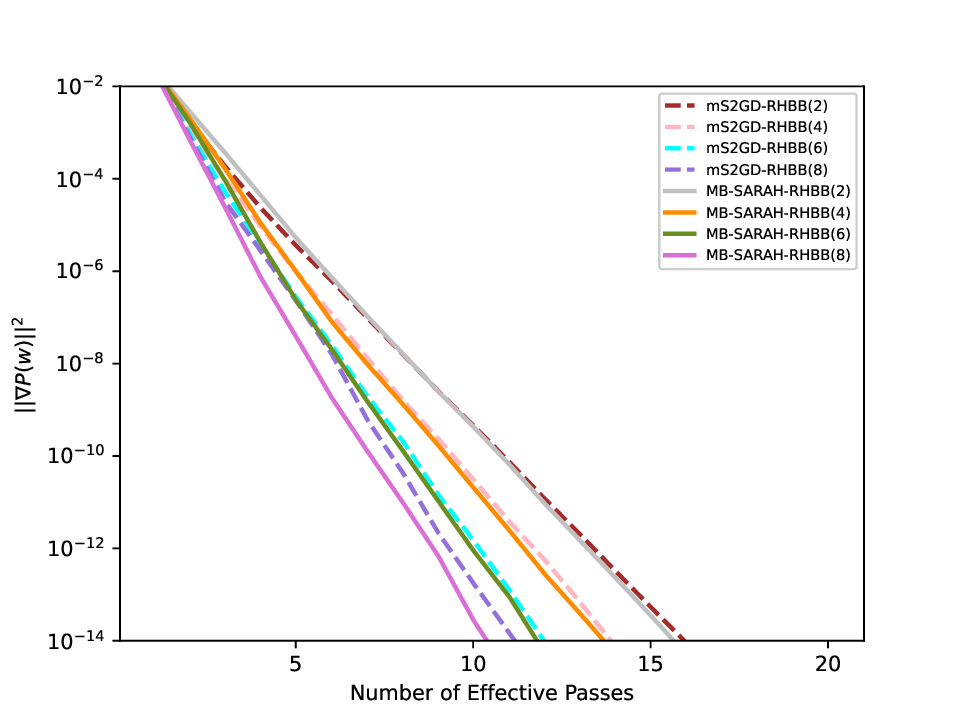}}
  		\subfigure[covtype]
  		{\includegraphics[width=0.327\textwidth]{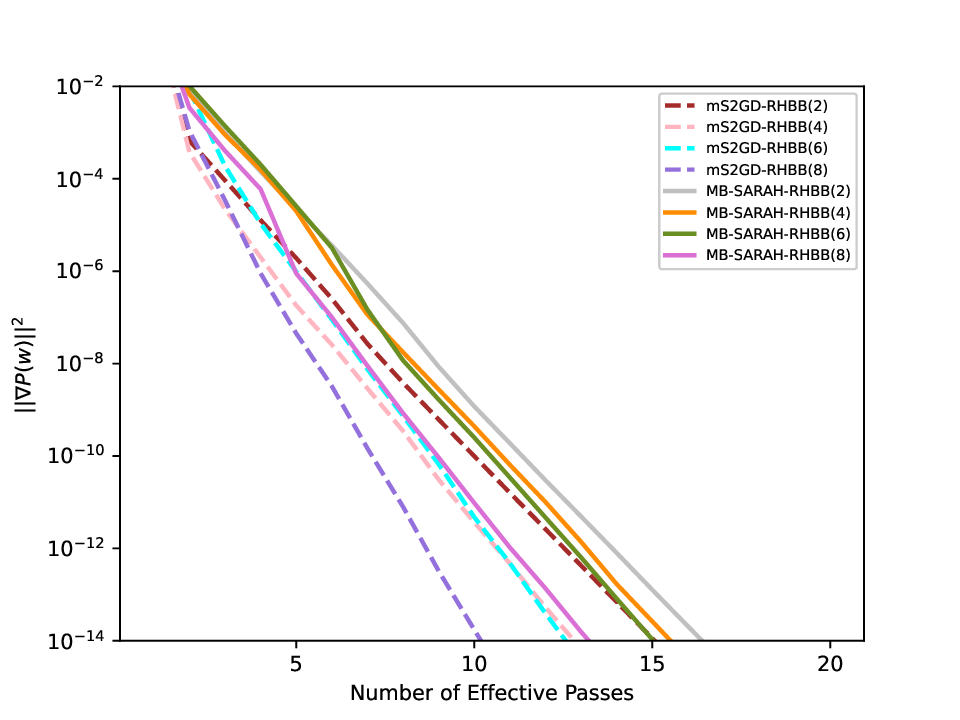}}
  		\subfigure[phishing]
  		{\includegraphics[width=0.327\textwidth]{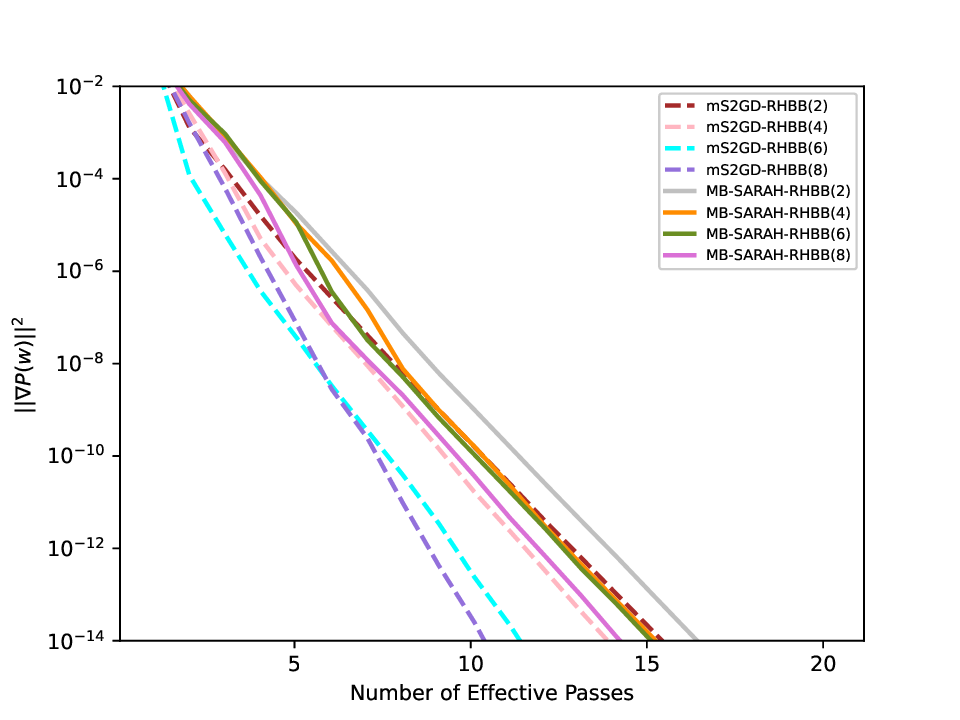}}
  		\subfigure[mushrooms]
  		{\includegraphics[width=0.327\textwidth]{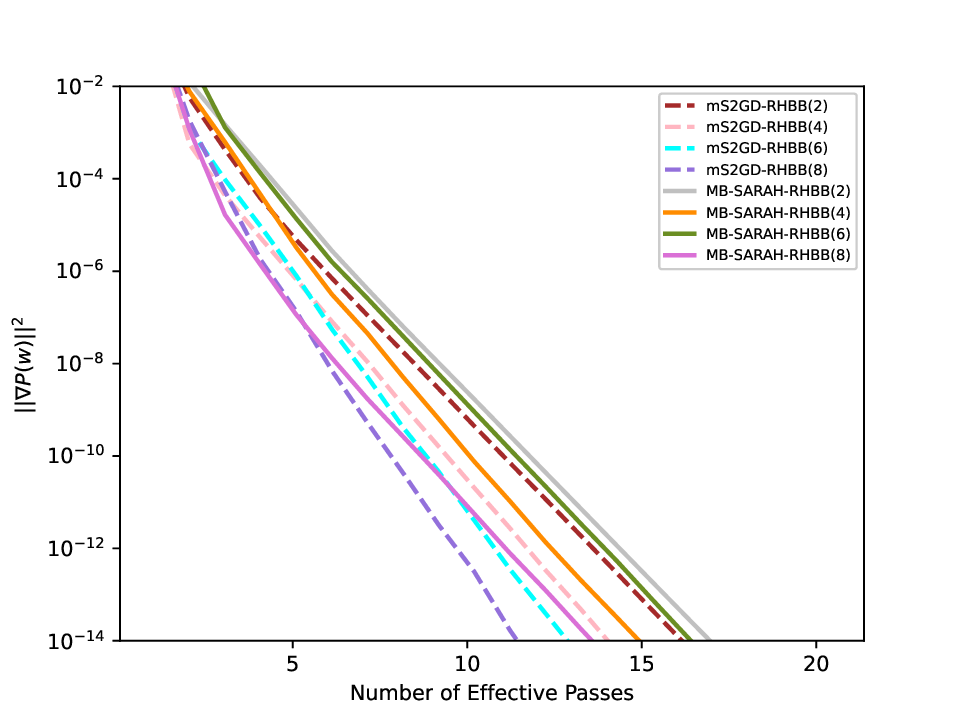}}
  		\caption{\footnotesize Comparisons of MB-SARAH-RHBB (solid lines) and mS2GD-RHBB (dash lines).}
  		\label{fig7}
  	\end{figure*}

\begin{figure*}[htbp]
  		\centering
  		\subfigure[a8a]
  		{\includegraphics[width=0.327\textwidth]{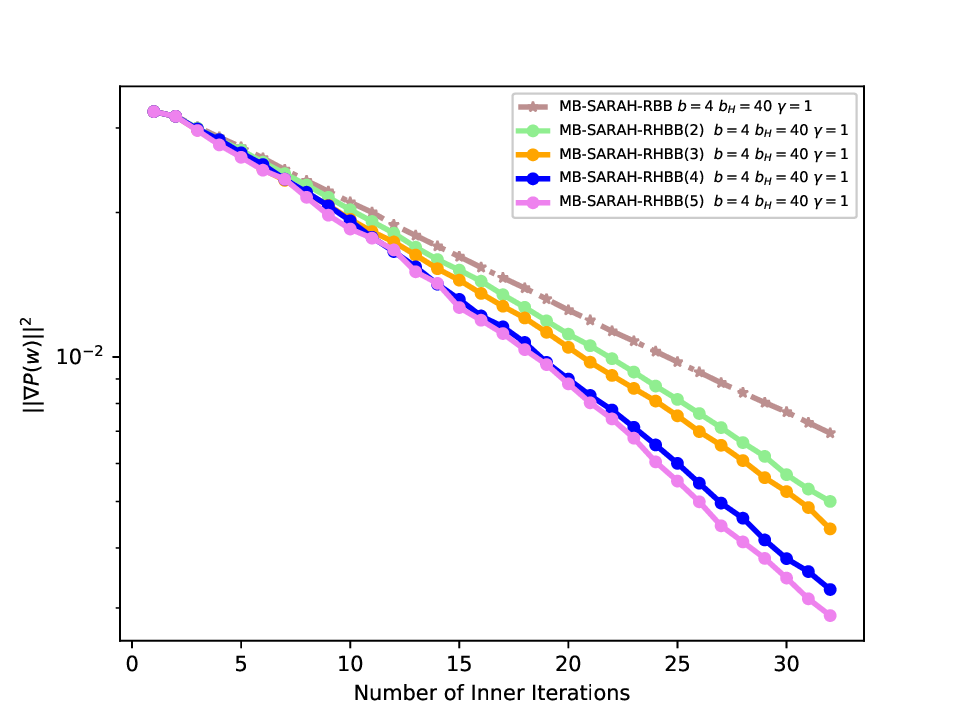}}
  		\subfigure[w8a]
  		{\includegraphics[width=0.327\textwidth]{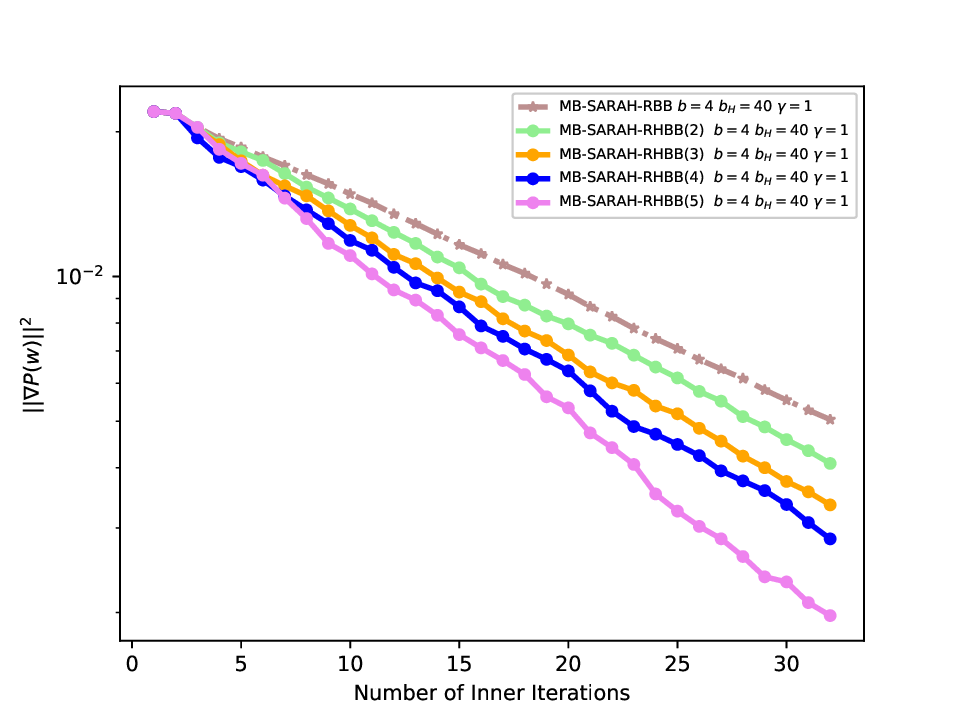}}
  		\subfigure[ijcnn1]
  		{\includegraphics[width=0.327\textwidth]{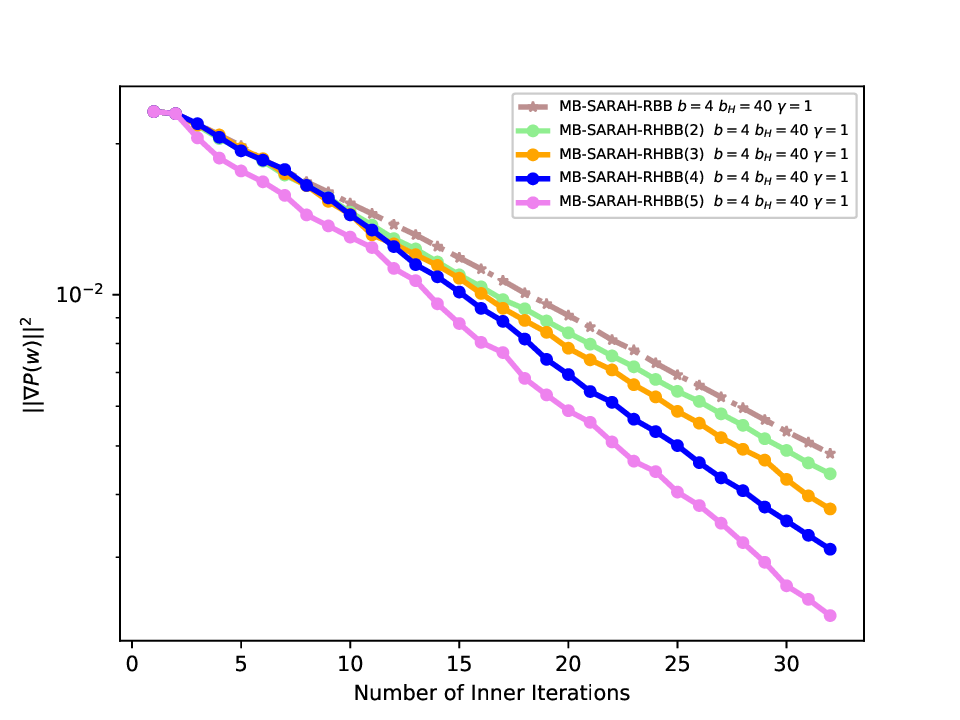}}
  		\subfigure[covtype]
  		{\includegraphics[width=0.327\textwidth]{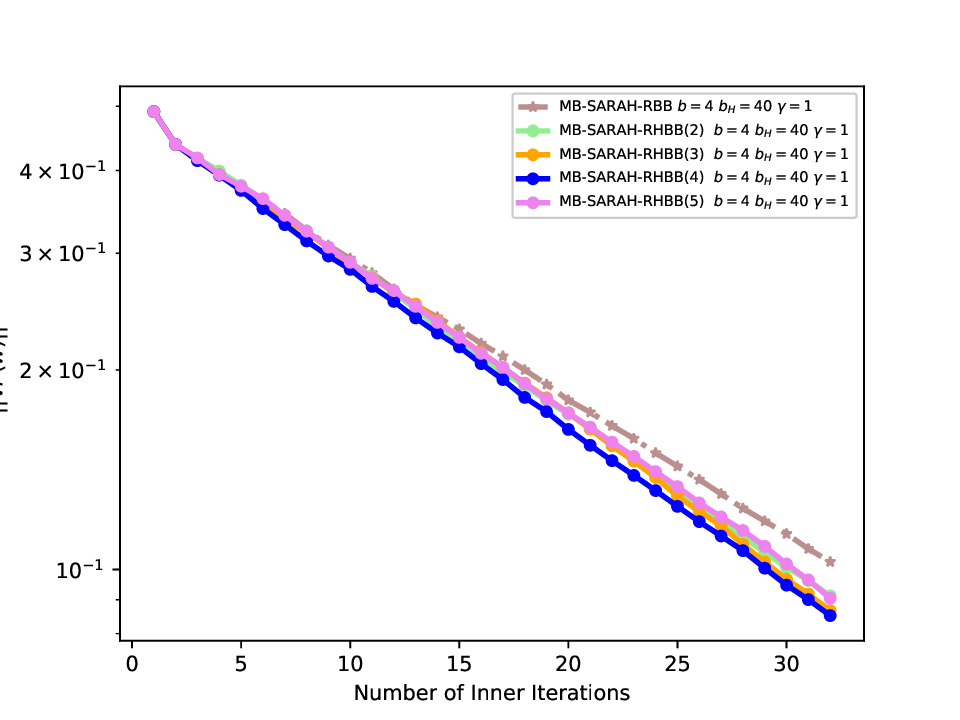}}
  		\subfigure[phishing]
  		{\includegraphics[width=0.327\textwidth]{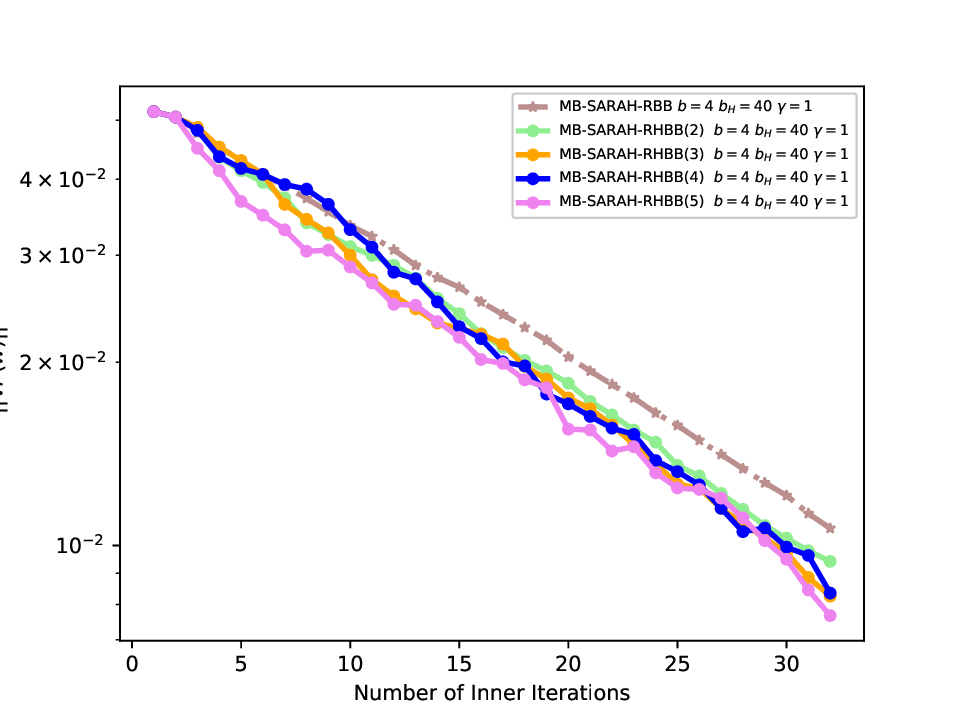}}
  		\subfigure[mushrooms]
  		{\includegraphics[width=0.327\textwidth]{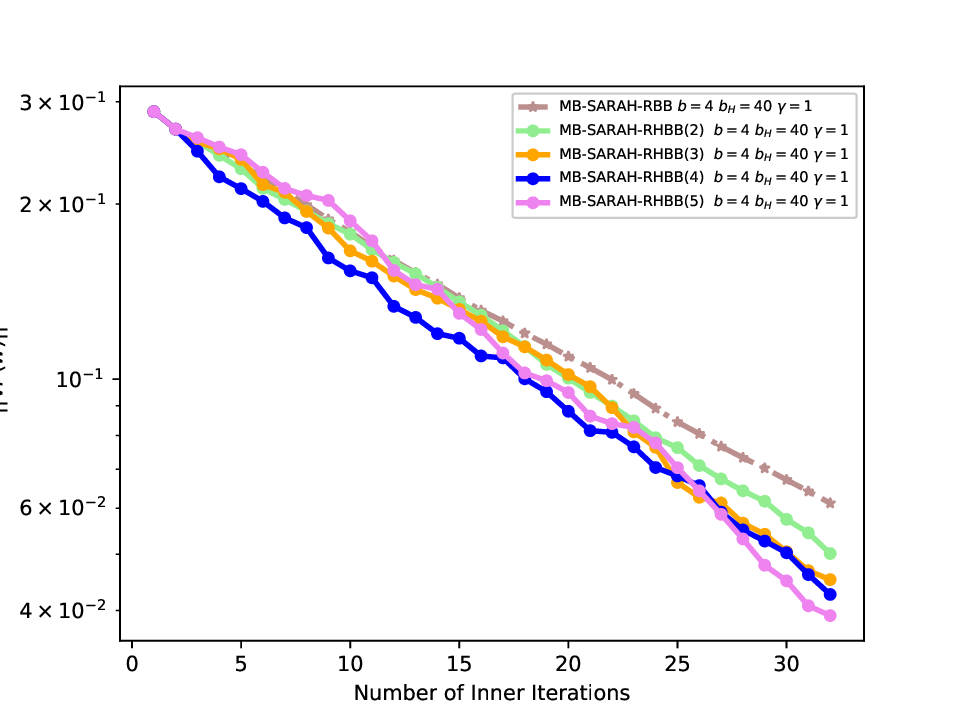}}
  		\caption{\footnotesize The performance of MB-SARAH-IN-RHBB and MB-SARAH-IN-RBB.}
  		\label{fig8}
  	\end{figure*}

\begin{figure*}[htbp]
  		\centering
  		\subfigure[a8a]
  		{\includegraphics[width=0.327\textwidth]{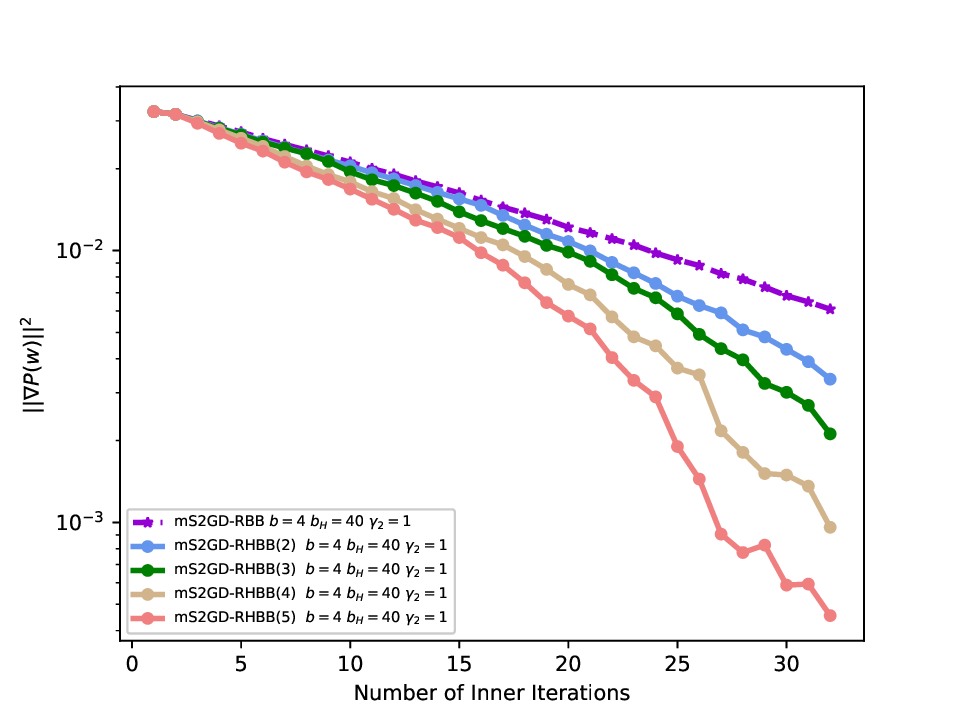}}
  		\subfigure[w8a]
  		{\includegraphics[width=0.327\textwidth]{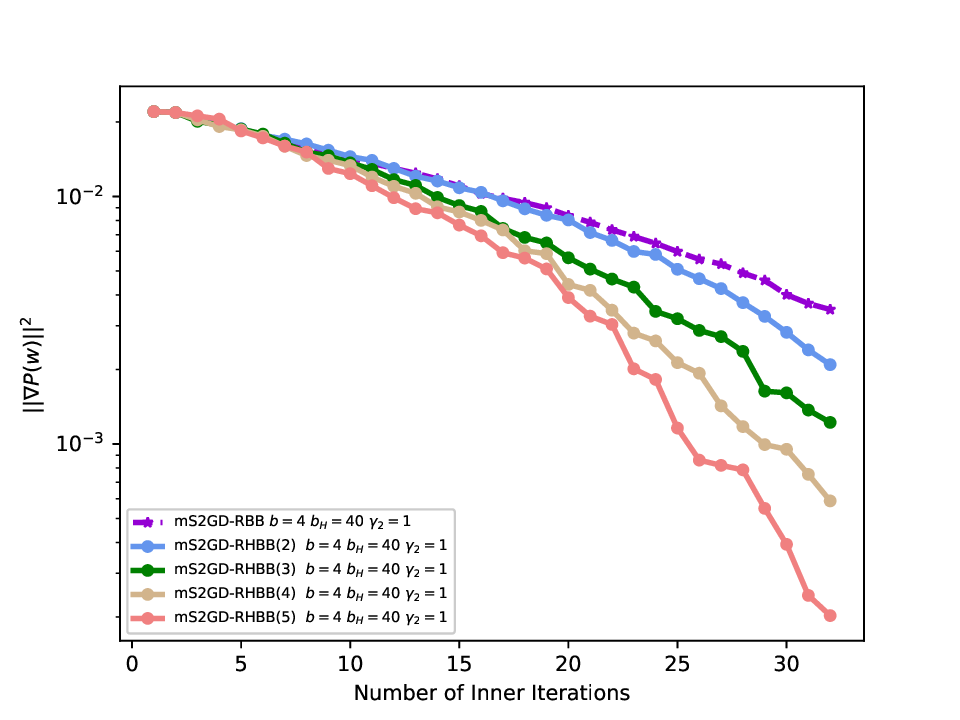}}
  		\subfigure[ijcnn1]
  		{\includegraphics[width=0.327\textwidth]{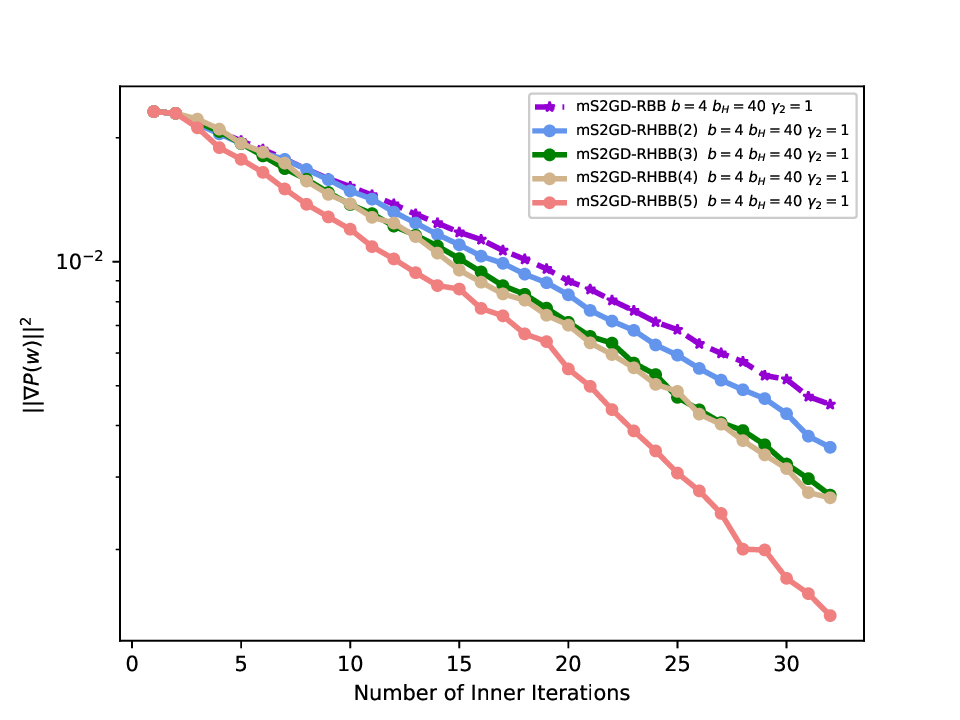}}
  		\subfigure[covtype]
  		{\includegraphics[width=0.327\textwidth]{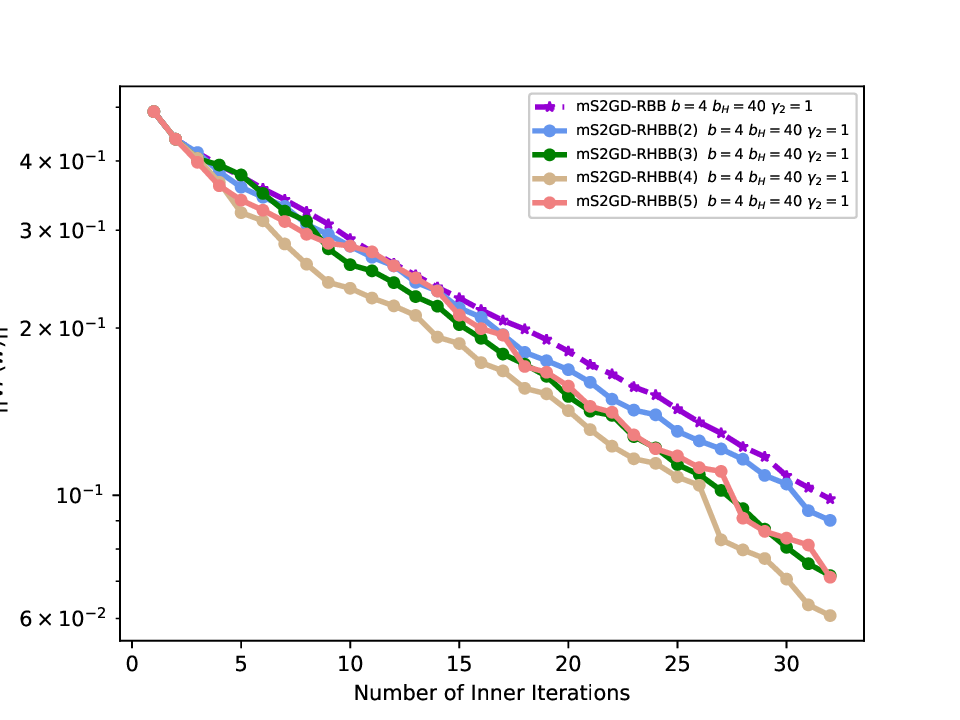}}
  		\subfigure[phishing]
  		{\includegraphics[width=0.327\textwidth]{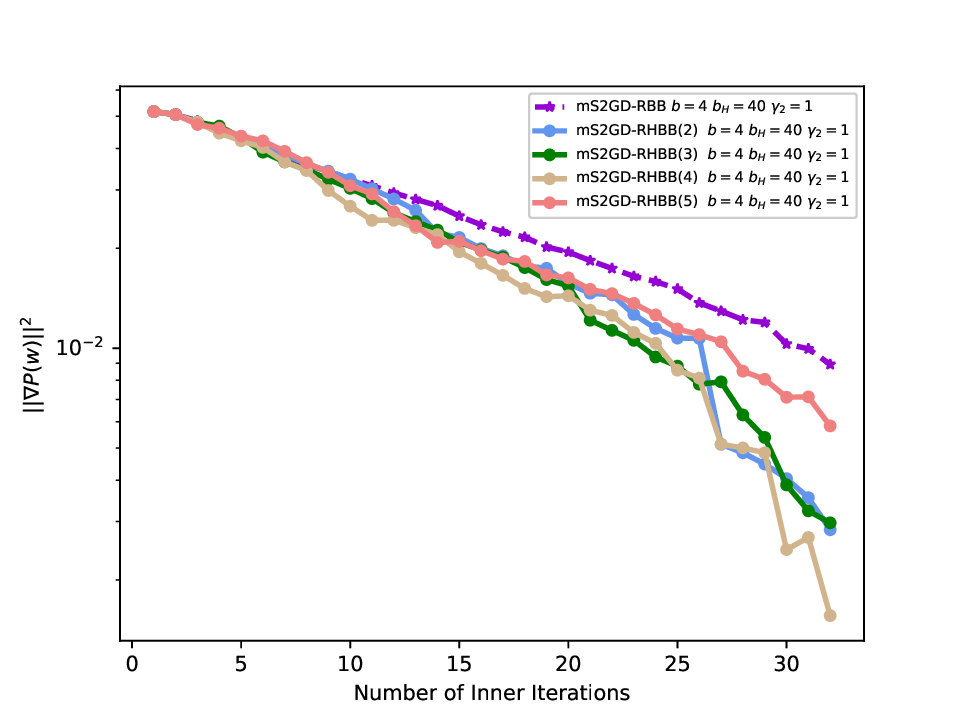}}
  		\subfigure[mushrooms]
  		{\includegraphics[width=0.327\textwidth]{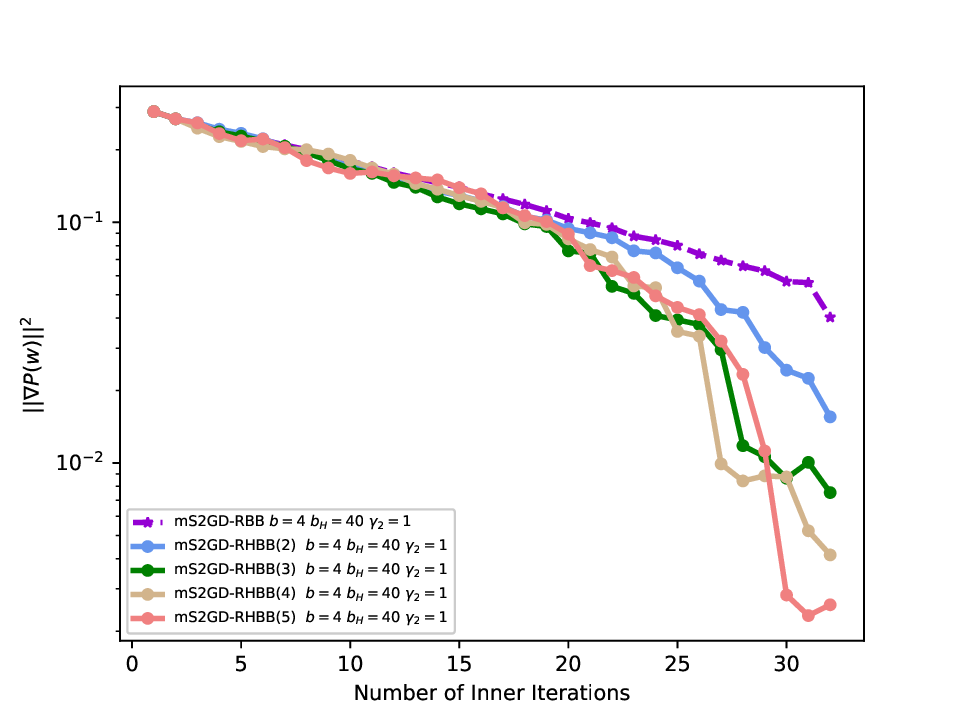}}
  		\caption{\footnotesize The performance of mS2GD-IN-RHBB and mS2GD-IN-RBB.}
  		\label{fig9}
  	\end{figure*}

  	Figs. \ref{fig2} - \ref{fig9} show the results of the properties of MB-SARAH-RHBB and mS2GD-RHBB. In all sub-figures, the horizontal axis denotes the number of effective passes, and the vertical axis represents the Euclidean norm of $\nabla P(\cdot)$.
  	
  	In Figs. \ref{fig2} - \ref{fig5}, we analyze the unvaried hedge effect by increasing the value of $\alpha$ either gradually or drastically. From Figs. \ref{fig2}, \ref{fig3}, we observe that the practical speeds of MB-SARAH-RHBB and mS2GD-RHBB are continuously improving with increasing $\alpha$ from $\{2, 3, 4, 5\}$. The following Figs. \ref{fig4}, \ref{fig5} indicate that the performance of the algorithms reaches a plateau as $\alpha$ becomes more aggressive from $\{10, 11, 12, 13\}$. Note that MB-SARAH-RHBB and mS2GD-RHBB outperform the original RBB-type algorithms consistently on all data sets.
  	
  	In Fig. \ref{fig6}, we analyze the constant step size sequences $\{\eta_0^s\}$, $\{\tilde{\eta}_0^s\}$ ($s\geq1$) that are applied in the deterministic steps. For reliability, we randomly tossed out a value of $\alpha=3$ to run the algorithms. In fact, we pick four unvaried sequences $\{0.05\}$, $\{0.1\}$, $\{0.5\}$, $\{1\}$, and four mingle sequences of ascending $mix1 = \{0.05, 0.1, 0.5, 1, ...\}$, descending $mix2 = \{1, 0.5, 0.1, 0.05, ...\}$, disordered $mix3 = \{0.5, 1, 0.05, 0.1, ...\}$ and disordered $mix4 = \{1, 0.05, 0.1, 0.5, ...\}$, as participants. The practical performance of MB-SARAH-RHBB and mS2GD-RHBB is not influenced, implying $\{\eta_0^s\}$ $\{\tilde{\eta}_0^s\}$ ($s\geq1$) are immaterial but provide sufficient curvature for the following RHBB/RHBB+ calculations.
  	
  	Fig. \ref{fig7} exhibits the comparisons between MB-SARAH-RHBB and mS2GD-RHBB, where we devise multiple comparison levels of $\alpha=2, 3, 4, 5$. MB-SARAH-RHBB outperforms mS2GD-RHBB on $ijcnn1$, mS2GD-RHBB performs better on $covtype$, $phishing$, $mushrooms$, $w8a$, and they performed equally well on $a8a$. In most cases, mS2GD-RHBB delivers superior performance.
  	
  	In order to analyze their properties within the inner loops, we discard the outer epoch and explore the performance of MB-SARAH-IN-RHBB and mS2GD-IN-RHBB in Figs. \ref{fig8}, \ref{fig9}. It's evident that mS2GD is more susceptible to the RHBB step sizes in the inner loops.
  	
\begin{figure*}[htbp]
  		\centering
  		\subfigure[MB-SARAH-RHBB(3)]
  		{\includegraphics[width=0.4\textwidth]{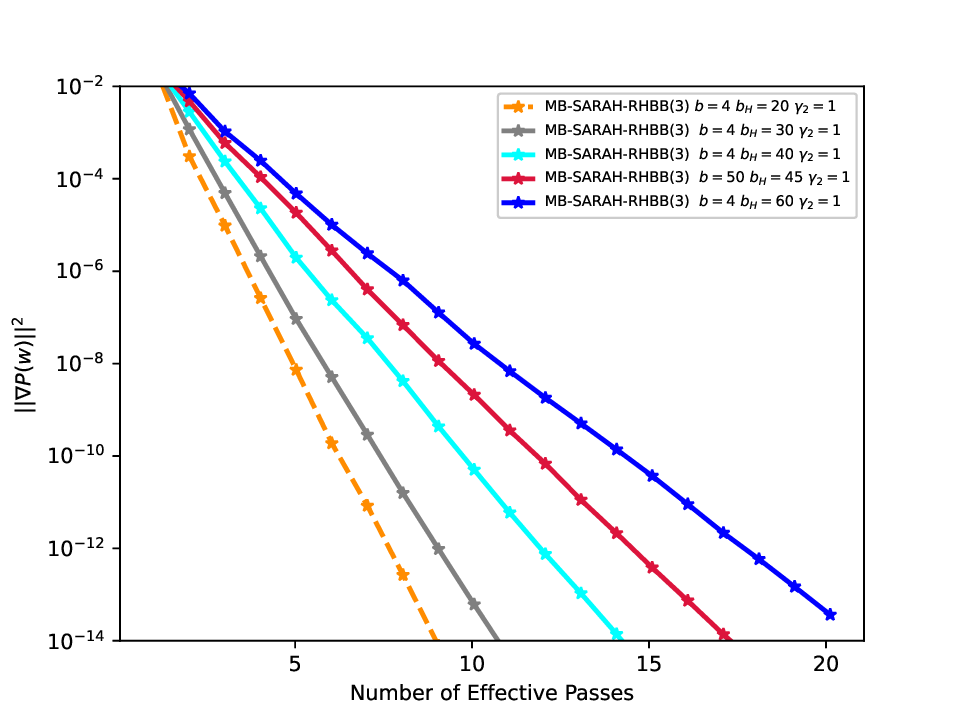}}
  		\subfigure[MB-SARAH-RHBB(6)]
  		{\includegraphics[width=0.4\textwidth]{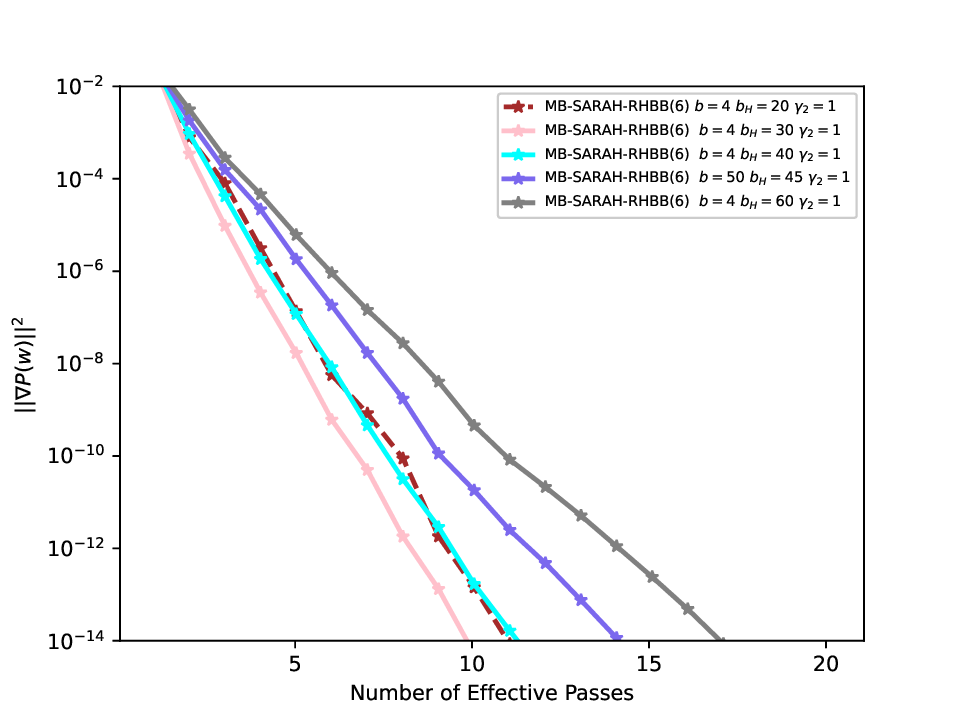}}
  		\caption{\footnotesize The performance of MB-SARAH-RHBB(3) and MB-SARAH-RHBB(6), under different unified $b_H$ on $a8a$.}
  		\label{fig10}
  	\end{figure*}

\begin{figure*}[htbp]
  		\centering
  		\subfigure[mS2GD-RHBB(3)]
  		{\includegraphics[width=0.4\textwidth]{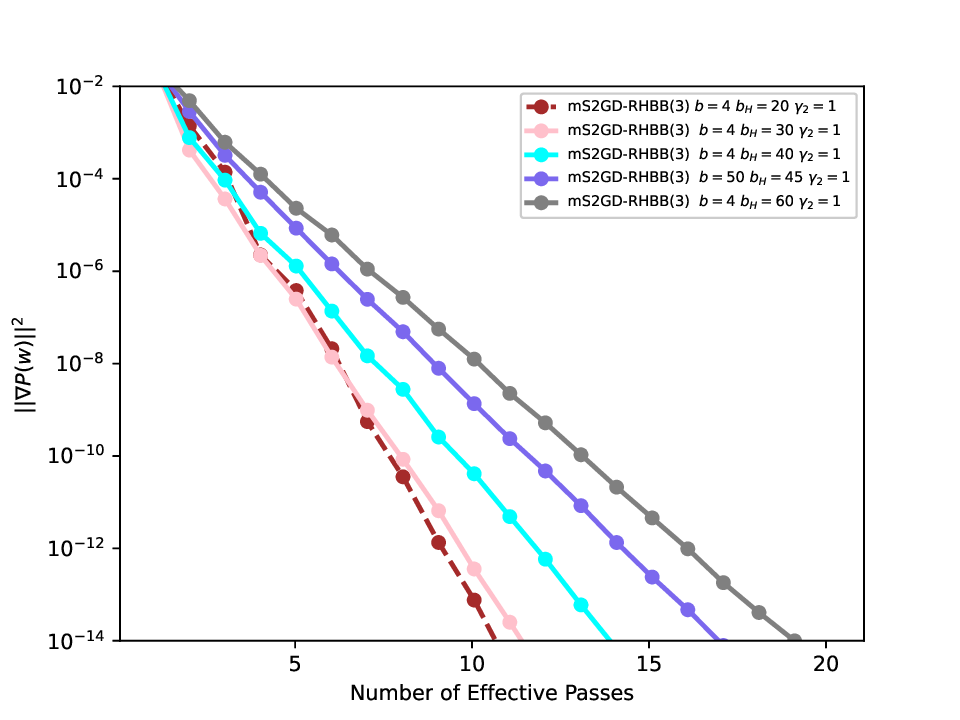}}
  		\subfigure[mS2GD-RHBB(6)]
  		{\includegraphics[width=0.4\textwidth]{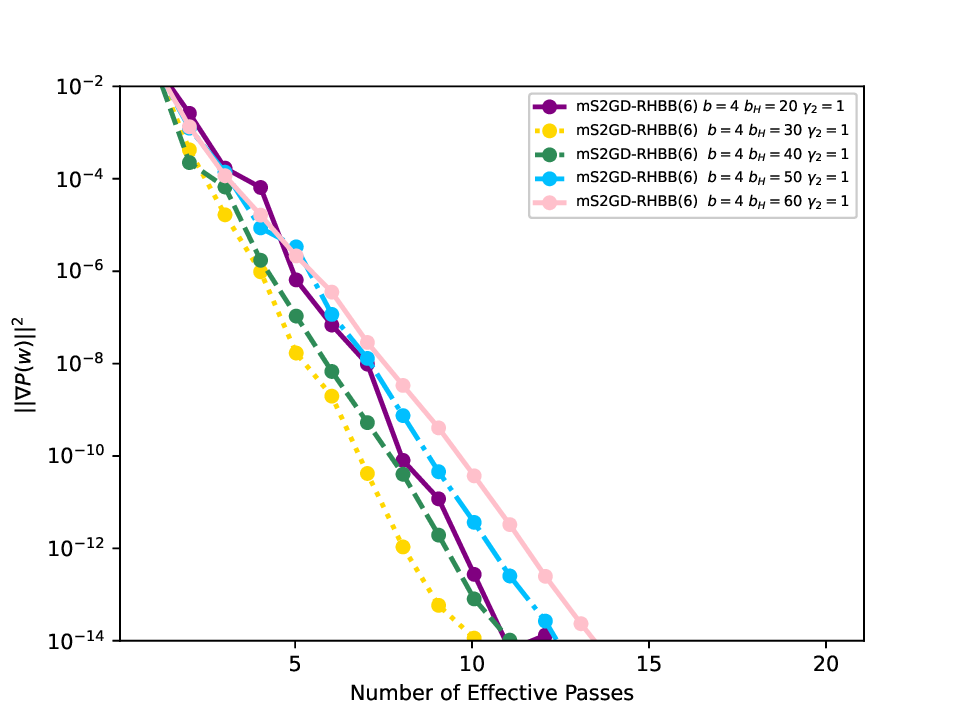}}
  		\caption{\footnotesize The performance of mS2GD-RHBB(3) and mS2GD-RHBB(6), under different unified $b_H$ on $a8a$.}
  		\label{fig11}
  	\end{figure*}

  	In Figs. \ref{fig10} - \ref{fig11}, we evaluate the performance under different unified batch sizes (i.e., the batch correction) of $b_H=20, 30, 40, 50, 60$ on the data set $a8a$. Note that MB-SARAH-RHBB and mS2GD-RHBB are both sensitive to the selection of $b_H$ (i.e., the batch correction).  
  	
  	\subsubsection{Non-adaptive MB-SARAH-RHBB+/mS2GD-RHBB+}
  	\ 
  	\vskip10pt
  	
  	\noindent\textbf{Parametric Settings:} We set $b=4$, the unified $b_H=40$ and sample the subsets $S_1$ and $S_2$ according to distributions $Q$, where $Q$ are configured by option I and option II. To mitigate the impact of $b_H=40$, we set $\tau=2$ in both option I and option II. To avoid potential over-utility, we conduct $\gamma=0.8$, $\gamma_2=0.8$. Notably, we opt $\alpha=3, 6, 8$ to represent different hedge scenarios.
  	
\begin{figure*}[htbp]
  		\centering
  		\subfigure[australian]
  		{\includegraphics[width=0.327\textwidth]{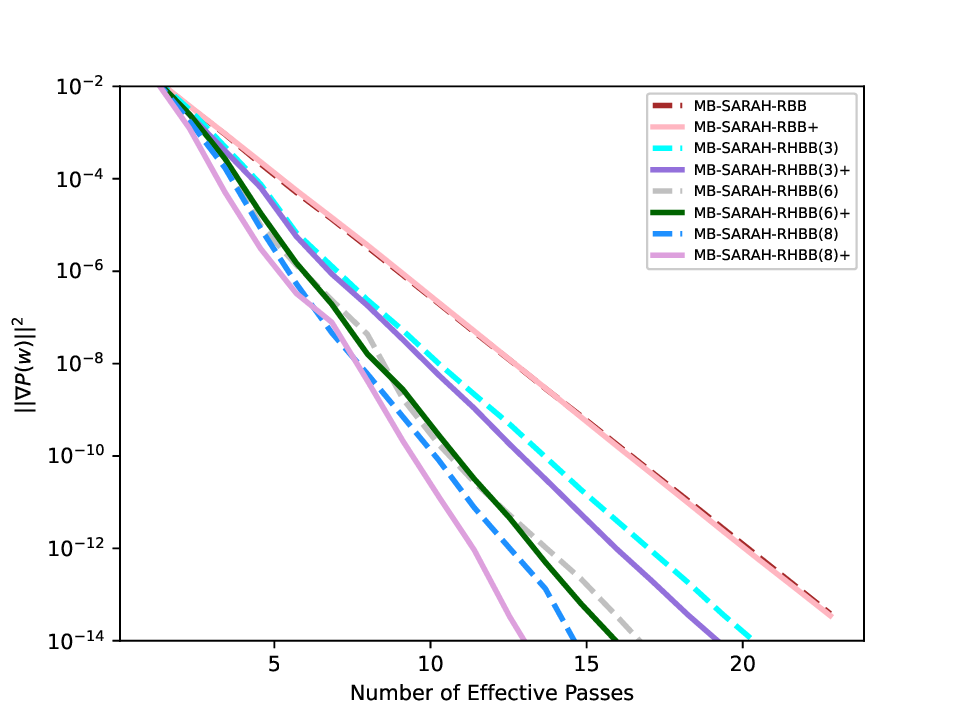}}
  		\subfigure[madelon]
  		{\includegraphics[width=0.327\textwidth]{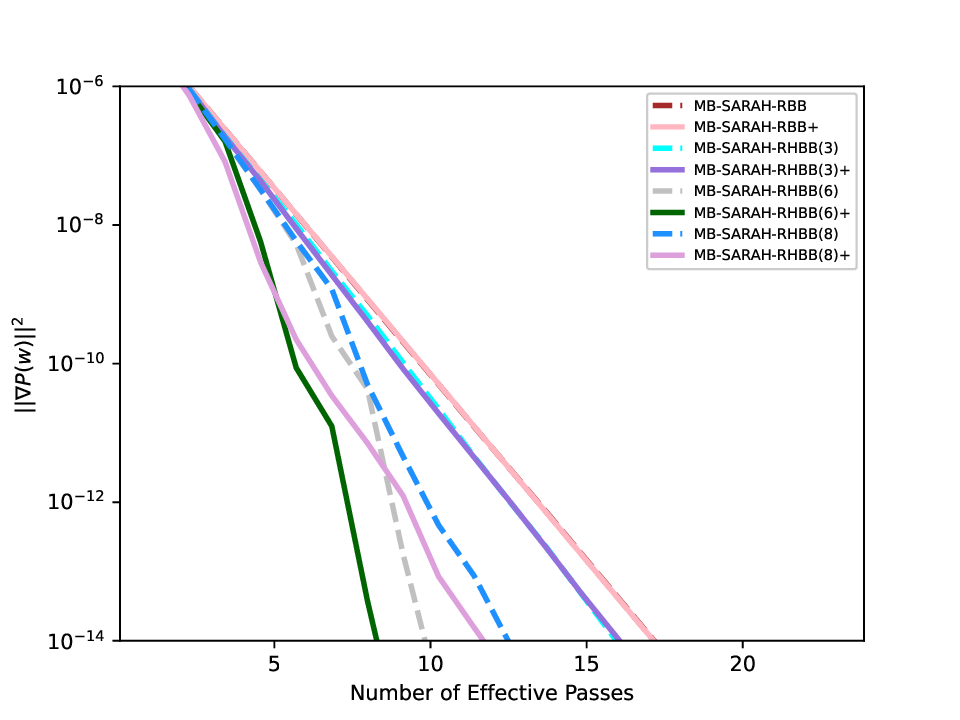}}
  		\subfigure[german]
  		{\includegraphics[width=0.327\textwidth]{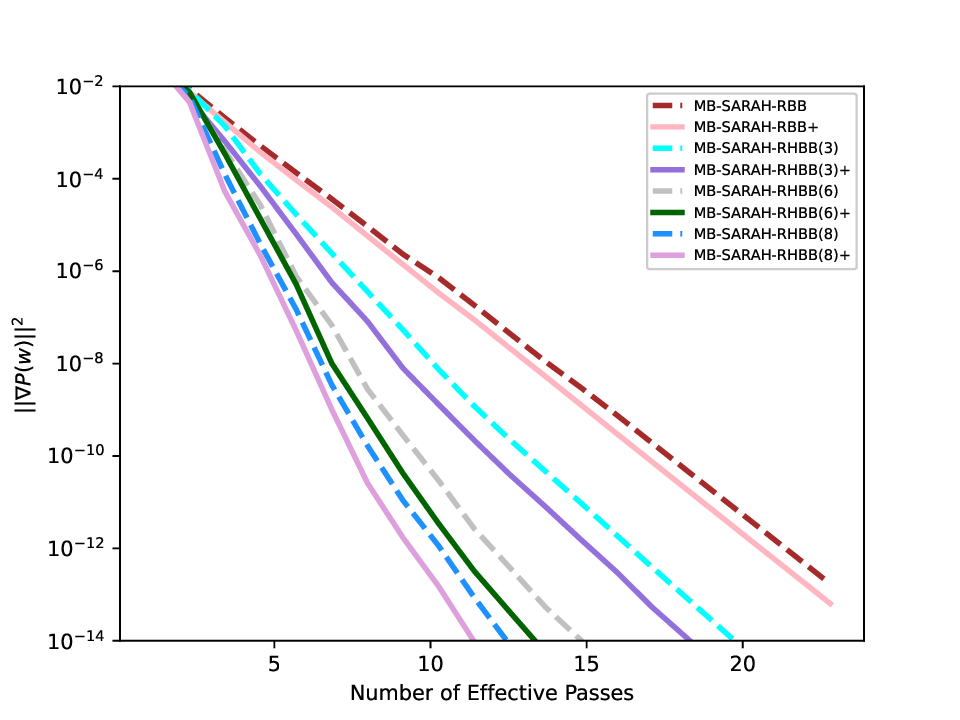}}
  		\caption{\footnotesize Comparisons of MB-SARAH-RHBB+ (solid lines) and MB-SARAH-RHBB (dash lines) under \textbf{option I}.}
  		\label{fig12}
  	\end{figure*}

\begin{figure*}[htbp]
  		\centering
  		\subfigure[australian]
  		{\includegraphics[width=0.327\textwidth]{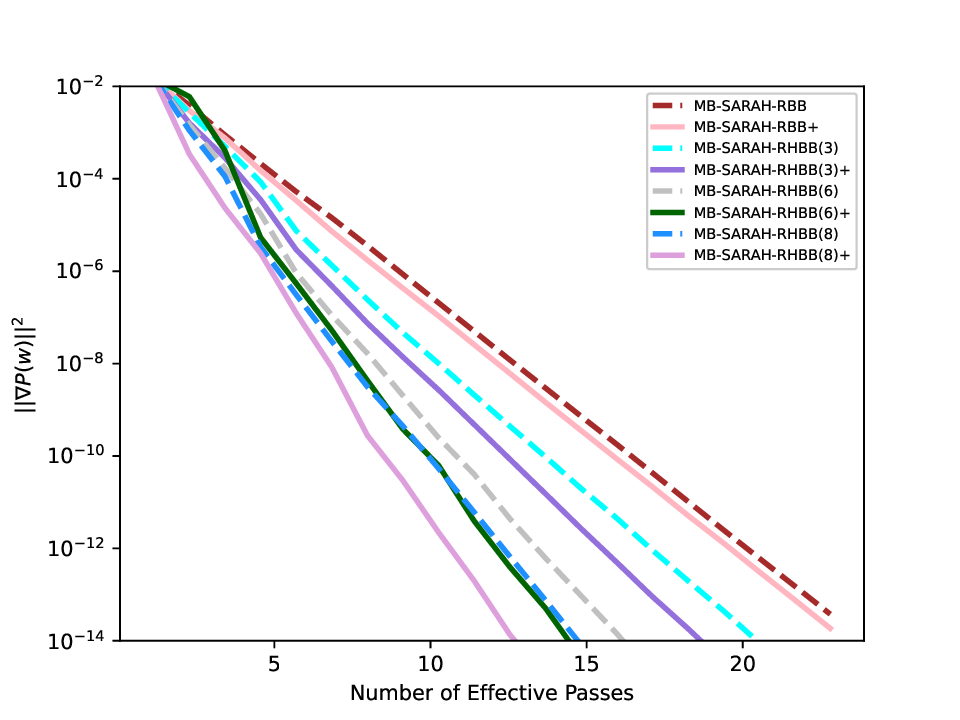}}
  		\subfigure[madelon]
  		{\includegraphics[width=0.327\textwidth]{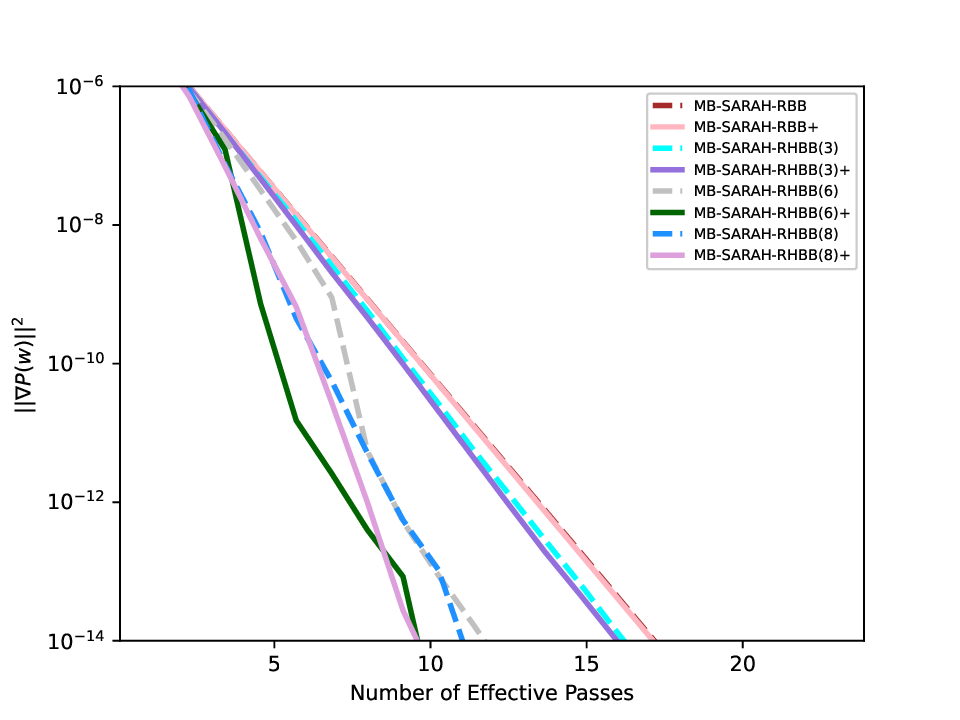}}
  		\subfigure[german]
  		{\includegraphics[width=0.327\textwidth]{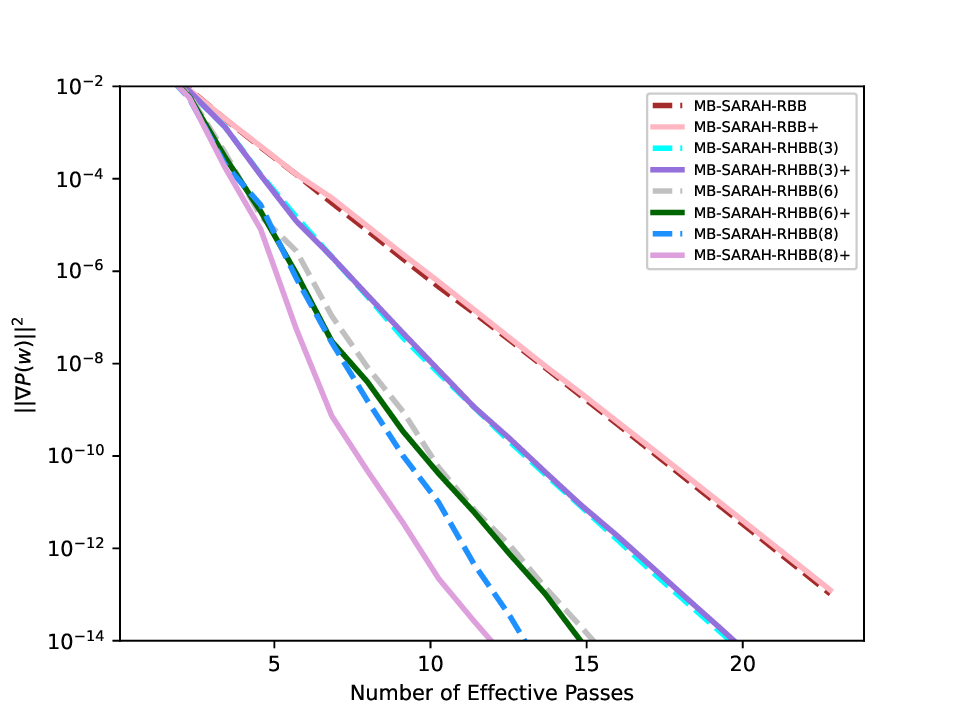}}
  		\caption{\footnotesize Comparisons of MB-SARAH-RHBB+ (solid lines) and MB-SARAH-RHBB (dash lines) under \textbf{option II}.}
  		\label{fig13}
  	\end{figure*}

\begin{figure*}[htbp]
  		\centering
  		\subfigure[australian]
  		{\includegraphics[width=0.327\textwidth]{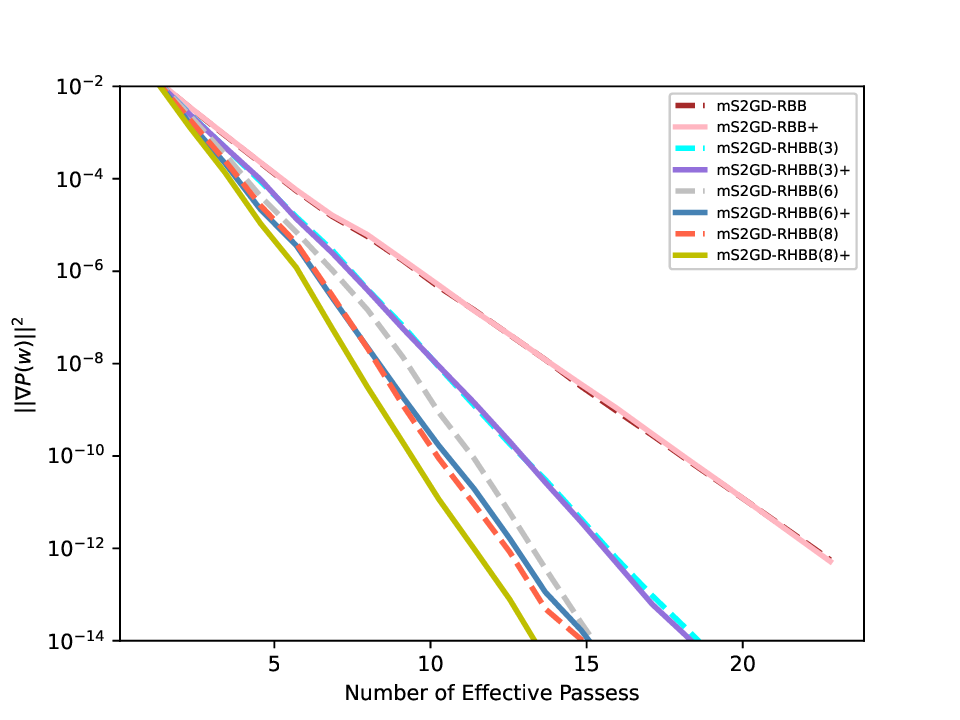}}
  		\subfigure[madelon]
  		{\includegraphics[width=0.327\textwidth]{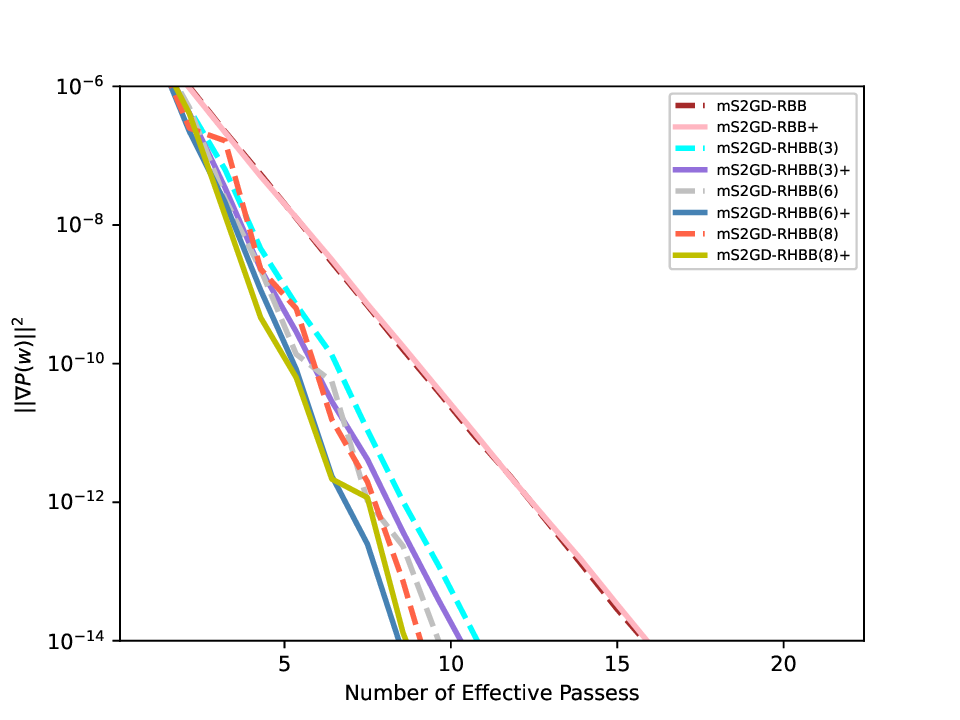}}
  		\subfigure[german.numer]
  		{\includegraphics[width=0.327\textwidth]{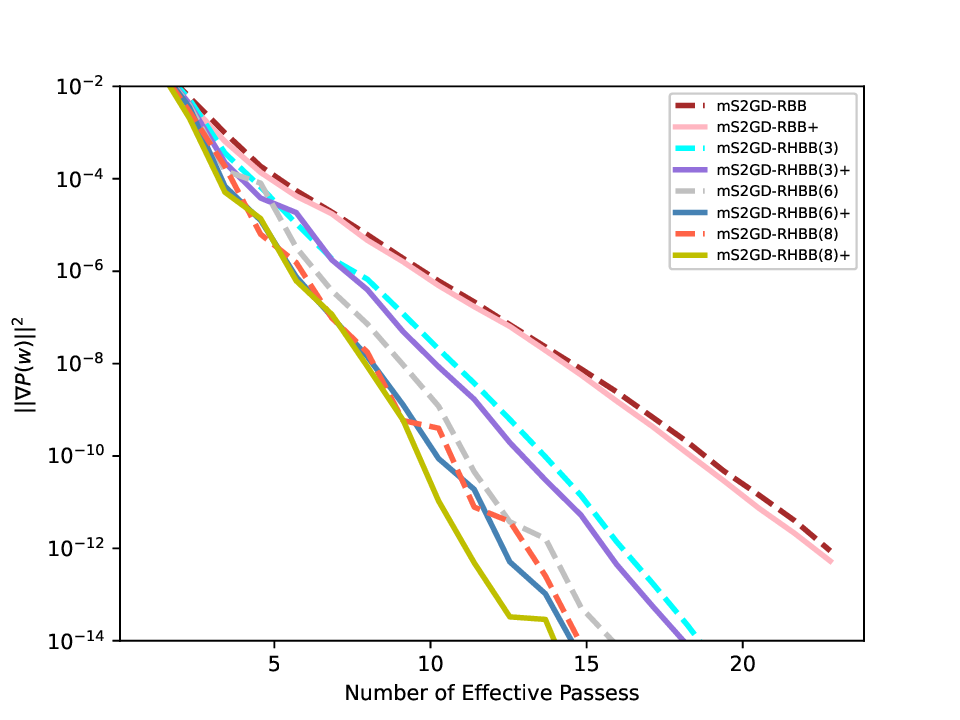}}
  		\caption{\footnotesize Comparisons of mS2GD-RHBB+ (solid lines) and mS2GD-RHBB (dash lines) under \textbf{option I}.}
  		\label{fig14}
  	\end{figure*}

\begin{figure*}[htbp]
  		\centering
  		\subfigure[australian]
  		{\includegraphics[width=0.327\textwidth]{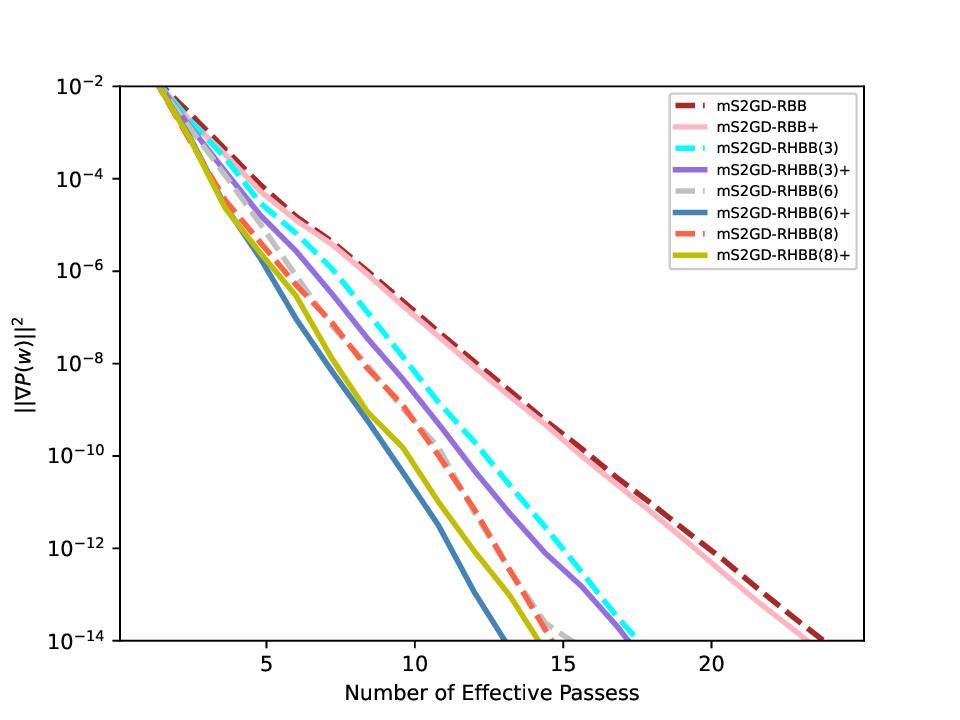}}
  		\subfigure[madelon]
  		{\includegraphics[width=0.327\textwidth]{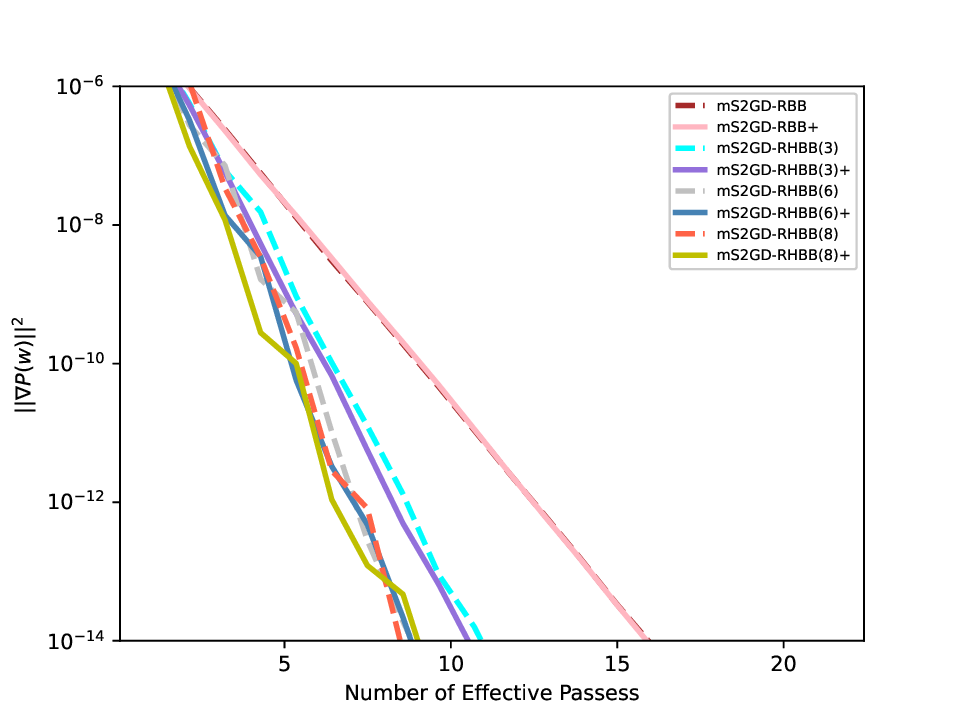}}
  		\subfigure[german.numer]
  		{\includegraphics[width=0.327\textwidth]{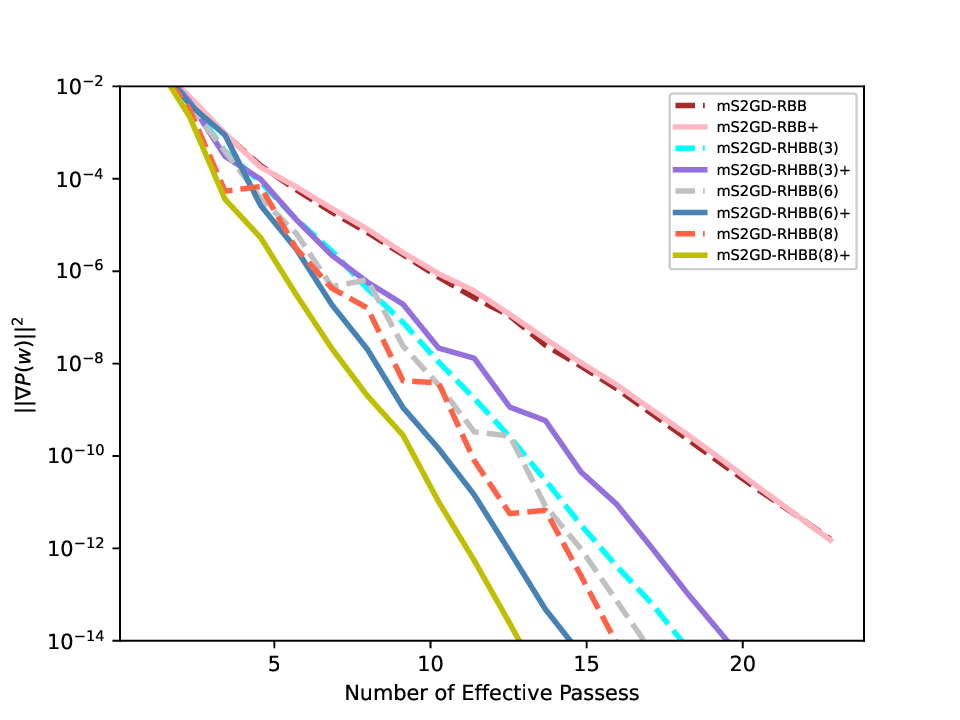}}
  		\caption{\footnotesize Comparisons of mS2GD-RHBB+ (solid lines) and mS2GD-RHBB (dash lines) under \textbf{option II}.}
  		\label{fig15}
  	\end{figure*}

  	Figs. \ref{fig12} - \ref{fig15} display the numerical results of MB-SARAH-RHBB+ and mS2GD-RHBB+. In the experiments, we select $\alpha=3, 6, 8$ to provide different hedge scenarios (slight, moderate and intense), under each of which we analyze the effect of the importance sampling technique in RHBB+.
  	
  	In Figs. \ref{fig12}, \ref{fig13} of MB-SARAH-RHBB+, when a conservative value of $\alpha=3$ is adopted (under the slight hedge scenario), both option I and option II seem to achieve limited improvement. However, under aggressive choices of $\alpha=6, 8$ (under the moderate and intense hedge scenarios), the importance sampling accelerates the convergence significantly.
  	
  	Figs. \ref{fig14}, \ref{fig15} show that the performance of mS2GD-RHBB+ has notable refinement when applying a large $\alpha=8$ (under the intense hedge scenario). Nonetheless, its numerical results on data set $madelon$ are not stable. Thereby, MB-SARAH-RHBB+ is more robust than mS2GD-RHBB+.
  	
  	In order to supply comprehensive illustrations, we further integrate the importance sampling into the original RBB rule and obtain the corresponding by-products of MB-SARAH-RBB+ and mS2GD-RBB+ algorithms. In Figs. \ref{fig12} - \ref{fig15}, we well include the comparisons between RBB and RBB+, verifying that the importance sampling can not yield improvement in the original RBB rule. Overall, the importance sampling in RHBB+ is more attuned to the large values of $\alpha$ (under the moderate and intense hedge scenarios).

  	\subsubsection{Comparison with other state-of-art methods}
  	\ 
  	\vskip10pt
  	\noindent  \textbf{Parametric Settings:} In MB-SARAH-RHBB and mS2GD-RHBB, we set $b=4$, the unified $b_H=40$ and sample the subsets $S$, $S_1$, $S_2$ according to uniform distribution. We employ the same $\gamma$ and $\gamma_2$ as in previous experiments to ensure consistent illustrations, which are $\gamma=1$, $\gamma_2=1$. In addition, we fine-tune each of the other algorithms as follow: \\
  	(1) \textbf{SVRG:} Accelerating stochastic gradient descent using predictive variance reduction \cite{johnson}. We employ the best-tuned constant step size in SVRG setting.\\
  	(2) \textbf{SVRG-BB:} Stochastic variance reduced algorithm (SVRG) with Barzilai and Borwein step size \cite{tan}.\\
  	(3) \textbf{mS2GD-BB:} Semi-stochastic algorithm (mS2GD) with Barzilai and Borwein step size \cite{yang3}, a batch version of SVRG-BB.\\
  	(4) \textbf{Acc-Prox-SVRG:} Accelerating variance reduced algorithm with a momentum Nesterov's structure \cite{nitanda}. We set $\eta=1$, $\delta = 1$, $b=100$, $m=\delta b$ and $\beta_k=\frac{b-2}{b+2}$ as suggested in \cite{nitanda}.\\
  	(5) \textbf{Acc-Prox-SVRG-BB:} Acc-Prox-SVRG stochastic algorithm with Barzilai and Borwein step size \cite{yang4}. We set the related parameters according to \cite{yang4}.\\
  	(6) \textbf{SARAH+:} An implementation version of SARAH \cite{nguyen}. Best hand-tuned constant step size was employed in the optimization process. \\
  	(7) \textbf{SVRG-ABB:} SVRG stochastic algorithm with adaptive Barzilai and Borwein step size. The adaptive parameter $k=0.5$ is set for the robustness.
  	
\begin{figure*}[htbp]
  		\centering
  		\subfigure[a8a]
  		{\includegraphics[width=0.45\textwidth]{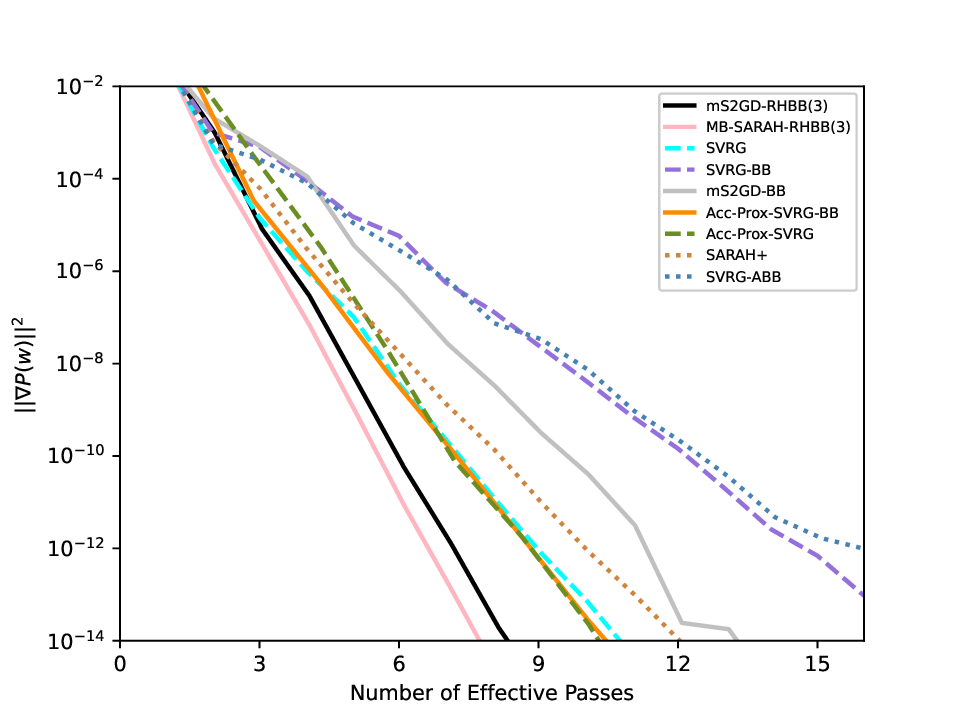}}
  		\subfigure[w8a]
  		{\includegraphics[width=0.45\textwidth]{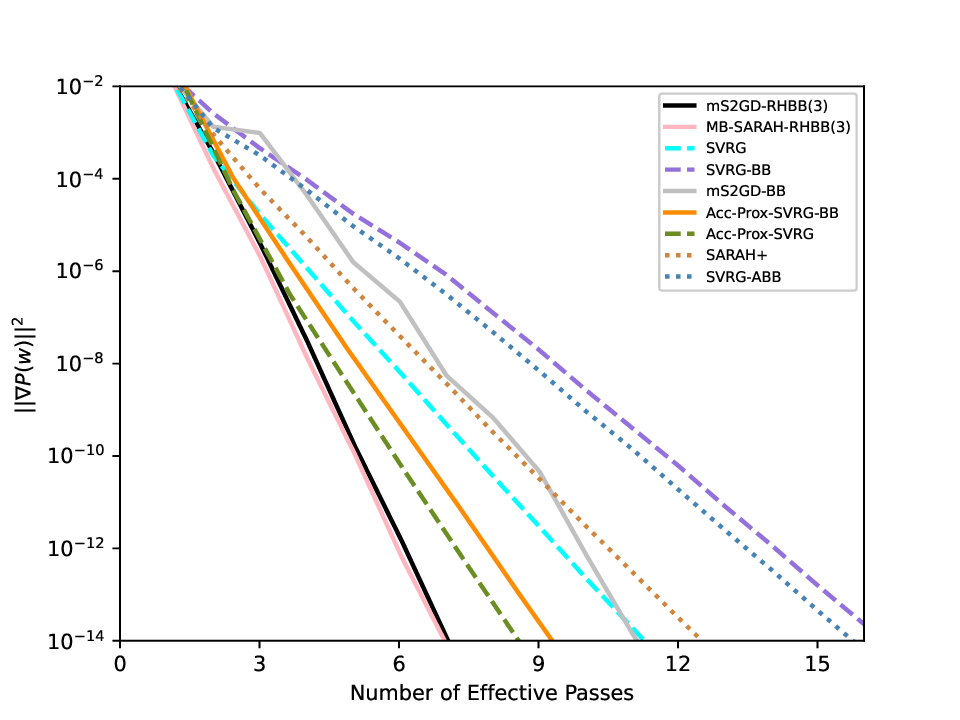}}
  		
  		\subfigure[ijcnn1]
  		{\includegraphics[width=0.45\textwidth]{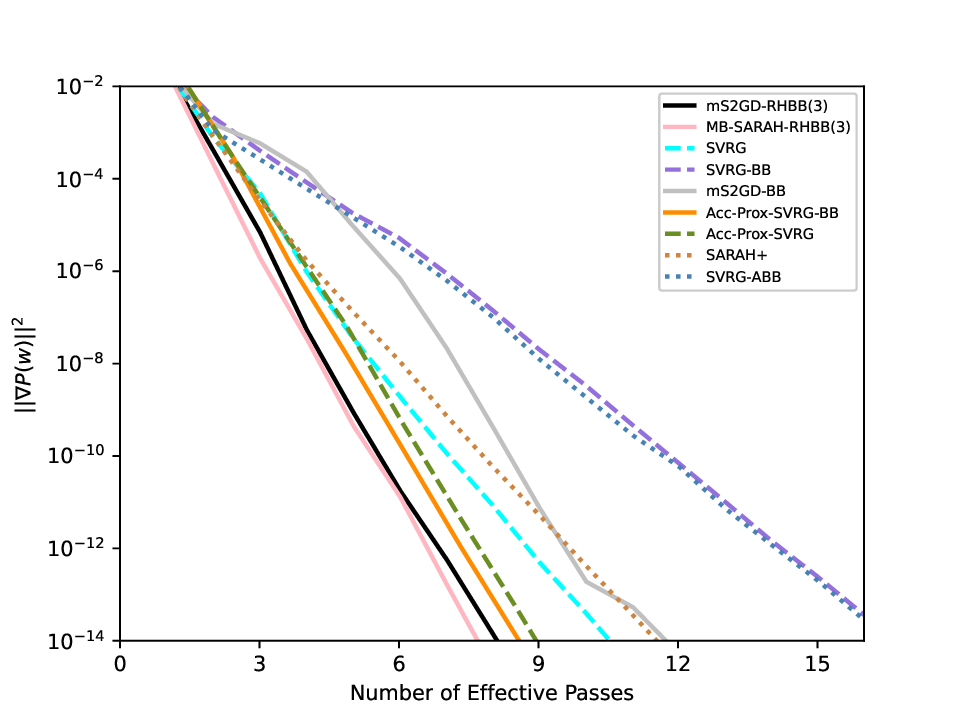}}
  		\subfigure[covtype]
  		{\includegraphics[width=0.45\textwidth]{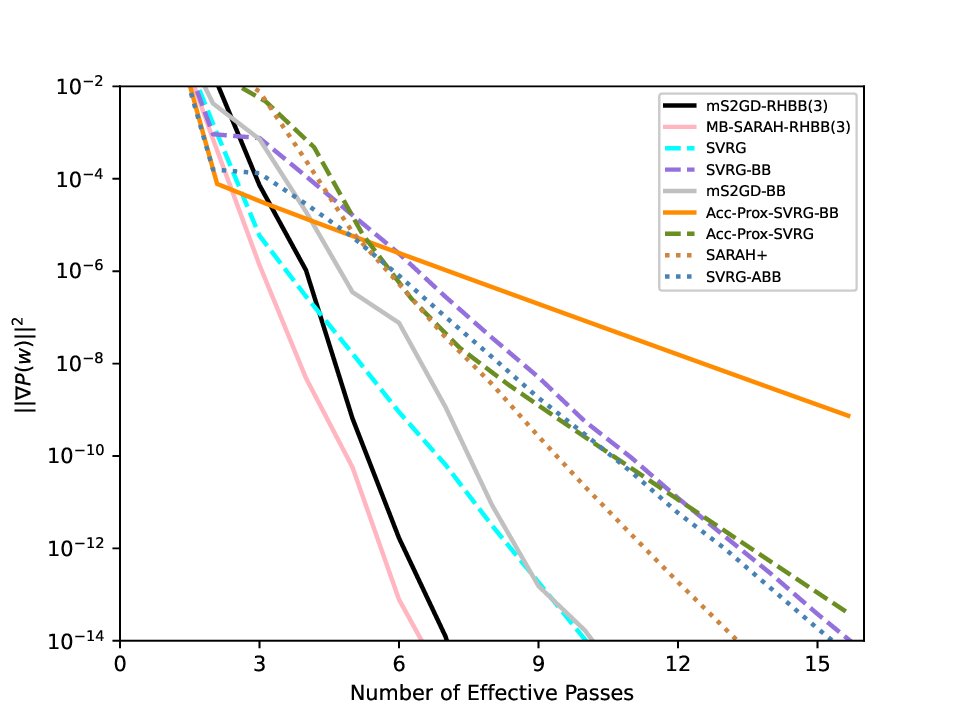}}
  		
  		\subfigure[phishing]
  		{\includegraphics[width=0.45\textwidth]{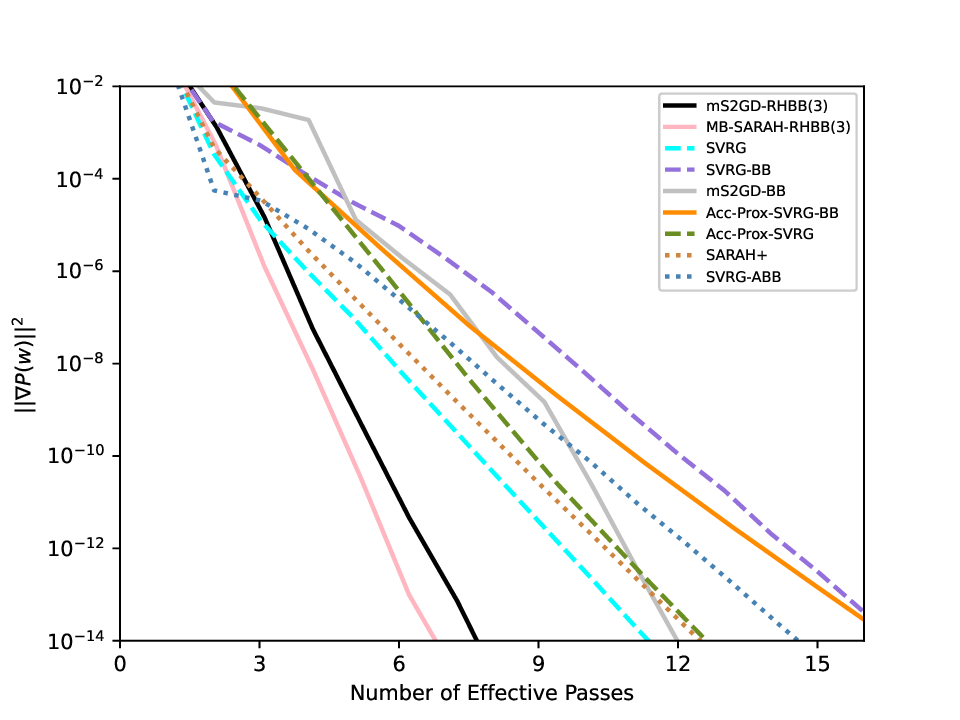}}
  		\subfigure[mushrooms]
  		{\includegraphics[width=0.45\textwidth]{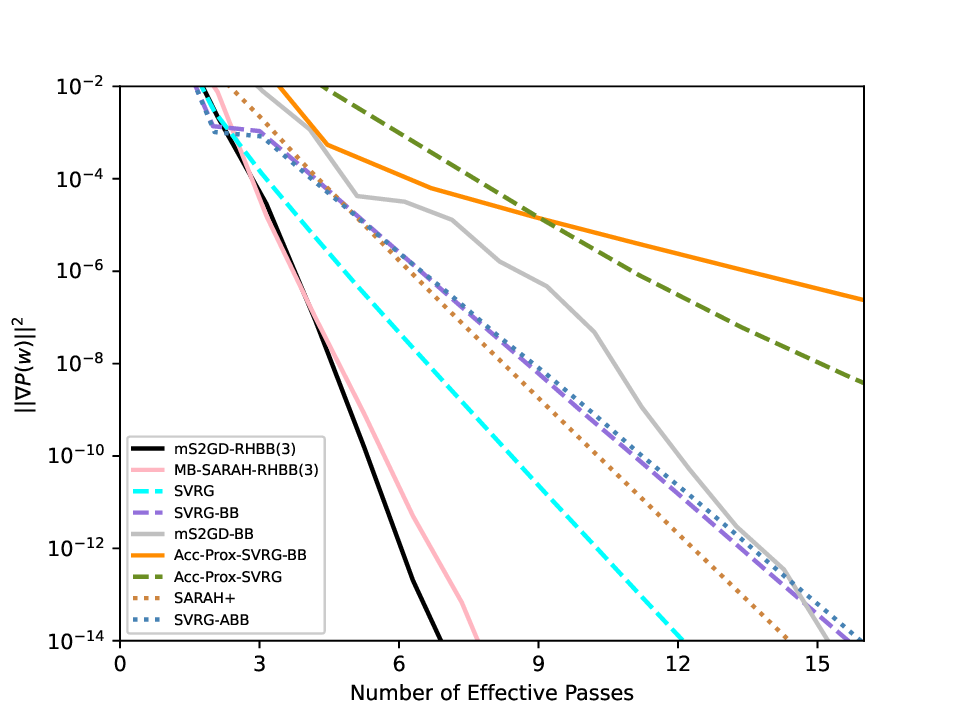}}
  		\caption{\footnotesize Comparisons of different algorithms.}
  		\label{fig16}
  	\end{figure*}

  	It's noted from Fig. \ref{fig16} that MB-SARAH-RHBB and mS2GD-RHBB, with the unvaried $\alpha=3$, outperform the other state-of-art methods consistently on all the six data sets. Certain algorithms may be competitive on $ijcnn1$, but soon expose their powerlessness on the others. 
  	
  	Referring back to Figs. \ref{fig2} - \ref{fig5}, there indeed exists an optimal bound on the hedge magnitude. We argue that further improvement can be achieved by controlling the hedge magnitude dynamically, instead of fixing $h(\cdot)=1$.
  	
  	\subsection{Experiment investigating for Adaptive Hedge Effect}
  	Now, we reveal the additional efficiency of our iterative adaptor $h(\cdot)$.
  	
  	Figs. \ref{fig2}, \ref{fig3} show that a properly-tuned $\alpha$ greatly accelerates the convergence, but the subsequent Figs. \ref{fig4}, \ref{fig5} show that the effective magnitude can be corrupted due to over-hedging or excessive enlargement. In view of it, we propose decreasing the value of $h(\cdot)$ along the updates: on one side to positively boost the convergence in early periods, on the other for avoiding over-aggressive steps around the global optimum. For conciseness, we devise the increment to be inversely proportional to the iterative indicators $\left(\sigma_1s+\sigma_2k\right)$, i.e.,  $h(\sigma_1s+\sigma_2k)=\frac{1+\sigma_1s+\sigma_2k}{\sigma_1s+\sigma_2k}$. Equally, we conduct experiments under $b_1=b_2=b_H$ first.
  	
  	In this subsection, we use notations with the suffix `pure' to represent algorithms under $h(\sigma_1s+\sigma_2k)=1$. Notations are summarized in \textbf{Table \ref{table3}}.
  	
  	\begin{table}
  		\setlength{\abovecaptionskip}{0pt}
  		\setlength{\belowcaptionskip}{1pt}
  		\caption{NOTATIONS DESCRIPTIONS}
  		\centering
  		\vspace{4pt}
  		\begin{tabular}{lclc}
  			\hline
  			\multicolumn{1}{l}{Notations}&\multicolumn{1}{l}{\qquad Hedge Bases} &\multicolumn{1}{l}{\qquad Step Sizes} &\multicolumn{1}{c}{\qquad Adaptivity}
  			\\ \hline\noalign{\smallskip}
  			MB-SARAH-RBB  & \qquad \ding{56} & \qquad RBB & \ding{52}\\
  			MB-SARAH-RBB+ & \qquad \ding{56} & \qquad RBB+ & \ding{52}\\
  			MB-SARAH-RHBB($\alpha$) & \qquad $\alpha$ & \qquad RHBB & \ding{52}\\
  			MB-SARAH-RHBB($\alpha$)+ & \qquad $\alpha$ & \qquad RHBB+ & \ding{52}\\
  			MB-SARAH-RHBB($\alpha$) - pure & \qquad $\alpha$ & \qquad RHBB & \ding{56}\\
  			MB-SARAH-RHBB($\alpha$)+ - pure & \qquad $\alpha$ & \qquad RHBB+ & \ding{56}\\
  			mS2GD-RBB  & \qquad \ding{56} & \qquad RBB & \ding{52}\\
  			mS2GD-RBB+ & \qquad \ding{56} & \qquad RBB+ & \ding{52}\\
  			mS2GD-RHBB($\alpha$) & \qquad $\alpha$ & \qquad RHBB & \ding{52}\\
  			mS2GD-RHBB($\alpha$)+ & \qquad $\alpha$ & \qquad RHBB+ & \ding{52}\\
  			mS2GD-RHBB($\alpha$) - pure & \qquad $\alpha$ & \qquad RHBB & \ding{56}\\
  			mS2GD-RHBB($\alpha$)+ - pure & \qquad $\alpha$ & \qquad RHBB+ & \ding{56}\\
  			\hline
  		\end{tabular}
  		\label{table3}
  	\end{table}

  	\subsubsection{Adaptive MB-SARAH-RHBB/mS2GD-RHBB}
  	\ 
  	\vskip10pt
  	\noindent  \textbf{Parametric Settings:} We set $b=4$, the unified $b_H=40$ and sample the subsets $S$, $S_1$, $S_2$ according to uniform distribution. We perform an extensive search for the adaptive pair $(\sigma_1, \sigma_2)$ with three different settings: $(0.6, 0.2), (0.7, 0.1), (0.4, 0.4)$ (we ensure $\sigma_1+\sigma_2=0.8<1$ to allow several quadratic accelerations in the early iterations.). To ensure the comparability across aspects, the hedge base $\alpha$ is opted within $\{2, 3, 4, 5\}$. Following guidelines from \cite{yang}, we choose $\gamma=1$. By considering a moderate trade-off in mS2GD-RHBB, we implement $\gamma_2=1$.
  	
\begin{figure*}[htbp]
  		\centering
  		\subfigure[a8a]
  		{\includegraphics[width=0.327\textwidth]{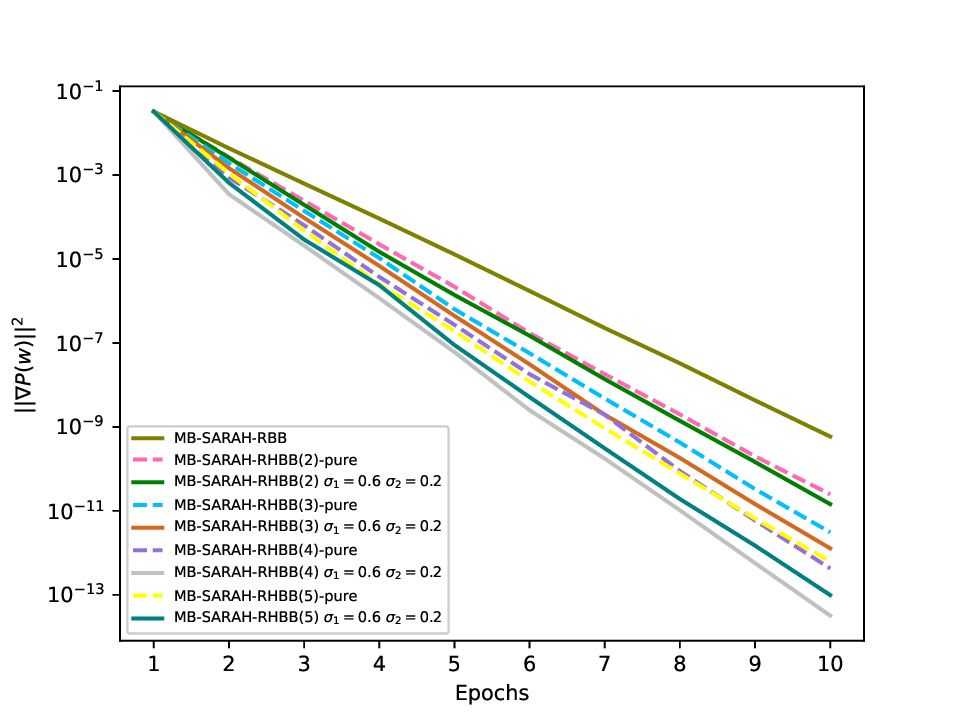}}
  		\subfigure[w8a]
  		{\includegraphics[width=0.327\textwidth]{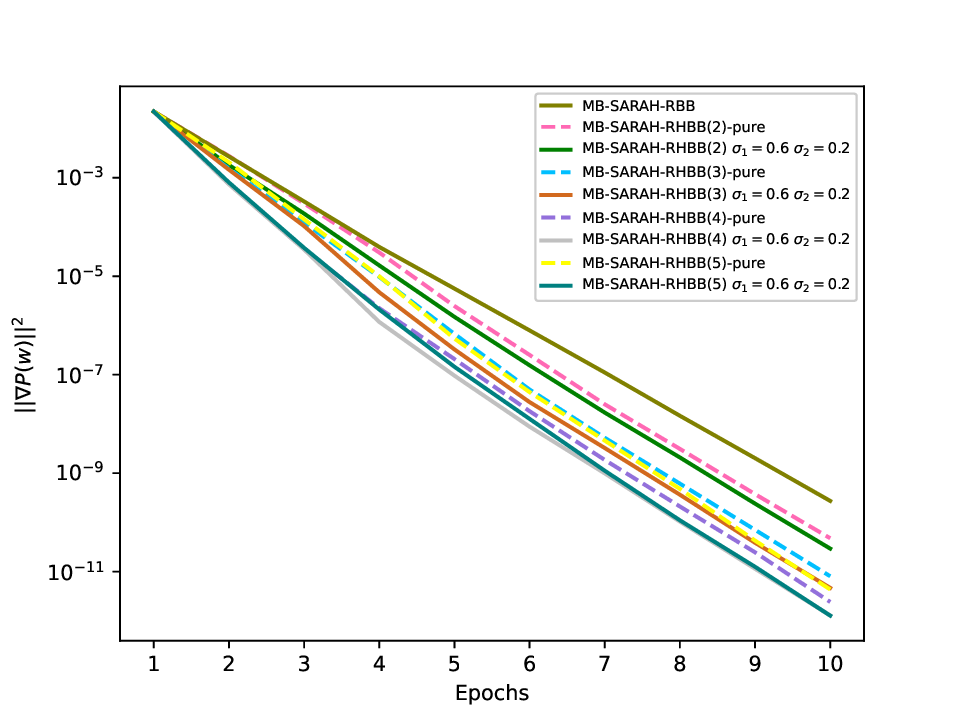}}
  		\subfigure[ijcnn1]
  		{\includegraphics[width=0.327\textwidth]{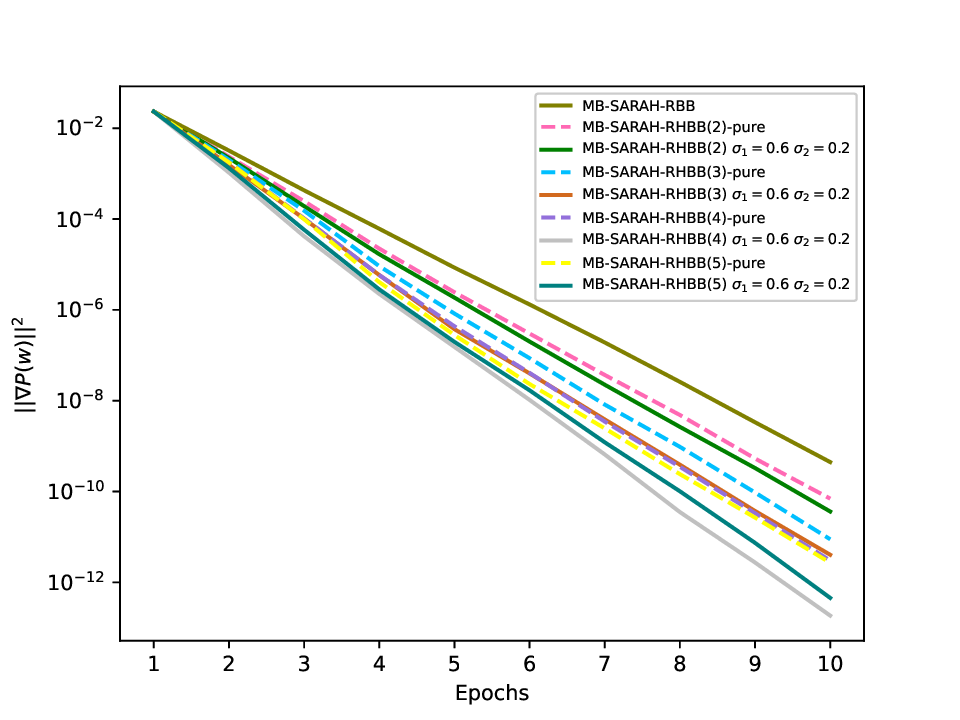}}
  		\caption{\footnotesize The performance of MB-SARAH-RBB, non-adaptive MB-SARAH-RHBB, and adaptive MB-SARAH-RHBB under $\sigma_1=0.6$, $\sigma_2=0.2$.}
  		\label{fig17}
  	\end{figure*}

\begin{figure*}[htbp]
  		\centering
  		\subfigure[a8a]
  		{\includegraphics[width=0.327\textwidth]{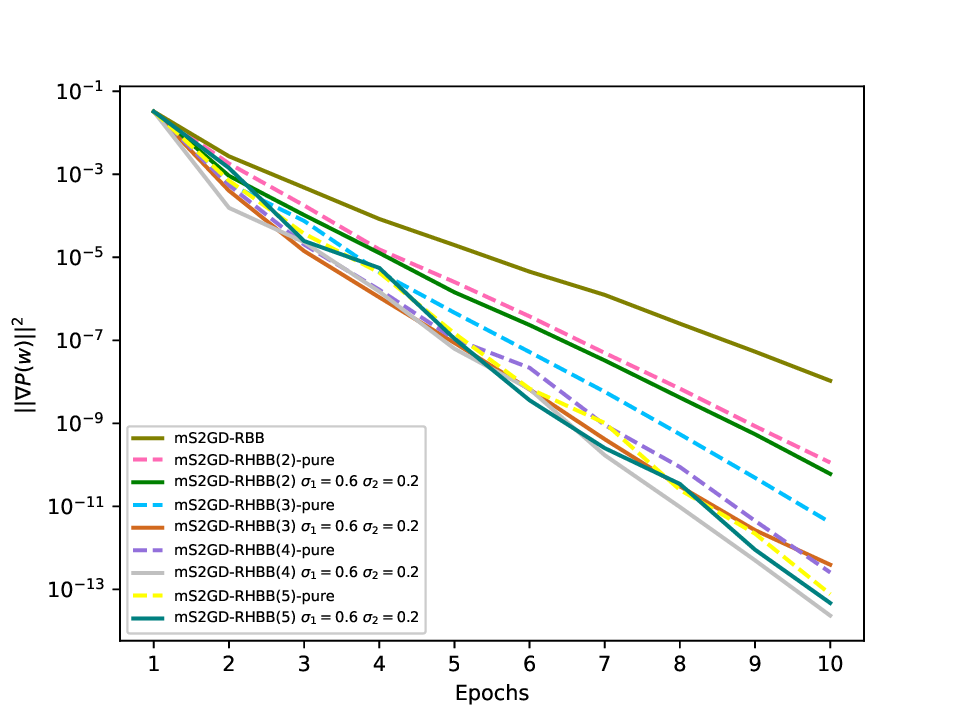}}
  		\subfigure[w8a]
  		{\includegraphics[width=0.327\textwidth]{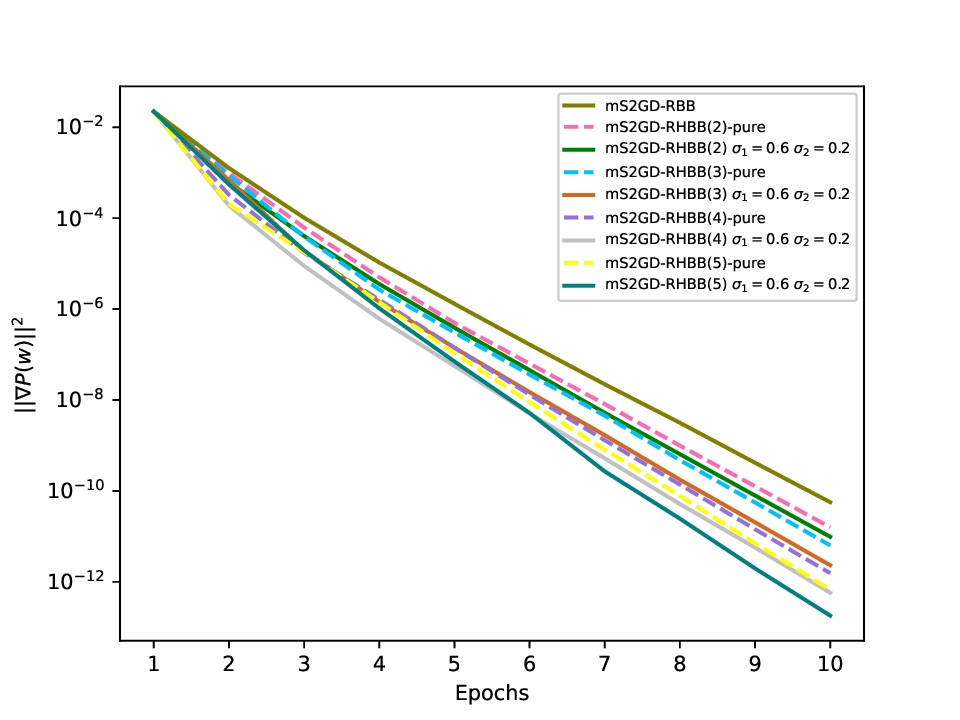}}
  		\subfigure[ijcnn1]
  		{\includegraphics[width=0.327\textwidth]{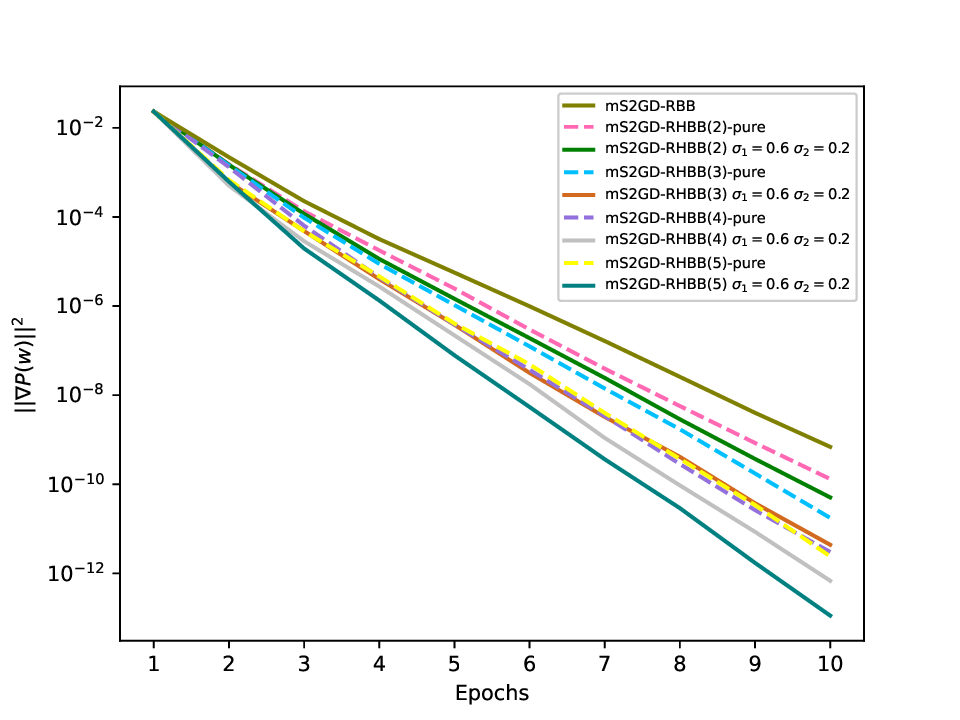}}
  		\caption{\footnotesize The performance of mS2GD-RBB, non-adaptive mS2GD-RHBB, and adaptive mS2GD-RHBB under $\sigma_1=0.6$, $\sigma_2=0.2$.}
  		\label{fig18}
  	\end{figure*}

\begin{figure*}[htbp]
  	\centering
  	\subfigure[a8a]
  	{\includegraphics[width=0.327\textwidth]{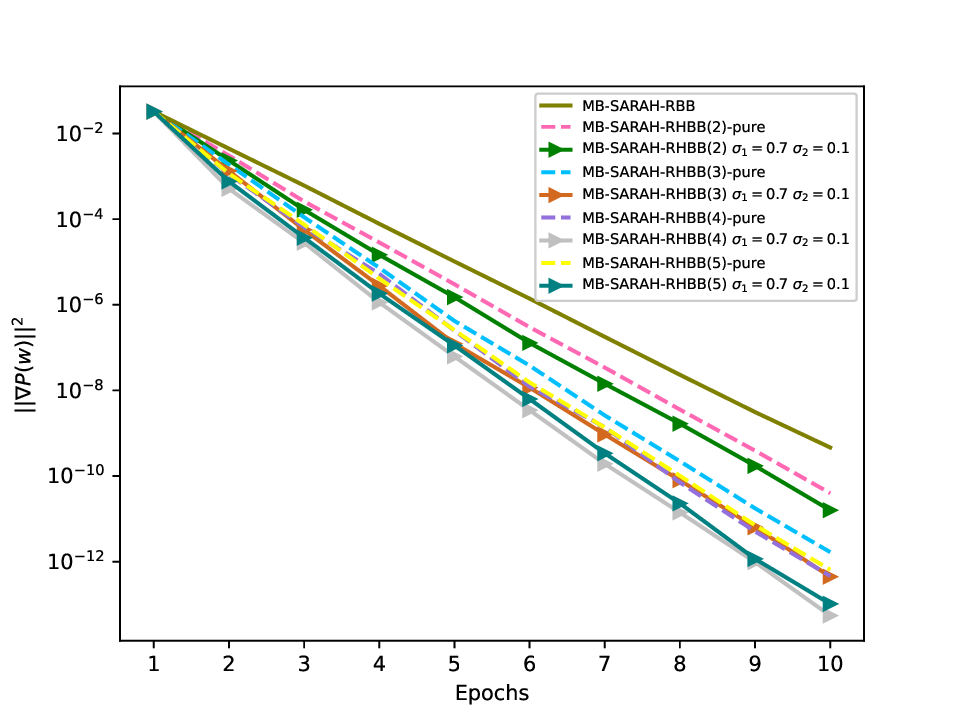}}
  	\subfigure[w8a]
  	{\includegraphics[width=0.327\textwidth]{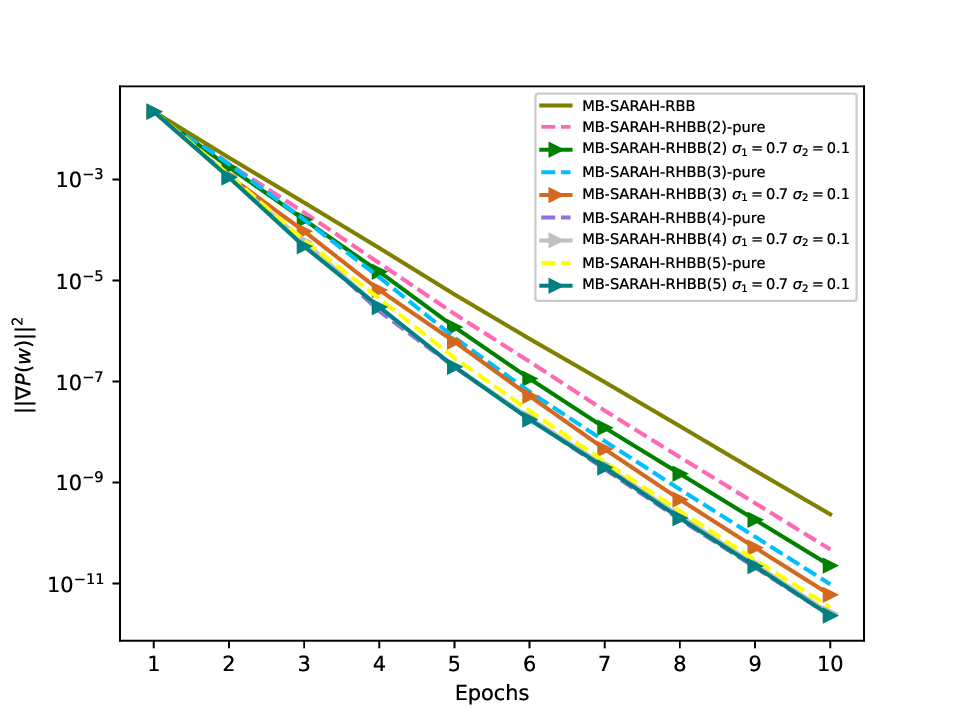}}
  	\subfigure[ijcnn1]
  	{\includegraphics[width=0.327\textwidth]{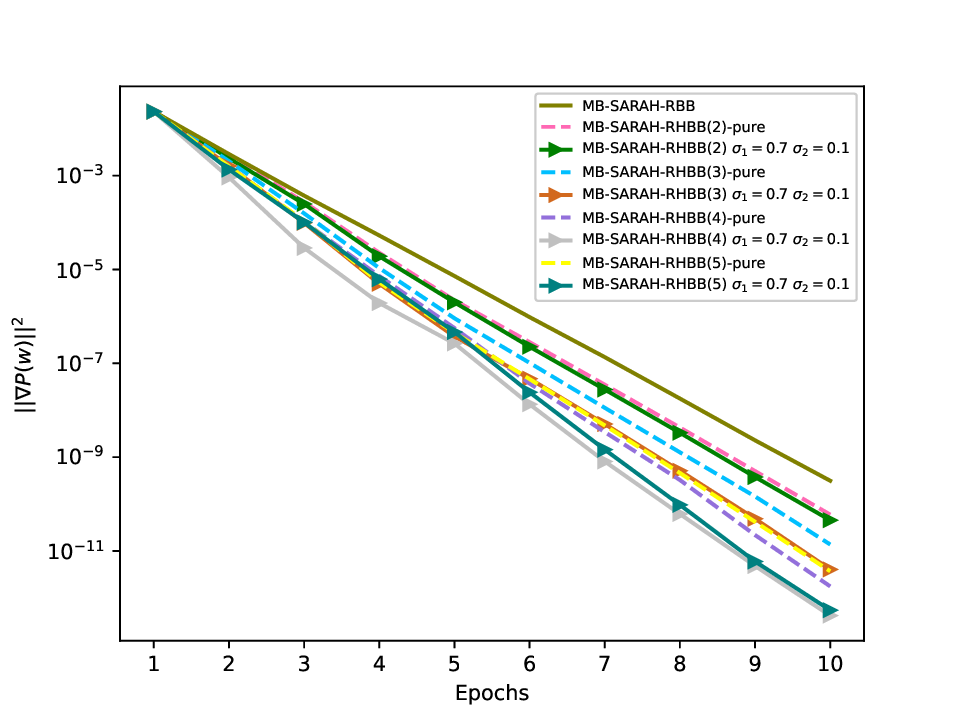}}
  	\caption{\footnotesize The performance of MB-SARAH-RBB, non-adaptive MB-SARAH-RHBB, and adaptive MB-SARAH-RHBB under $\sigma_1=0.7$, $\sigma_2=0.1$.}
  	\label{fig19}
  \end{figure*}

\begin{figure*}[htbp]
  	\centering
  	\subfigure[a8a]
  	{\includegraphics[width=0.327\textwidth]{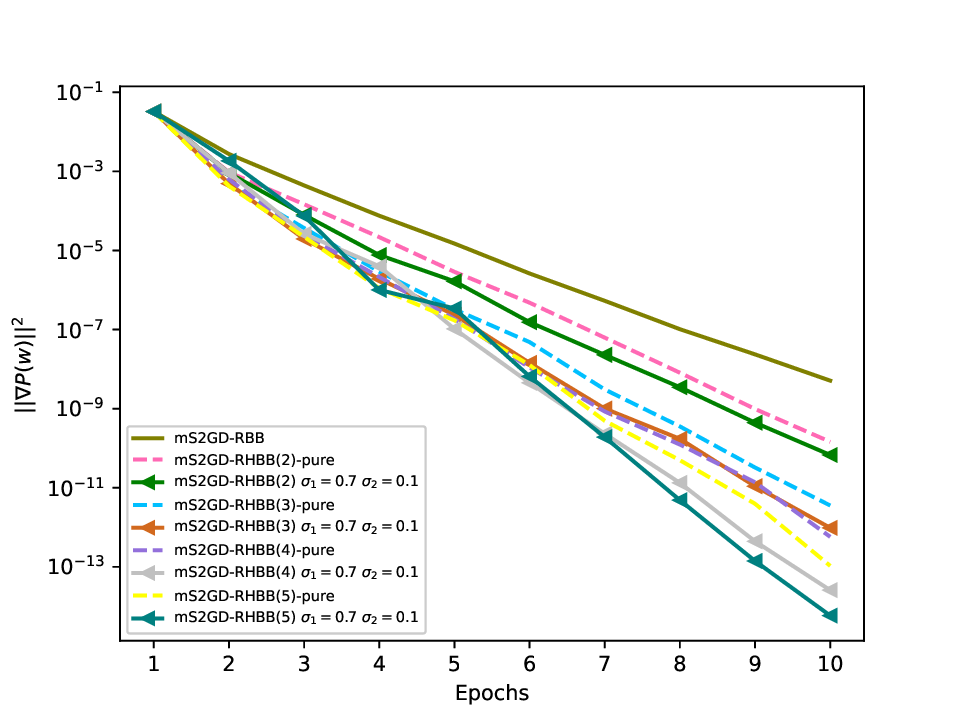}}
  	\subfigure[w8a]
  	{\includegraphics[width=0.327\textwidth]{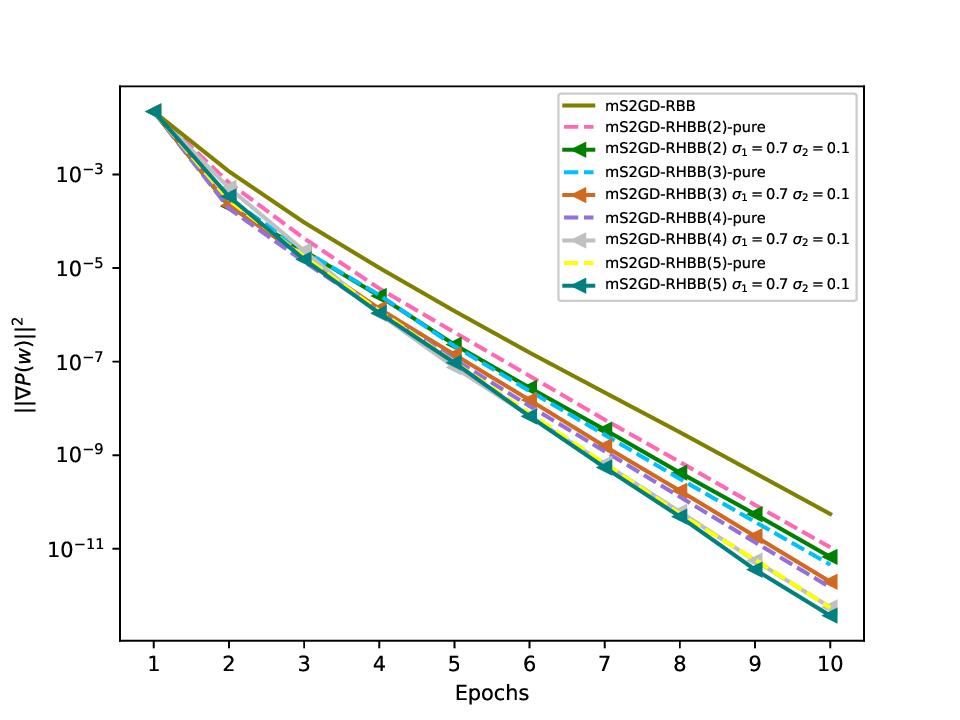}}
  	\subfigure[ijcnn1]
  	{\includegraphics[width=0.327\textwidth]{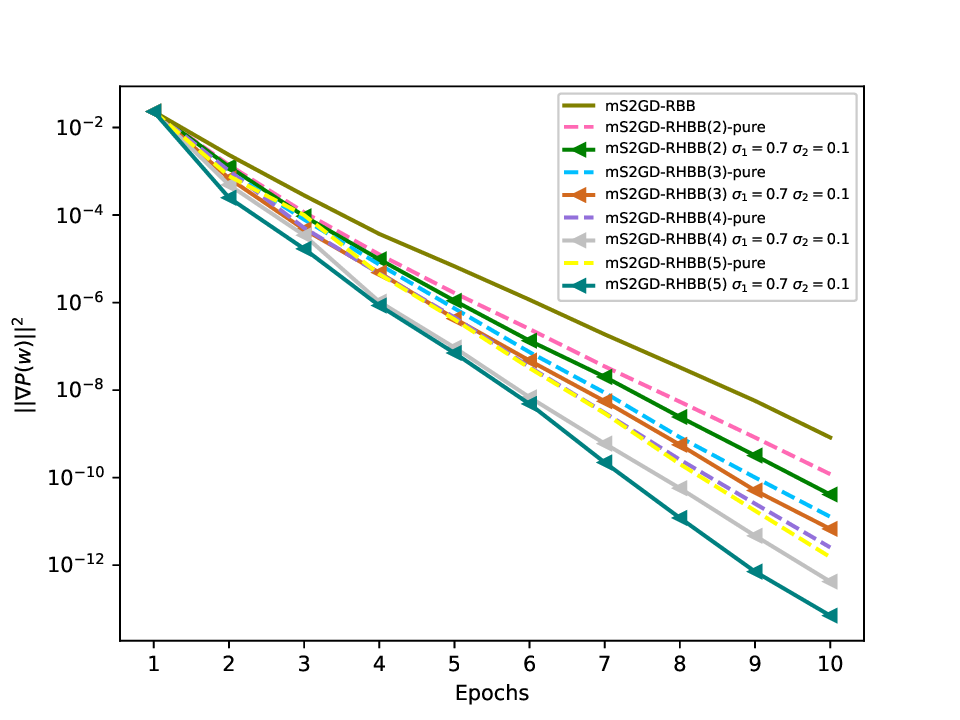}}
  	\caption{\footnotesize The performance of mS2GD-RBB, non-adaptive mS2GD-RHBB, and adaptive mS2GD-RHBB under $\sigma_1=0.7$, $\sigma_2=0.1$.}
  	\label{fig20}
  \end{figure*}

\begin{figure*}[htbp]
  		\centering
  		\subfigure[a8a]
  		{\includegraphics[width=0.327\textwidth]{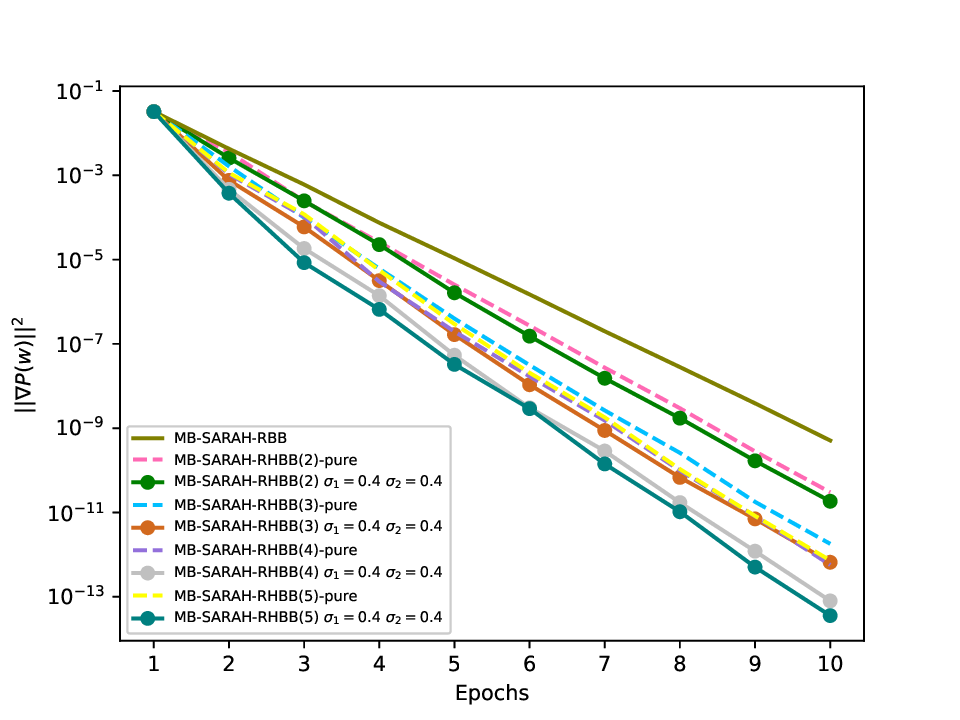}}
  		\subfigure[w8a]
  		{\includegraphics[width=0.327\textwidth]{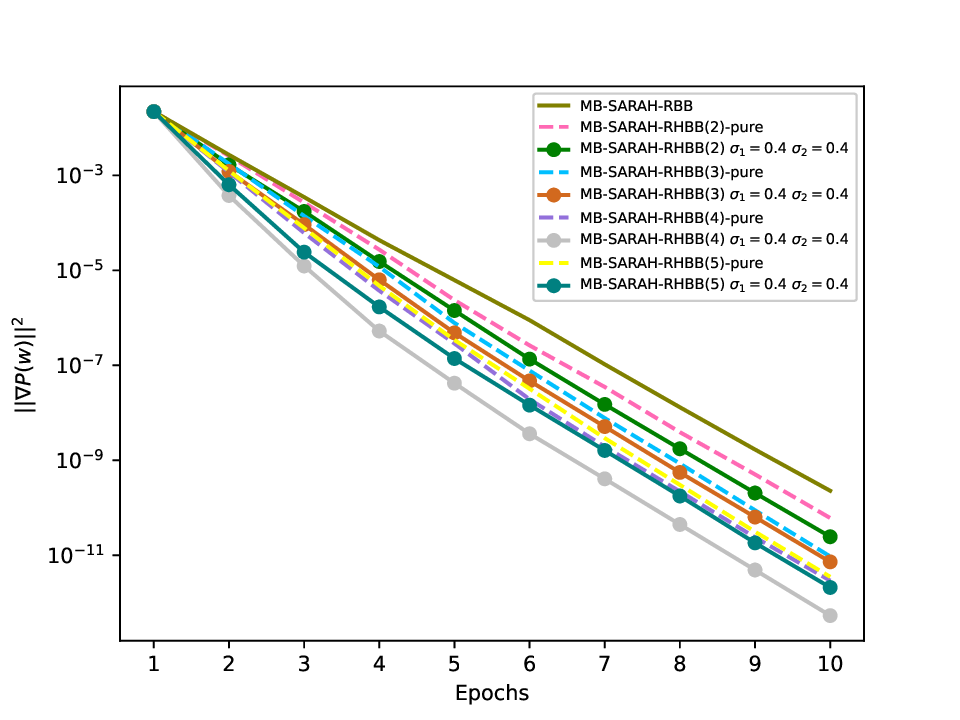}}
  		\subfigure[ijcnn1]
  		{\includegraphics[width=0.327\textwidth]{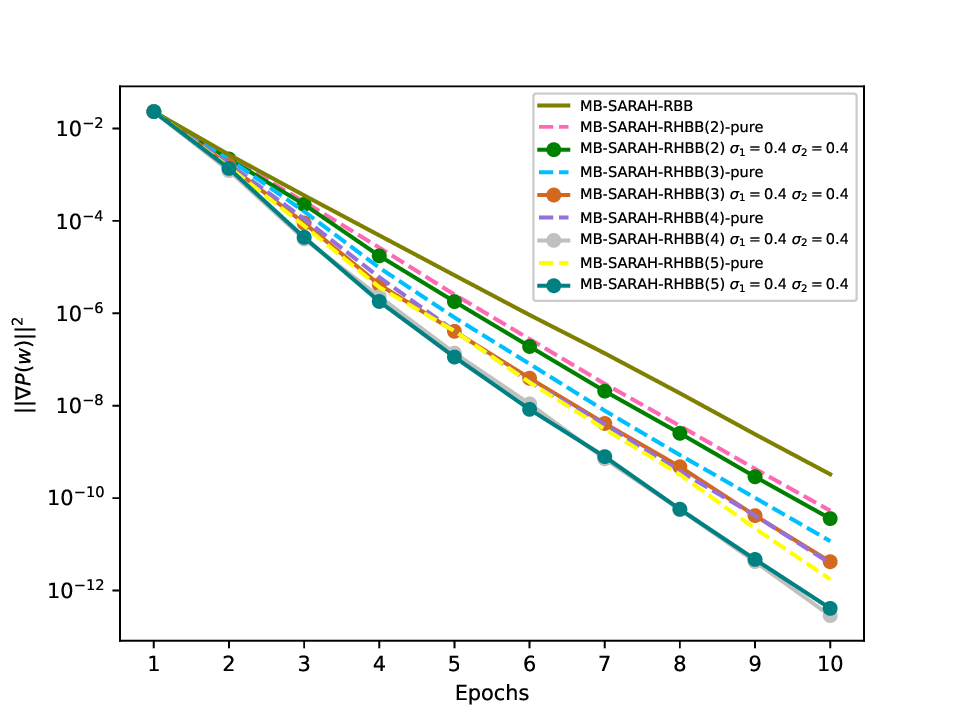}}
  		\caption{\footnotesize The performance of MB-SARAH-RBB, non-adaptive MB-SARAH-RHBB, and adaptive MB-SARAH-RHBB under $\sigma_1=0.4$, $\sigma_2=0.4$.}
  		\label{fig21}
  	\end{figure*}

\begin{figure*}[htbp]
  		\centering
  		\subfigure[a8a]
  		{\includegraphics[width=0.327\textwidth]{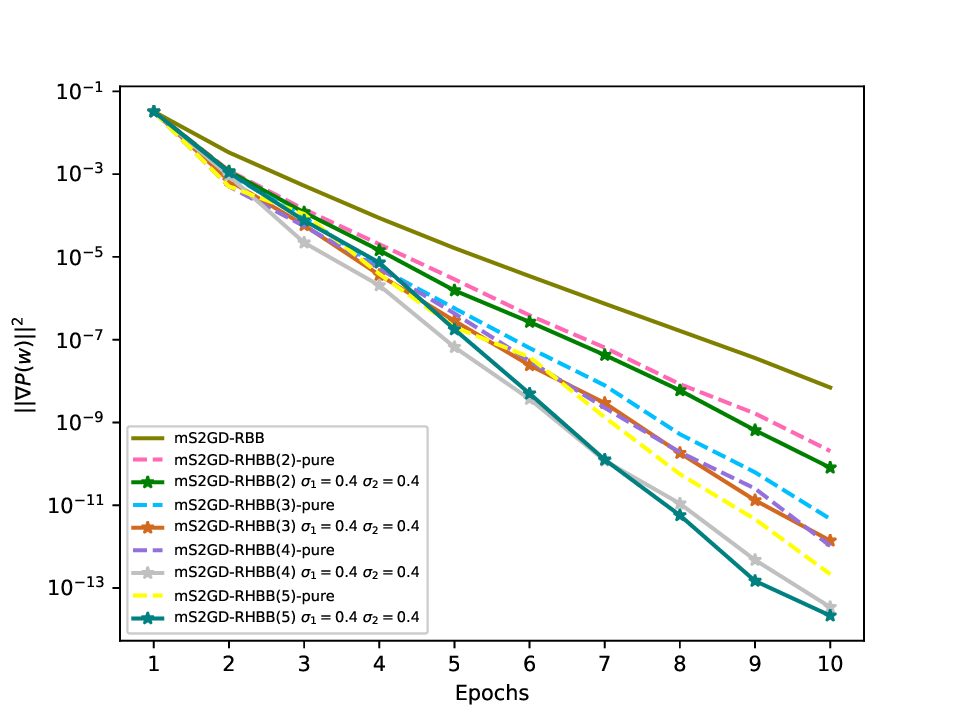}}
  		\subfigure[w8a]
  		{\includegraphics[width=0.327\textwidth]{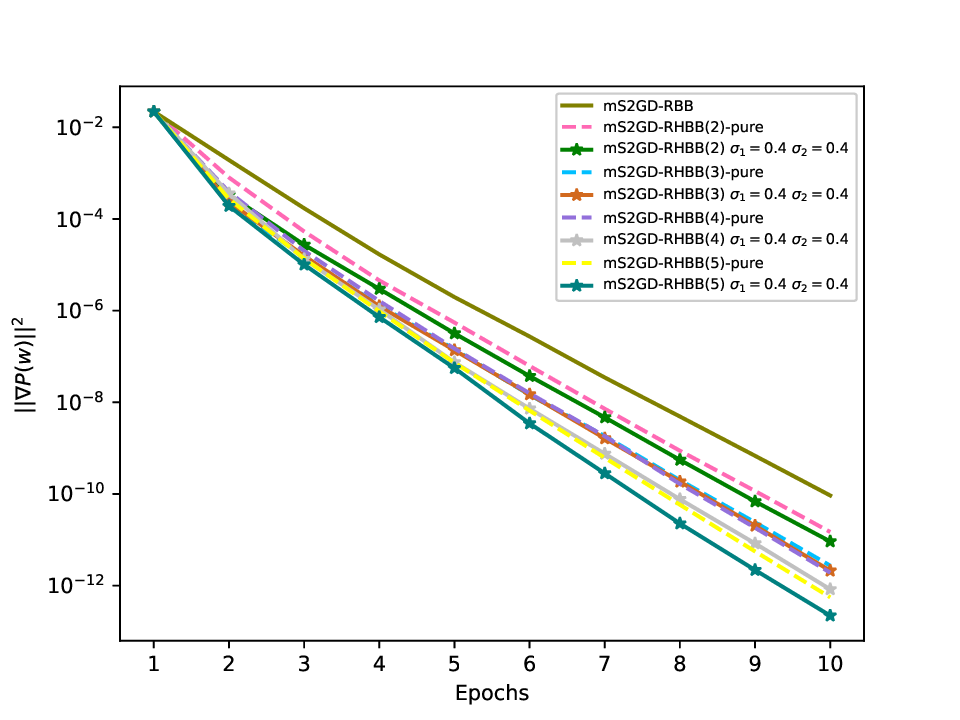}}
  		\subfigure[ijcnn1]
  		{\includegraphics[width=0.327\textwidth]{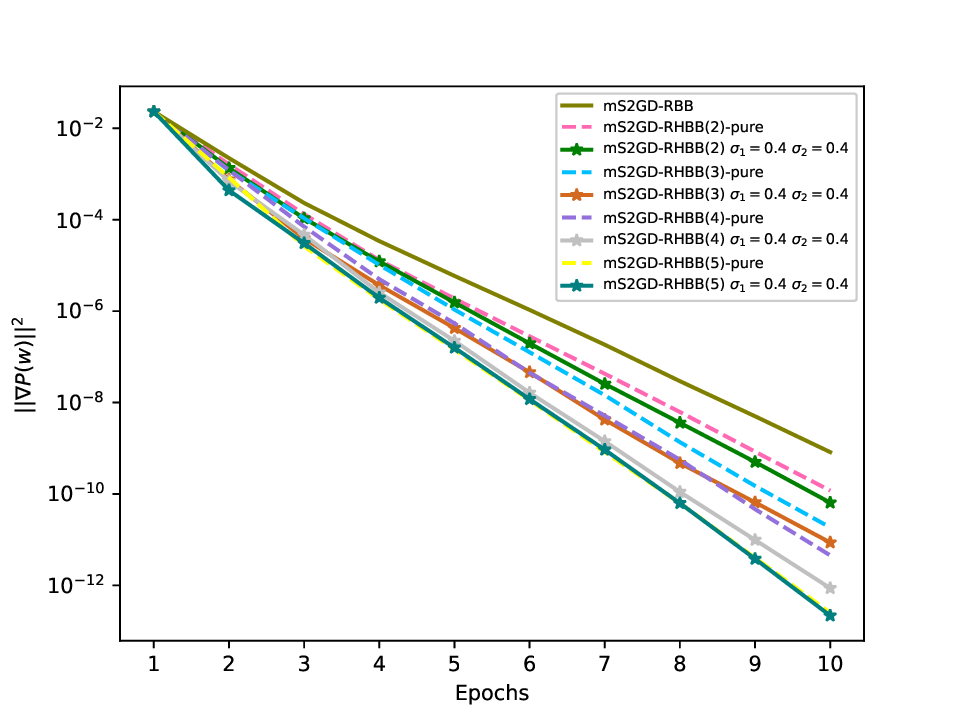}}
  		\caption{\footnotesize The performance of mS2GD-RBB, non-adaptive mS2GD-RHBB, and adaptive mS2GD-RHBB under $\sigma_1=0.4$, $\sigma_2=0.4$.}
  		\label{fig22}
  	\end{figure*}

  	In Figs. \ref{fig17} - \ref{fig22}, we compare the RBB rule, the non-adaptive RHBB rule and the adaptive RHBB rule in terms of the evolution of $\|\nabla P(\cdot)\|^2$. Note that the adaptive pair $(\sigma_1$, $\sigma_2)$ in the iterative adaptor $h(\cdot)$ is explored triply, with $(0.6, 0.2)$ in Fig. \ref{fig17}, \ref{fig18}, $(0.7, 0.1)$ in Fig. \ref{fig19}, \ref{fig20} and $(0.4, 0.4)$ in Fig. \ref{fig21}, \ref{fig22}. 
  	
  	Figs. \ref{fig17}, \ref{fig19}, \ref{fig21} show that the adaptive MB-SARAH-RHBB consistently outperforms the `non-adaptive' MB-SARAH-RHBB and surpass the original MB-SARAH-RBB by a large margin. Similar and consistent results can be seen for the adaptive mS2GD-RHBB as illustrated in Figs. \ref{fig18}, \ref{fig20}, \ref{fig22}.
  	
  	In most cases, $h(\cdot)$ provides significant speedup for algorithms using the hedge base $\alpha=4, 5$. It's also worth noticing that under $\alpha=4, 5$, the associated algorithms perform equally well on $w8a$. Combining the previous results from Figs \ref{fig4},  \ref{fig5}, it is evident that the optimal hedge magnitude of RHBB is slightly lower on $w8a$ than on the other data sets. 
  
    The performance of algorithms varies significantly over the hyper-parameter pair $(\sigma_1$, $\sigma_2)$, however, we do not require a strict guideline for the selection, our three casual and moderate choices have proven to be sufficiently effective.
  	
  	Indeed, we fix an exponential adaptor at the start for the convenience in paper. A simple incremental function has already resulted in noticeable improvements. We believe that if more research is conducted (e.g. using sigmoid increments of $h(\cdot)$ or choosing the non-exponential adaptors), our algorithms can be accelerated further.
  	
  	\subsubsection{Adaptive MB-SARAH-RHBB+/mS2GD-RHBB+}
  	\ 
  	\vskip10pt
  	\noindent  \textbf{Parametric Settings:} We set $b=4$, the unified $b_H=40$ and sample the subsets $S_1$ and $S_2$ according to distributions $Q$, where $Q$ are configured by option I and option II. Analogously, we set $\tau=2$ in both option I and option II. To avoid possible over-utility, we implement $\gamma=0.8$ and $\gamma_2=0.8$. We employ the first values in the adaptive pair as $(\sigma_1, \sigma_2)=(0.6, 0.2)$. Eventually, the hedge base $\alpha$ is selected within the same set $\{2, 3, 4, 5\}$.
  	
\begin{figure*}[htbp]
  		\centering
  		\subfigure[australian]
  		{\includegraphics[width=0.327\textwidth]{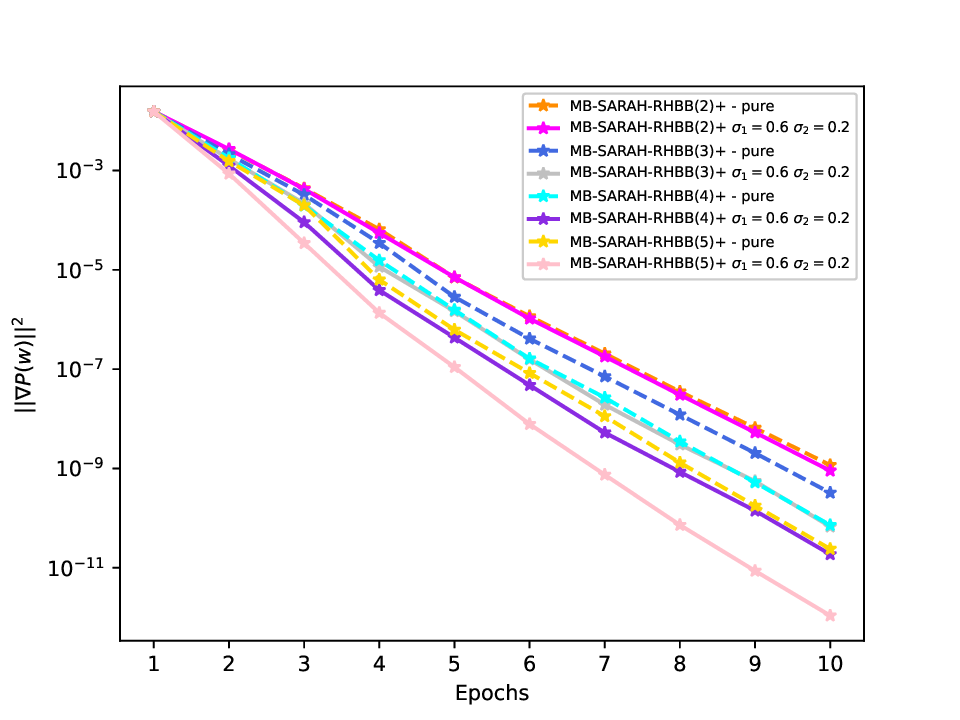}}
  		\subfigure[madelon]
  		{\includegraphics[width=0.327\textwidth]{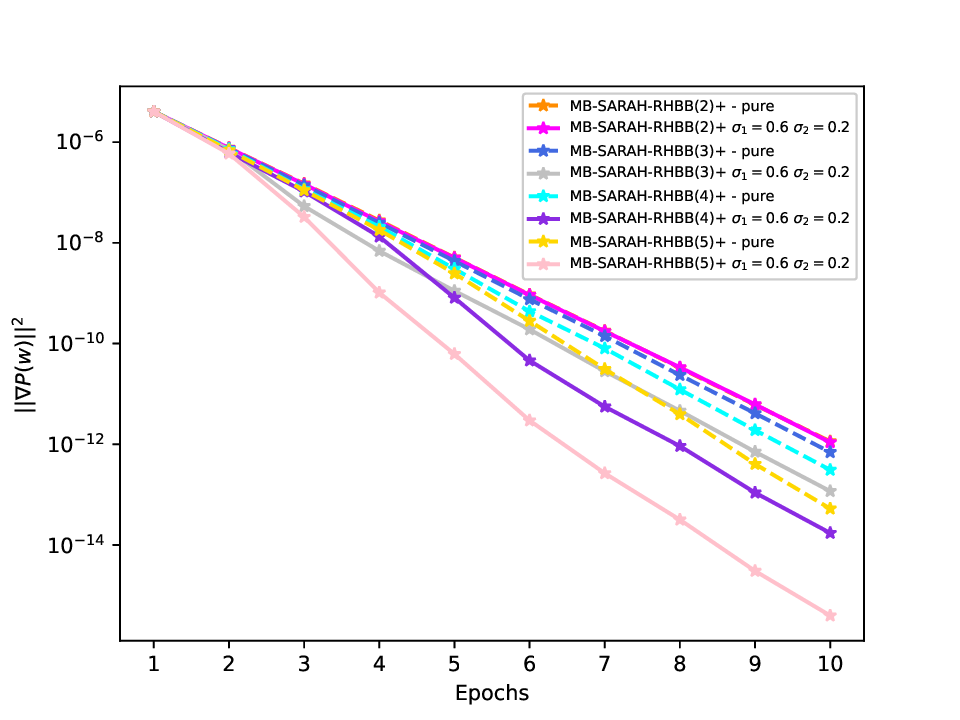}}
  		\subfigure[german.numer]
  		{\includegraphics[width=0.327\textwidth]{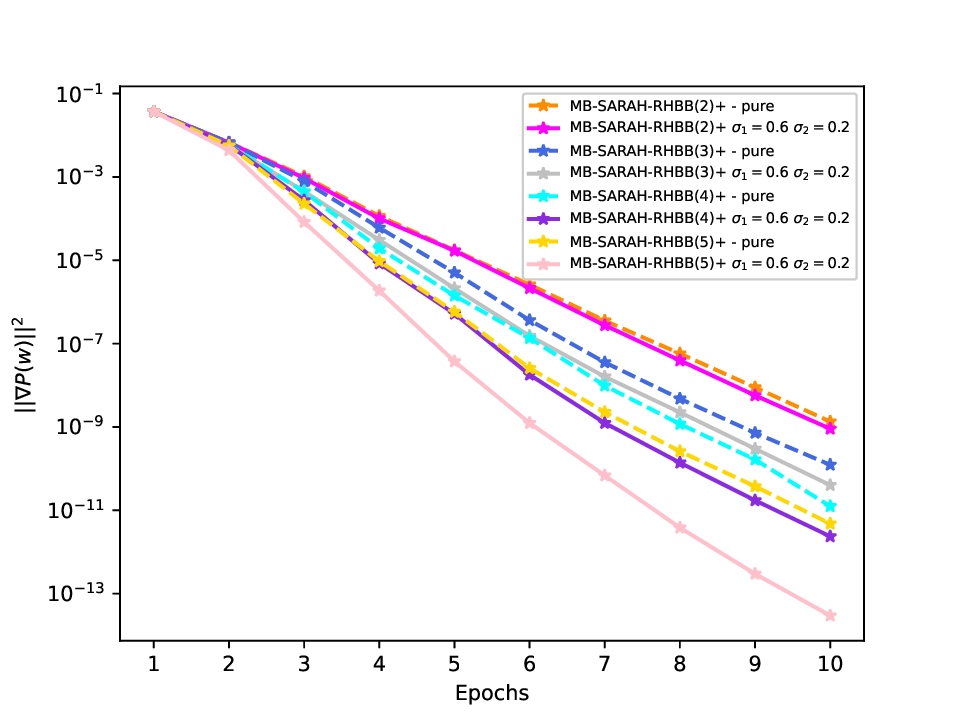}}
  		\caption{\footnotesize The performance of MB-SARAH-RHBB+ with $\sigma_1=0.6$, $\sigma_2=0.2$ and non-adaptive MB-SARAH-RHBB+. $Q$ is configured under \textbf{option I}.}
  		\label{fig23}
  	\end{figure*}

\begin{figure*}[htbp]
  		\centering
  		\subfigure[australian]
  		{\includegraphics[width=0.327\textwidth]{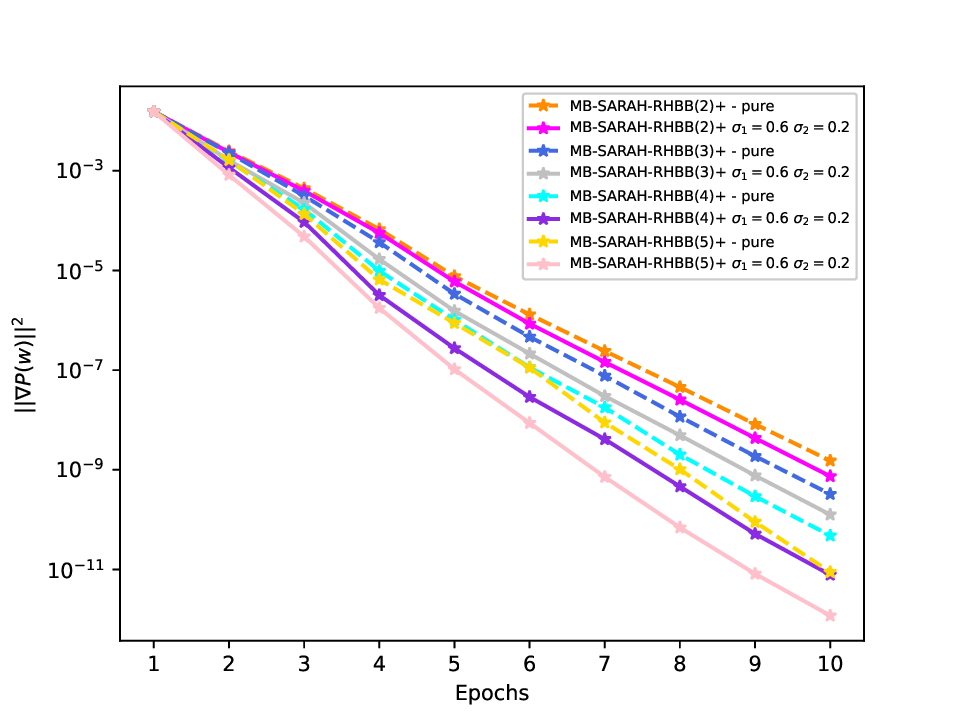}}
  		\subfigure[madelon]
  		{\includegraphics[width=0.327\textwidth]{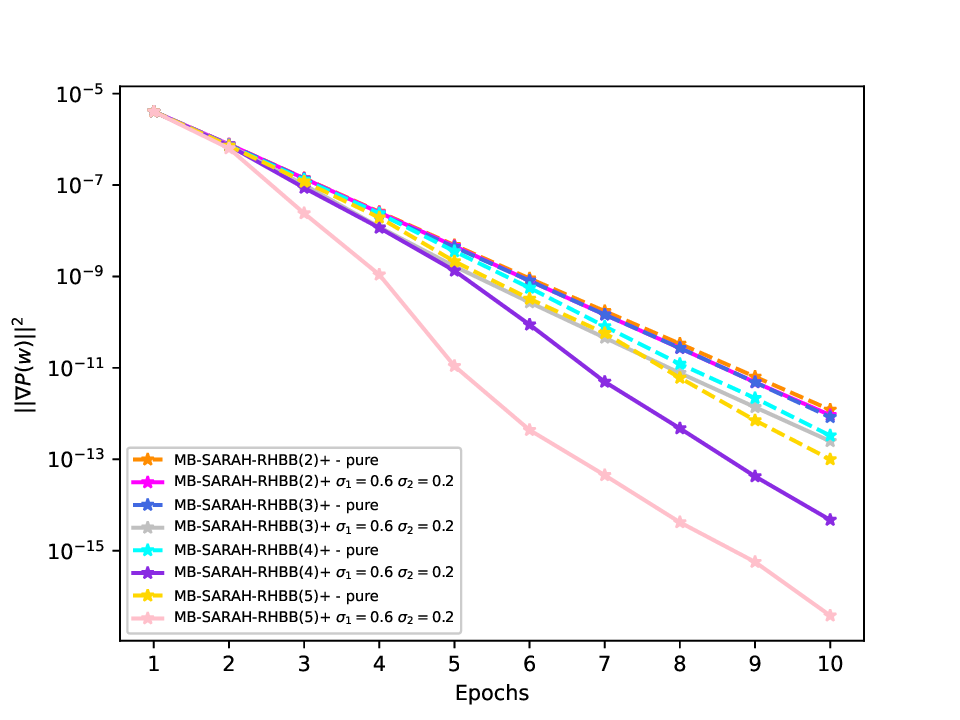}}
  		\subfigure[german.numer]
  		{\includegraphics[width=0.327\textwidth]{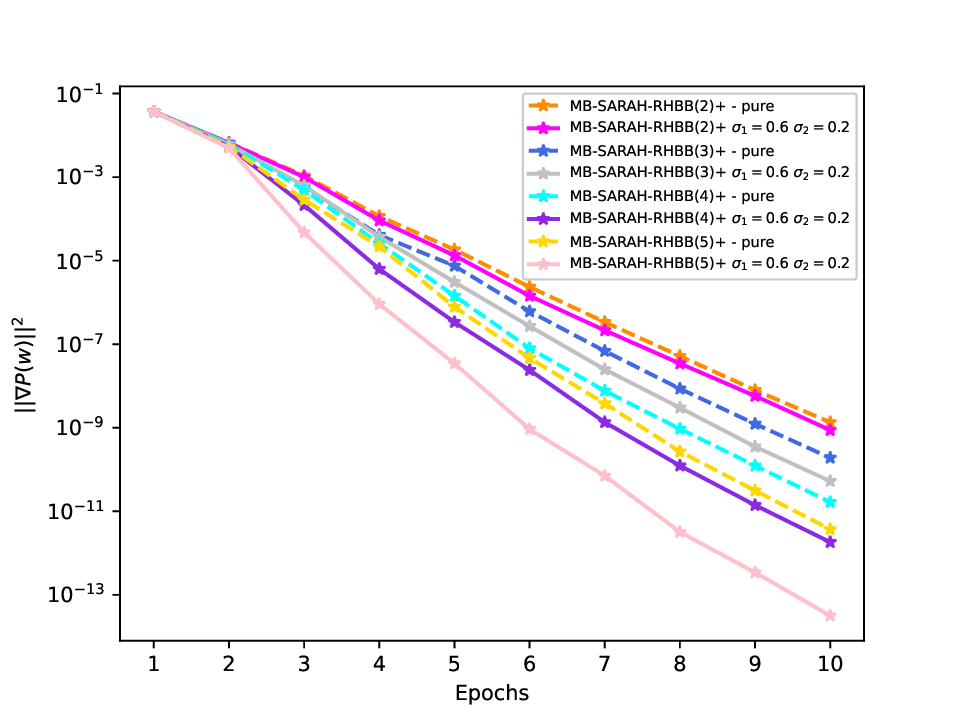}}
  		\caption{\footnotesize The performance of MB-SARAH-RHBB+ with $\sigma_1=0.6$, $\sigma_2=0.2$ and non-adaptive MB-SARAH-RHBB+. $Q$ is configured under \textbf{option II}.}
  		\label{fig24}
  	\end{figure*}

\begin{figure*}[htbp]
  		\centering
  		\subfigure[australian]
  		{\includegraphics[width=0.327\textwidth]{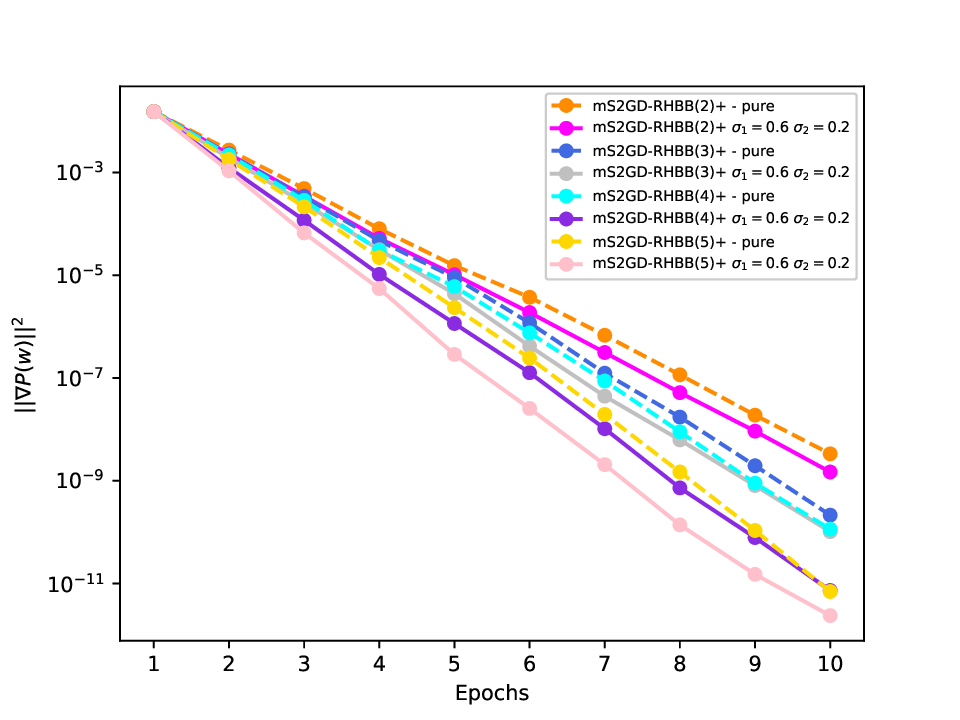}}
  		\subfigure[madelon]
  		{\includegraphics[width=0.327\textwidth]{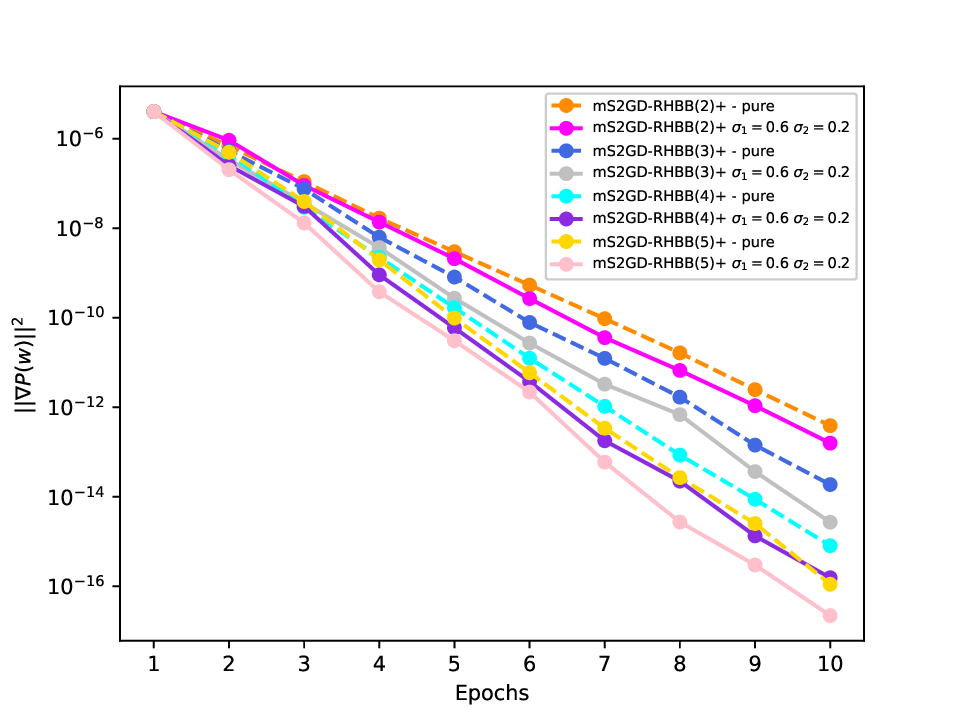}}
  		\subfigure[german.numer]
  		{\includegraphics[width=0.327\textwidth]{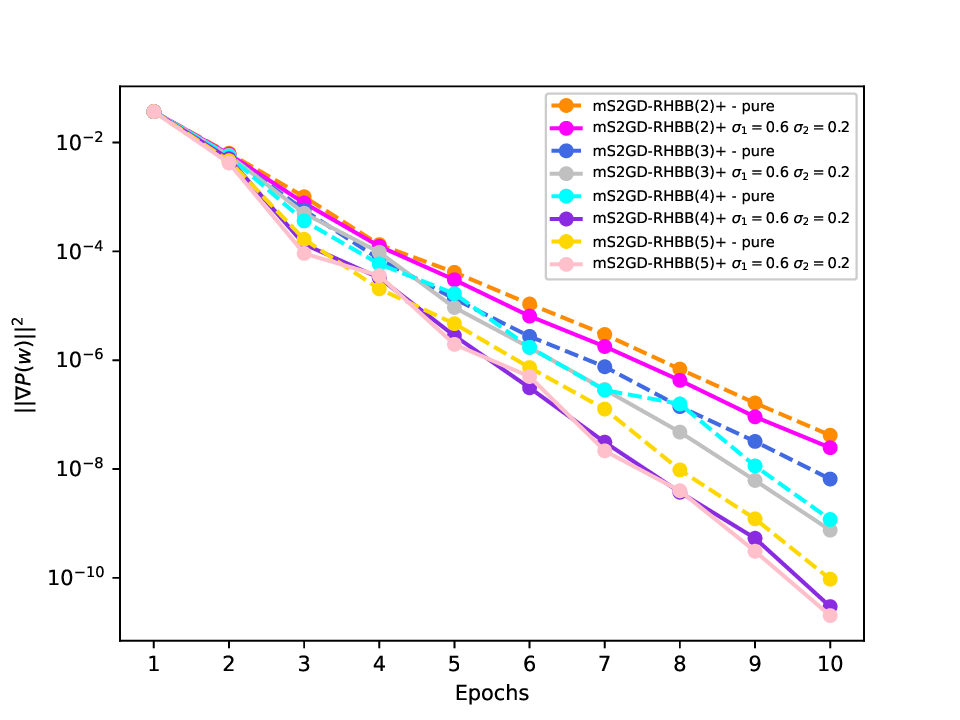}}
  		\caption{\footnotesize The performance of mS2GD-RHBB+ with $\sigma_1=0.6$, $\sigma_2=0.2$ and non-adaptive mS2GD-RHBB+. $Q$ is configured under \textbf{option I}.}
  		\label{fig25}
  	\end{figure*}

\begin{figure*}[htbp]
  		\centering
  		\subfigure[australian]
  		{\includegraphics[width=0.327\textwidth]{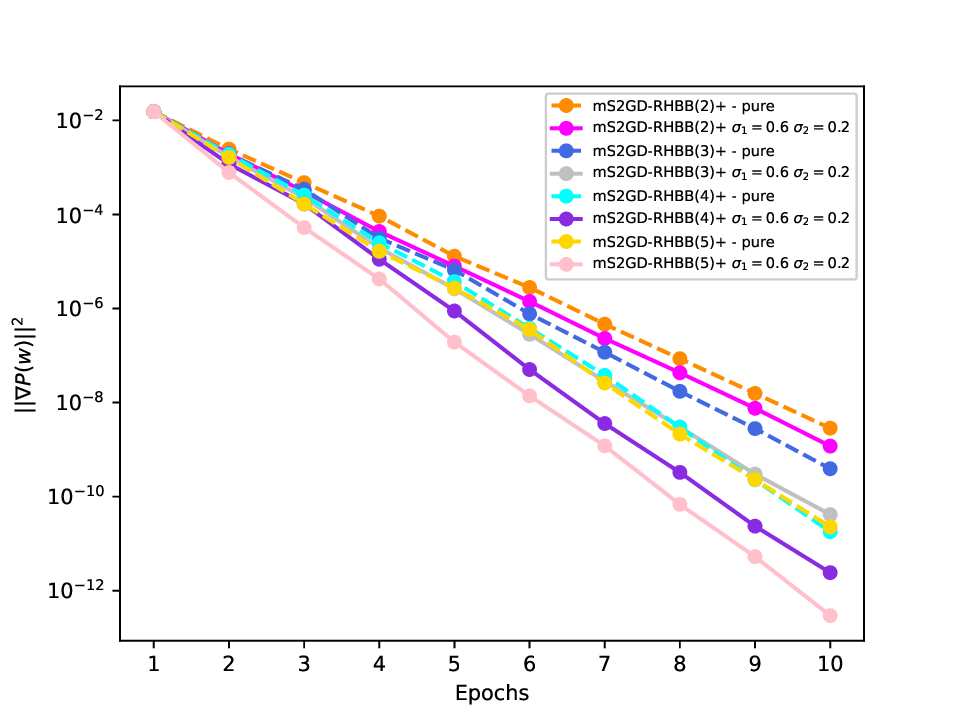}}
  		\subfigure[madelon]
  		{\includegraphics[width=0.327\textwidth]{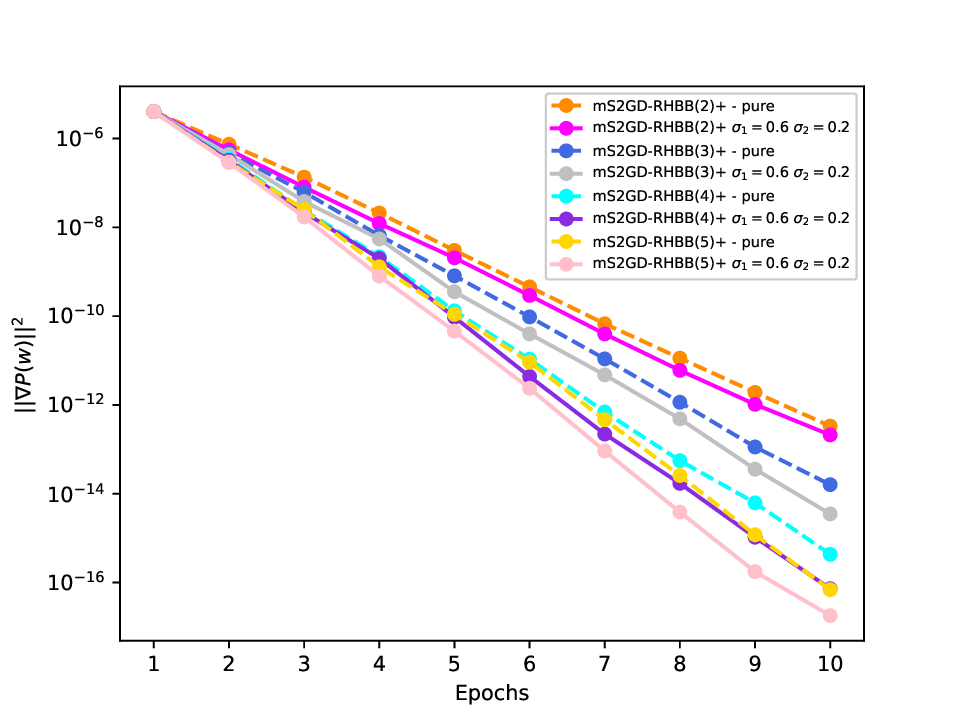}}
  		\subfigure[german.numer]
  		{\includegraphics[width=0.327\textwidth]{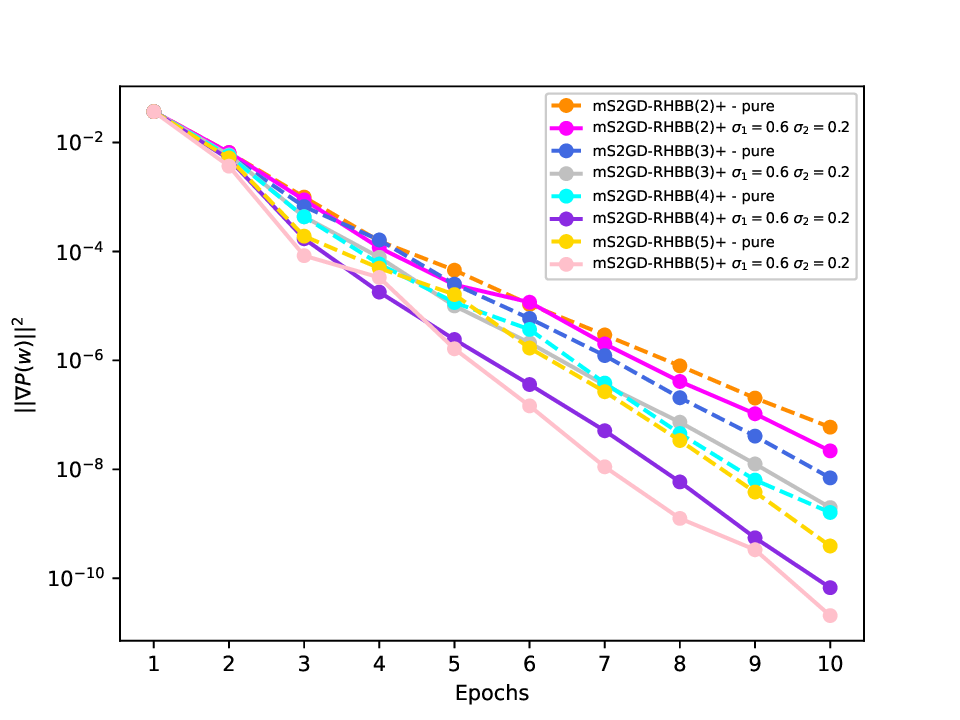}}
  		\caption{\footnotesize The performance of mS2GD-RHBB+ with $\sigma_1=0.6$, $\sigma_2=0.2$ and non-adaptive mS2GD-RHBB+. $Q$ is configured under \textbf{option II}.}
  		\label{fig26}
  	\end{figure*}

  	Fig. \ref{fig23}, \ref{fig24} exhibit that the practical speeds of adaptive MB-SARAH-RHBB+ are faster than `non-adaptive' MB-SARAH-RHBB+, under different hedge bases and different distribution options. Especially when $\alpha=5$, the overall improvements are particularly remarkable, which begins from the start of iterations. We have the corresponding results for adaptive mS2GD-RHBB+ in Fig. \ref{fig25}, \ref{fig26}, showing that adaptive mS2GD-RHBB+ outperforms `non-adaptive' mS2GD-RHBB+ consistently. Therefore, it's reasonable as well advisable to equip an iterative adaptor to achieve additional accelerations in early periods. 
  	
    Fig. \ref{fig23} - \ref{fig26} corroborate the previous conclusions that the importance sampling is more responsive under intense hedge scenarios ( relatively large $\alpha$).
    
  	Massive results in this subsection have suggested that $h(\cdot)$ is instrumental in completing an efficient step size rule, for it addresses the defect of inflexibility in stochastic algorithms. The consistent performance implies the accordance between the importance sampling and the iterative scaling, dispelling potential concerns about discrepancies in overall optimization.
  	
  	The current distributions (option I and option II) are especially productive on $phishing$, $mushrooms$ and $german.numer$. Still, practitioners can configure particular distributions to match targeted sets to their needs.
  	
  	\subsubsection{Comparison with other state-of-art methods}
  	\ 
  	\vskip10pt
  	\noindent  \textbf{Parametric Settings:} In MB-SARAH-RHBB and mS2GD-RHBB, we set $b=4$, the unified $b_H=40$ and sample subsets $S$, $S_1$, $S_2$ according to uniform distribution. We set $\gamma=1$, $\gamma_2=1$ and decide the adaptive pair $(\sigma_1,\sigma_2)=(0.8, 0.8)$ for a fresh try. Eventually, we opt $\alpha=3$ as a gentle hedge base.  

\begin{figure*}[htbp]
  		\centering
  		\subfigure[a8a]
  		{\includegraphics[width=0.45\textwidth]{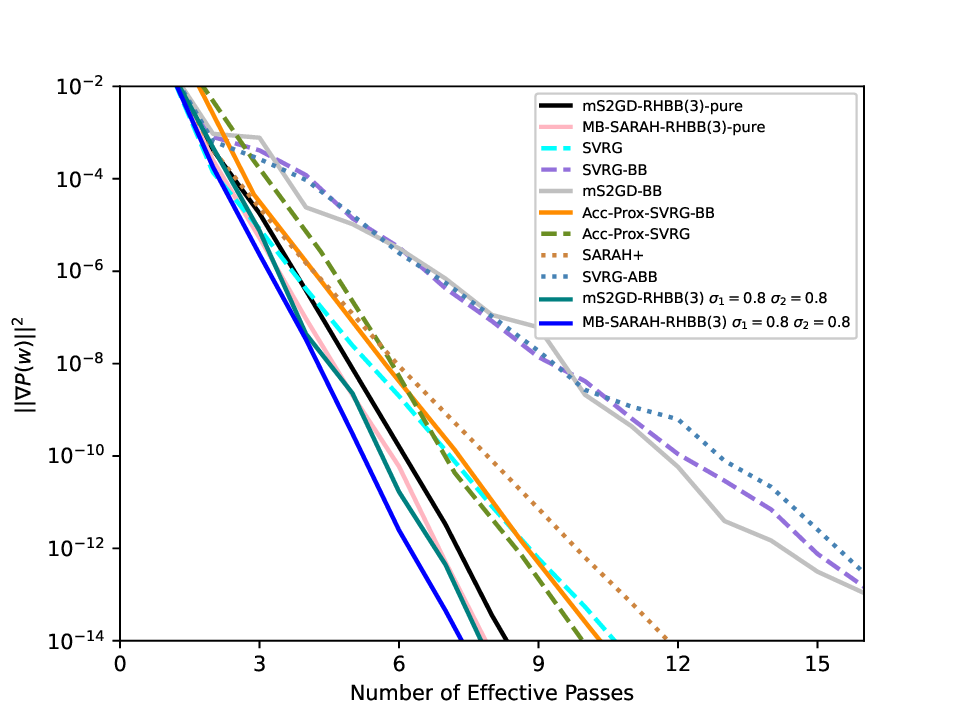}}
  		\subfigure[w8a]
  		{\includegraphics[width=0.45\textwidth]{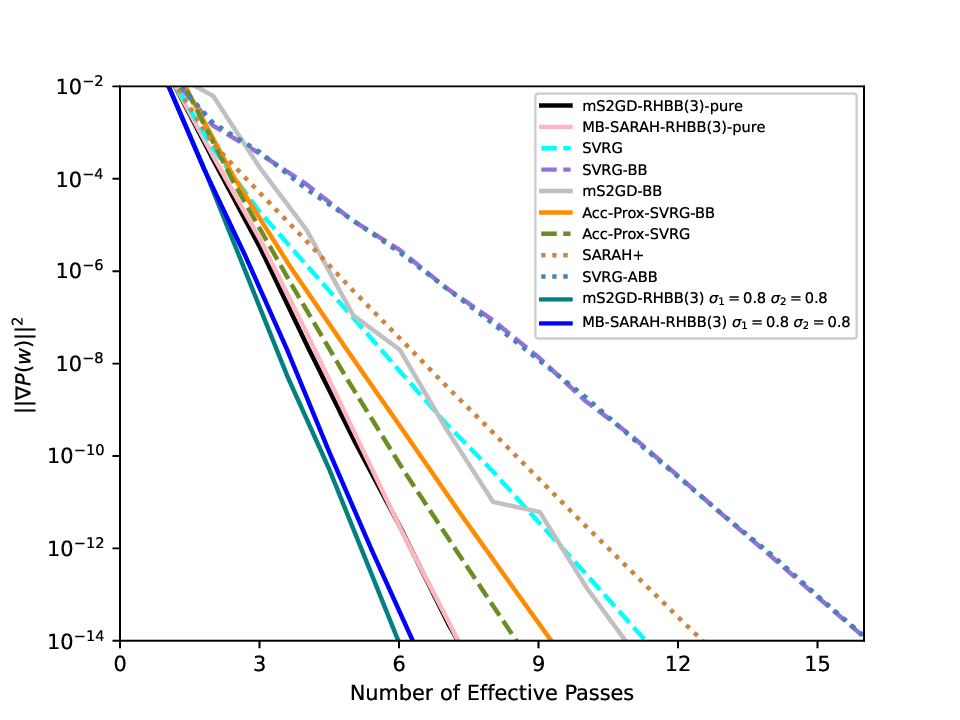}}
  		
  		\subfigure[ijcnn1]
  		{\includegraphics[width=0.45\textwidth]{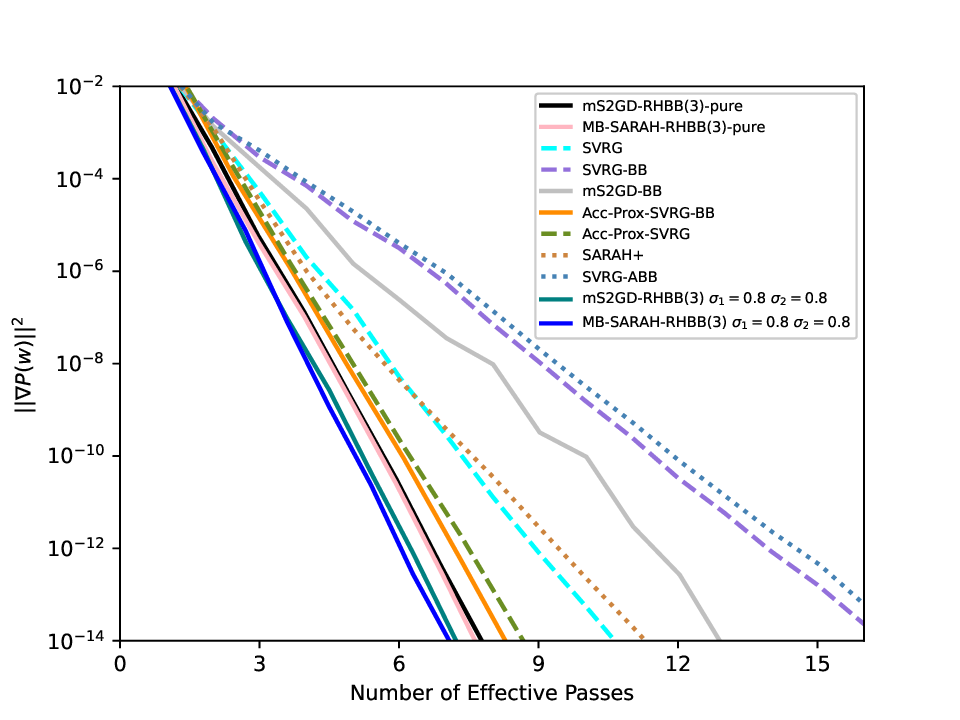}}
  		\subfigure[covtype]
  		{\includegraphics[width=0.45\textwidth]{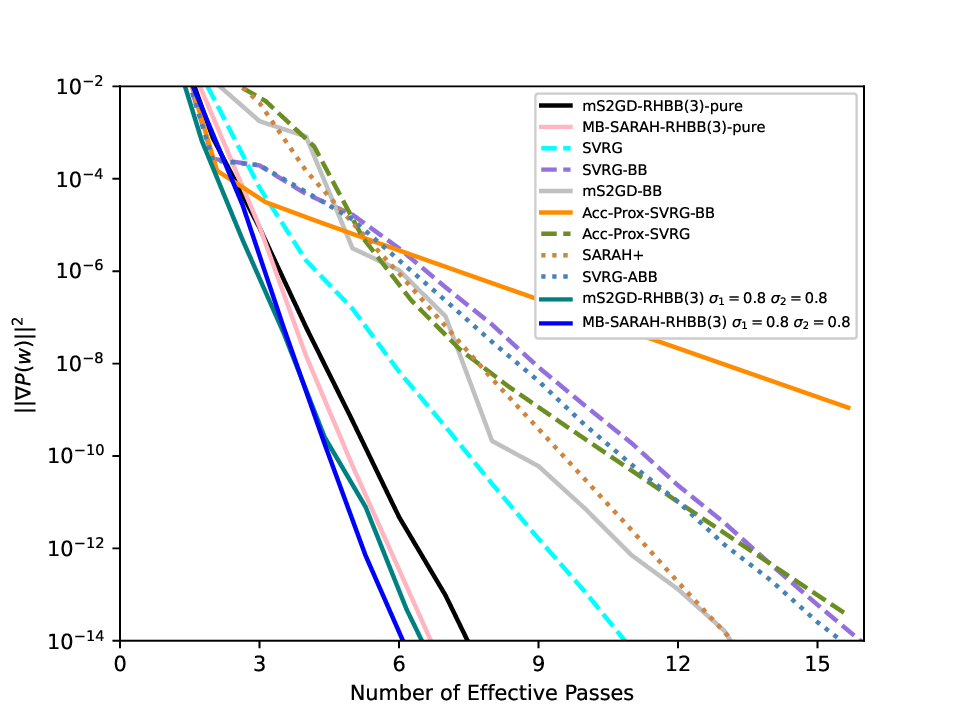}}
  		
  		\subfigure[phishing]
  		{\includegraphics[width=0.45\textwidth]{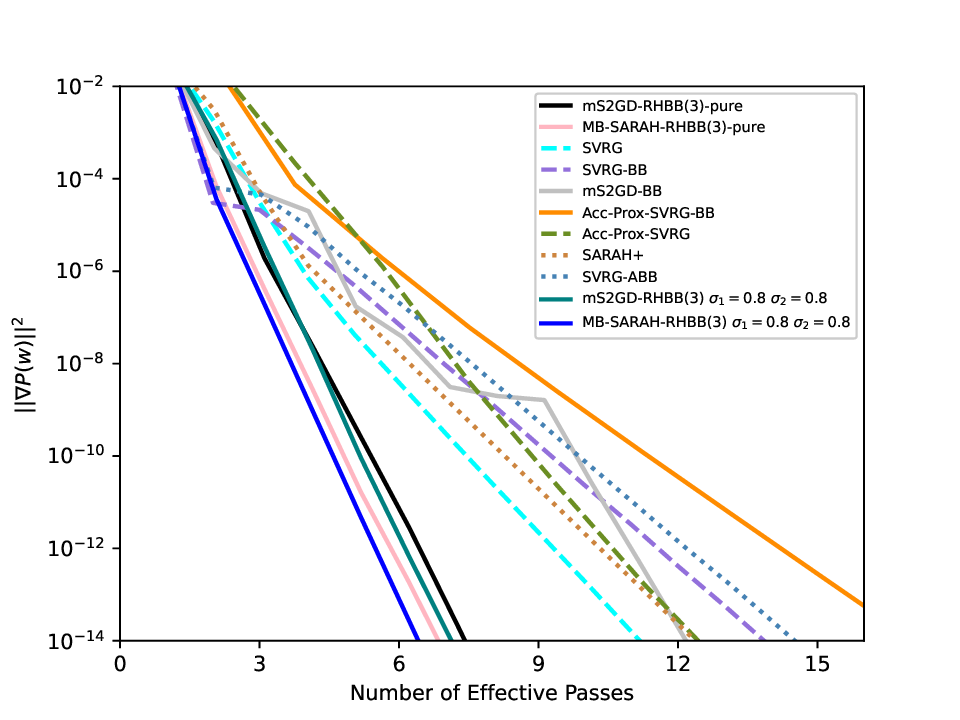}}
  		\subfigure[mushrooms]
  		{\includegraphics[width=0.45\textwidth]{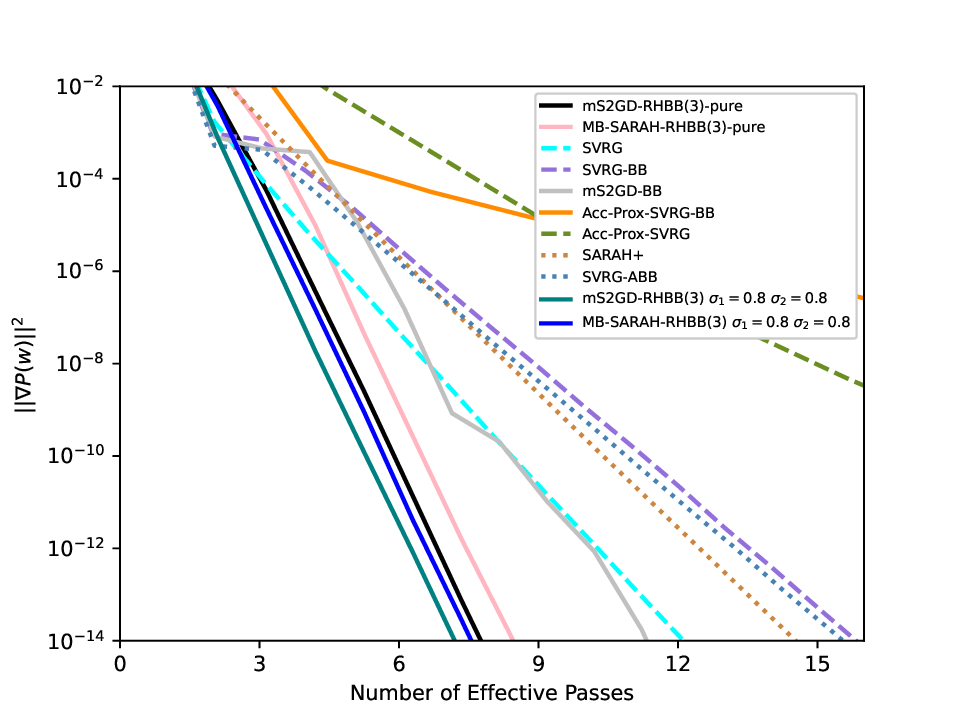}}
  		\caption{\footnotesize Comparisons of different algorithms.}
  		\label{fig27}
  	\end{figure*}

  	As can be seen from Fig \ref{fig27}, all of our adaptive and `non-adaptive' methods outperform various state-of-the-art algorithms.

  	\subsection{Investigation on Batch Sizes}
\begin{figure*}[htbp]
  		\centering
  		\subfigure[non-adaptive]
  		{\includegraphics[width=0.4\textwidth]{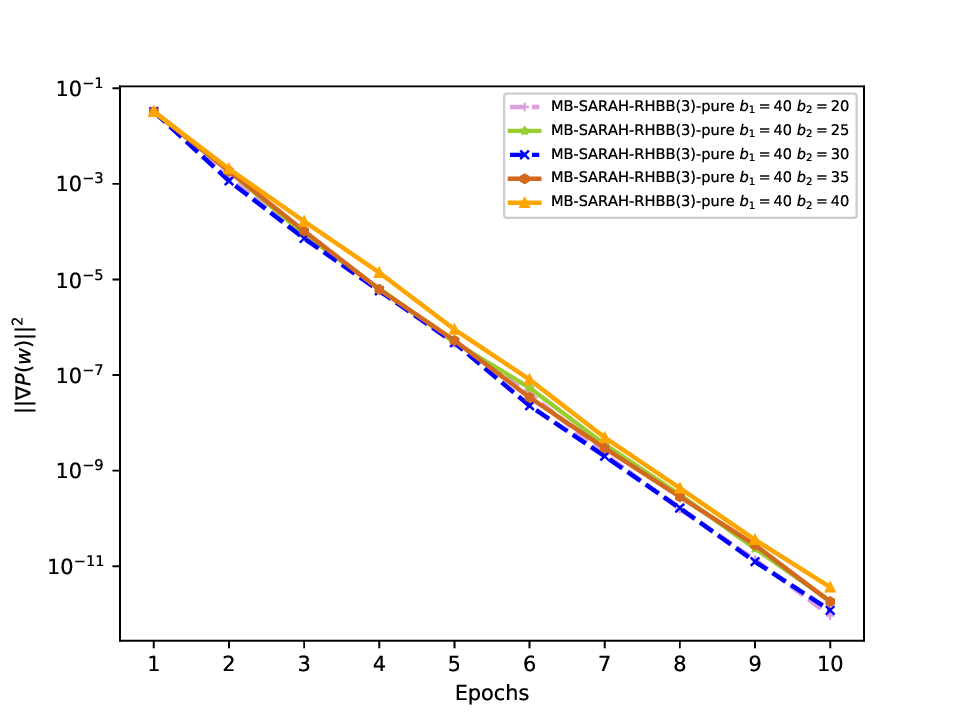}}
  		\subfigure[adaptive]
  		{\includegraphics[width=0.4\textwidth]{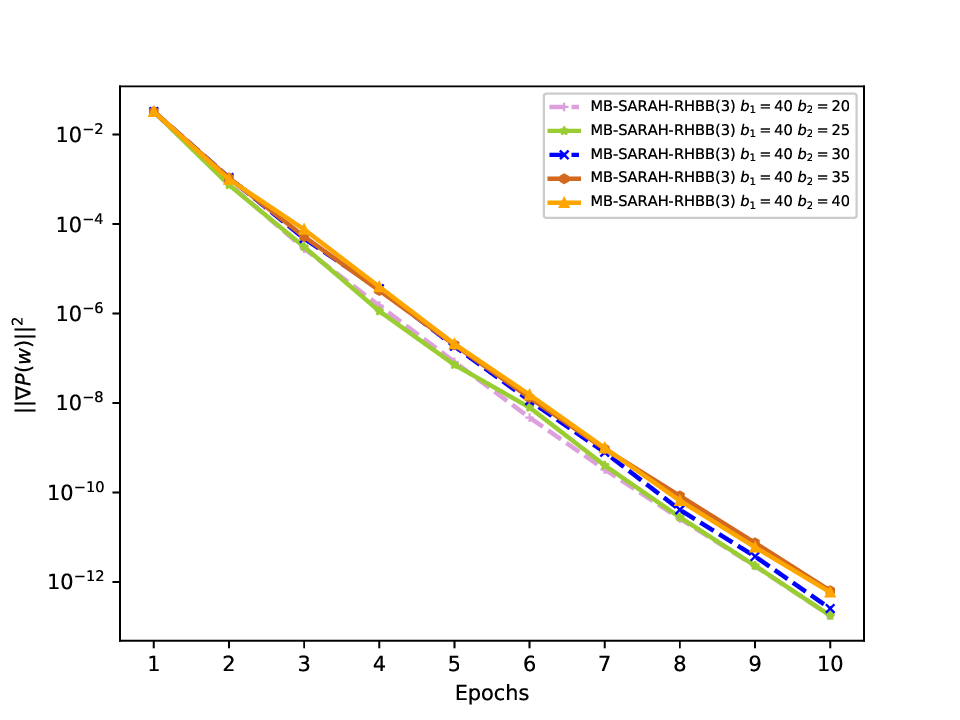}}
  		\caption{\footnotesize The performance of non-adaptive MB-SARAH-RHBB(3) and adaptive MB-SARAH-RHBB(3) ($\sigma_1=0.6, \sigma_2=0.2$), with different batch sizes on $a8a$.}
  		\label{fig28}
  	\end{figure*}

\begin{figure*}[htbp]
  		\centering
  		\subfigure[non-adaptive]
  		{\includegraphics[width=0.4\textwidth]{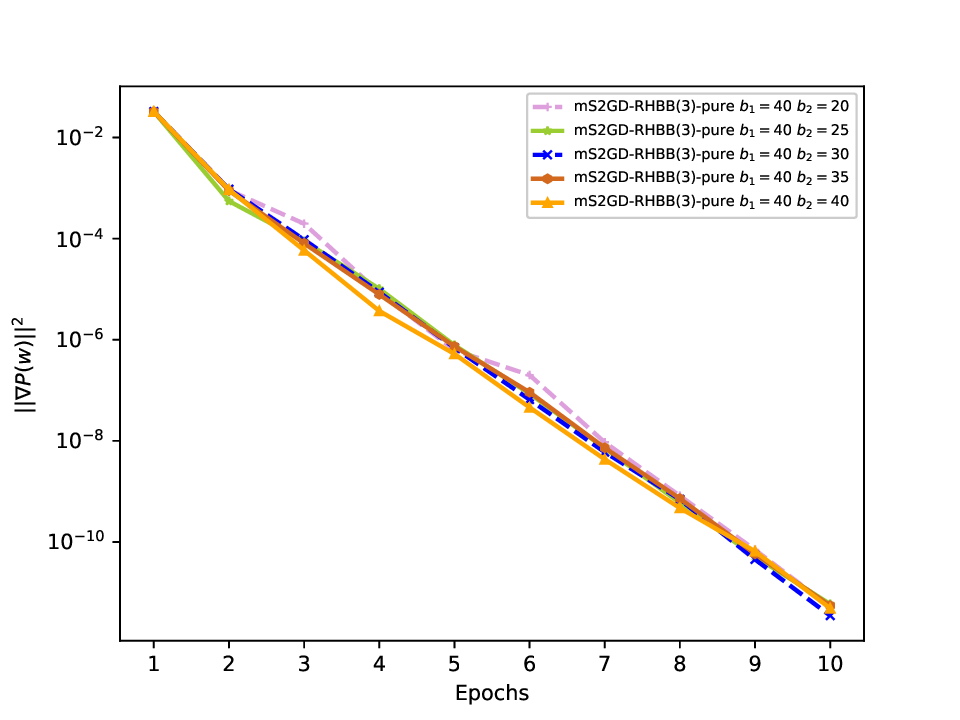}}
  		\subfigure[adaptive]
  		{\includegraphics[width=0.4\textwidth]{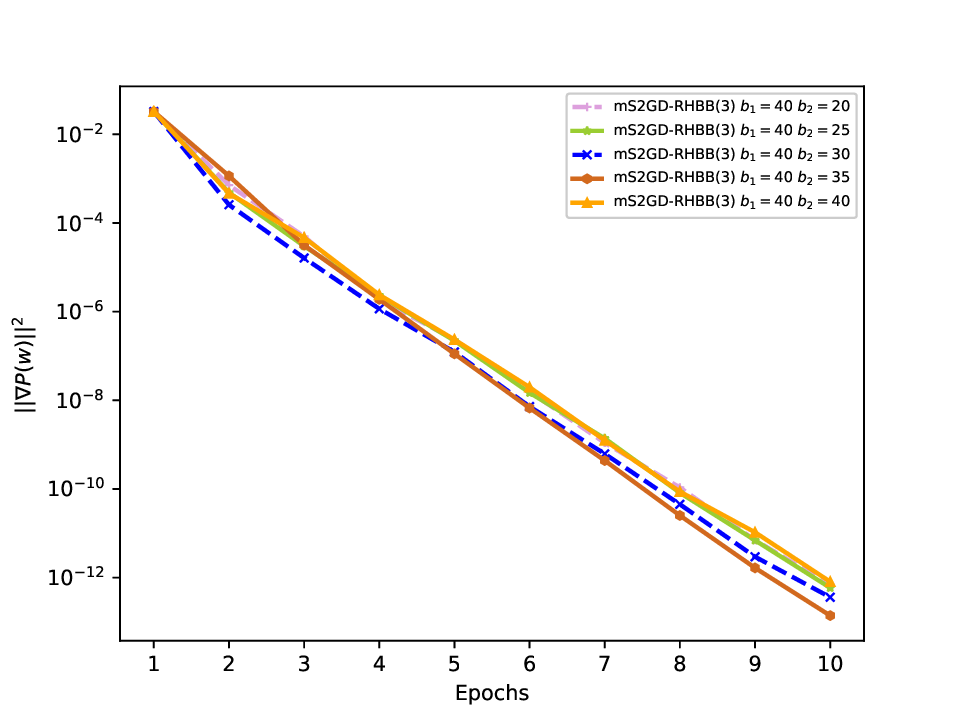}}
  		\caption{\footnotesize The performance of non-adaptive mS2GD-RHBB(3) and adaptive mS2GD-RHBB(3) ($\sigma_1=0.6, \sigma_2=0.2$), with different batch sizes on $a8a$.}
  		\label{fig29}
  	\end{figure*}

\begin{figure*}[htbp]
  		\centering
  		\subfigure[non-adaptive]
  		{\includegraphics[width=0.4\textwidth]{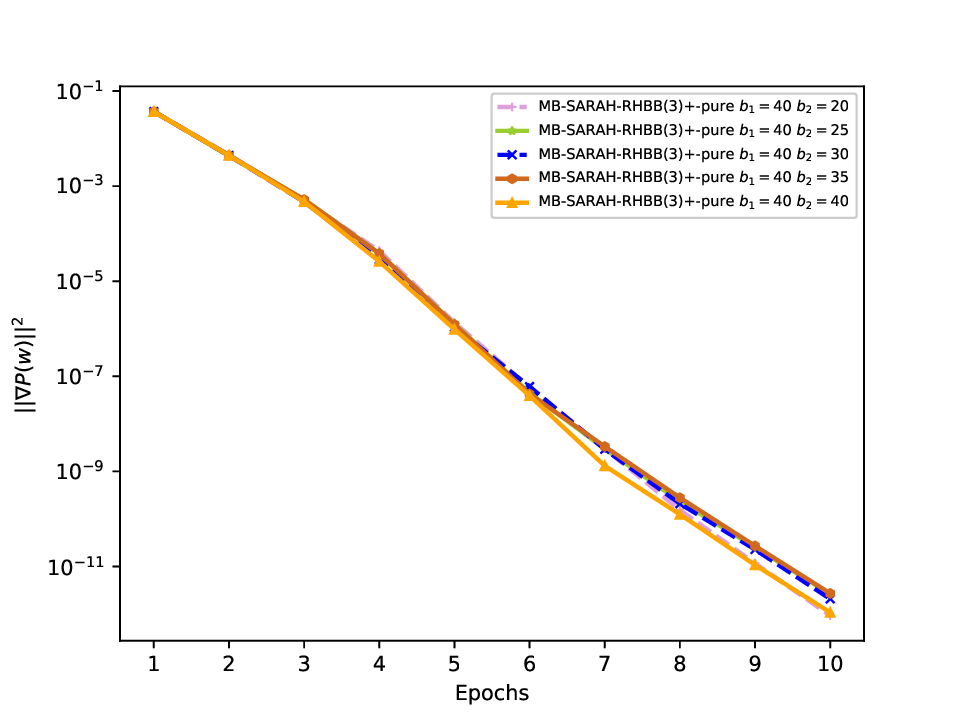}}
  		\subfigure[adaptive]
  		{\includegraphics[width=0.4\textwidth]{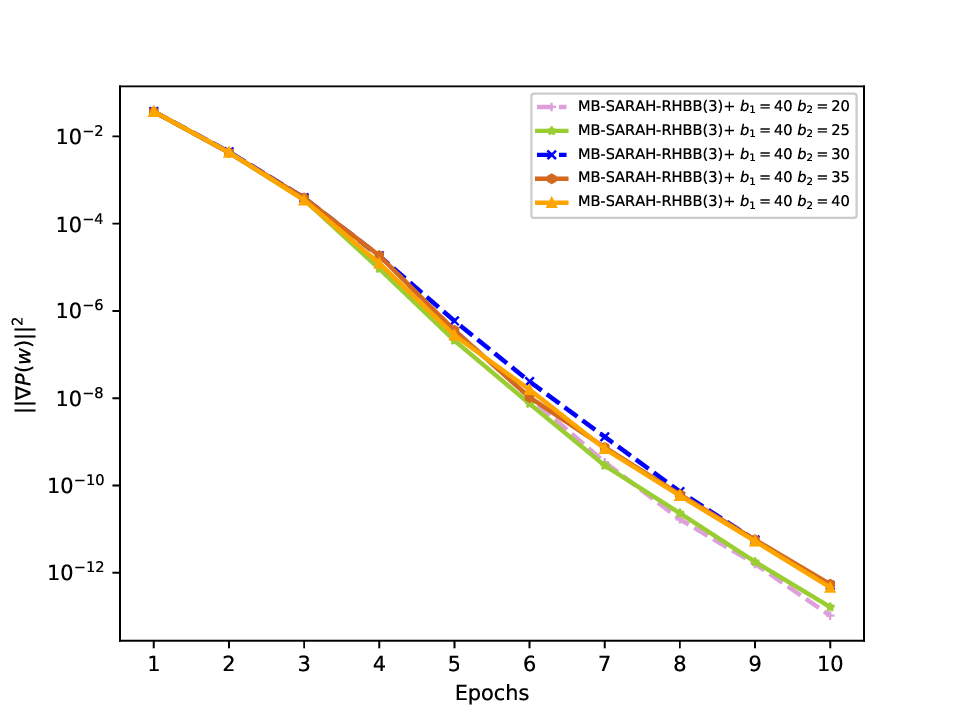}}
  		\caption{\footnotesize  The performance of non-adaptive MB-SARAH-RHBB(3)+ and adaptive MB-SARAH-RHBB(3)+ ($\sigma_1=0.6, \sigma_2=0.2$), with different batch sizes on $a8a$. $Q$ is configured under \textbf{option I}.}
  		\label{fig30}
  	\end{figure*}

\begin{figure*}[htbp]
  		\centering
  		\subfigure[non-adaptive]
  		{\includegraphics[width=0.4\textwidth]{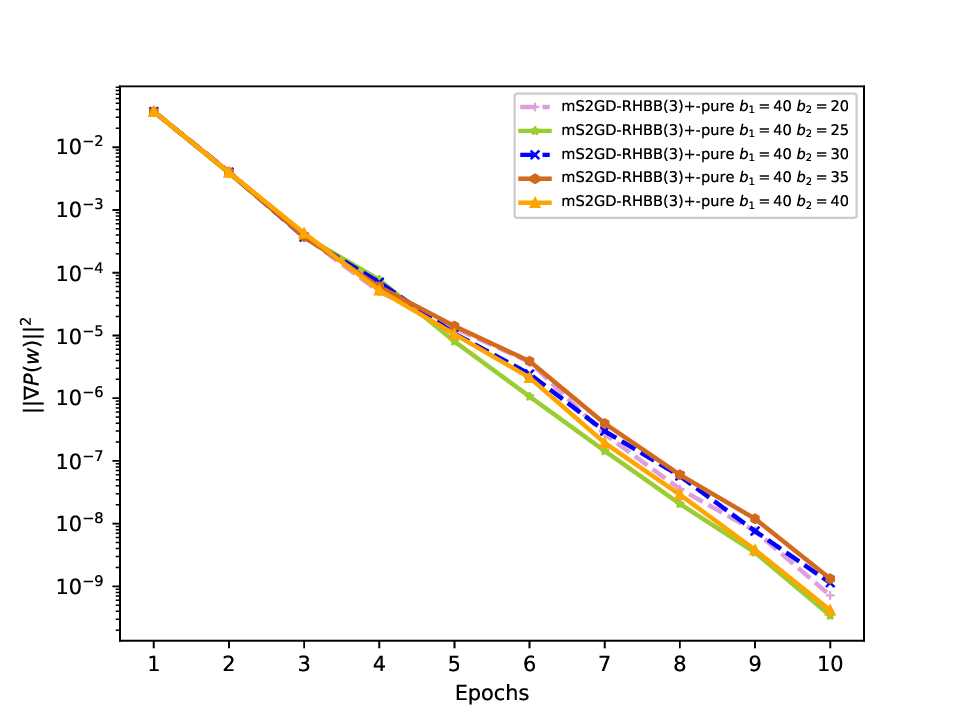}}
  		\subfigure[adaptive]
  		{\includegraphics[width=0.4\textwidth]{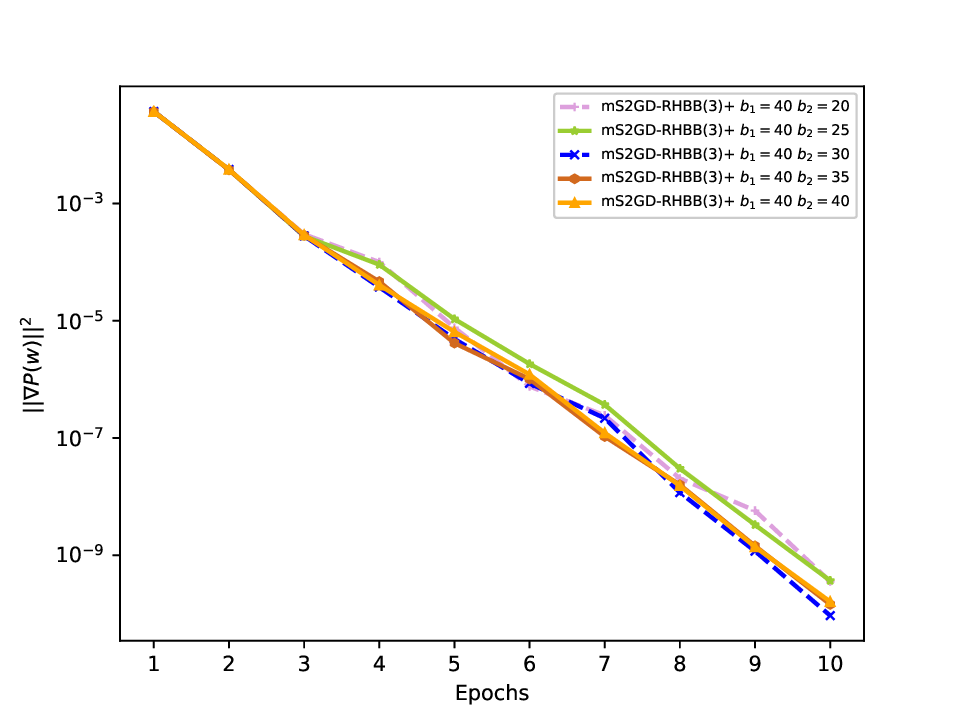}}
  		\caption{\footnotesize The performance of non-adaptive mS2GD-RHBB(3)+ and adaptive mS2GD-RHBB(3)+ ($\sigma_1=0.6, \sigma_2=0.2$), with different batch sizes on $a8a$. $Q$ is configured under \textbf{option I}}
  		\label{fig31}
  	\end{figure*}

\begin{figure*}[htbp]
  		\centering
  		\subfigure[non-adaptive]
  		{\includegraphics[width=0.4\textwidth]{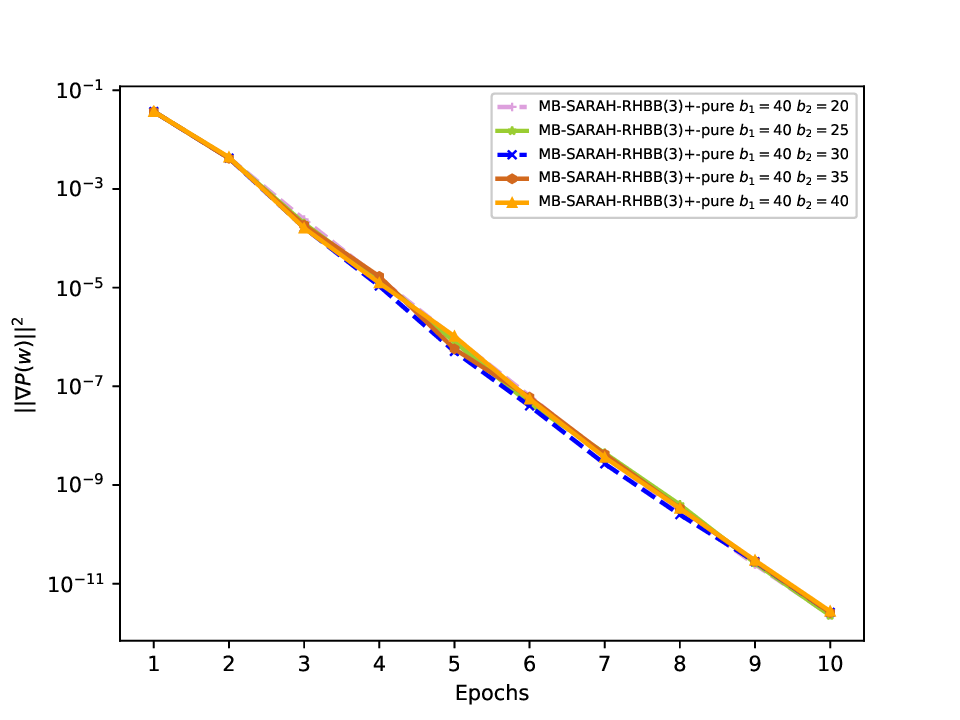}}
  		\subfigure[adaptive]
  		{\includegraphics[width=0.4\textwidth]{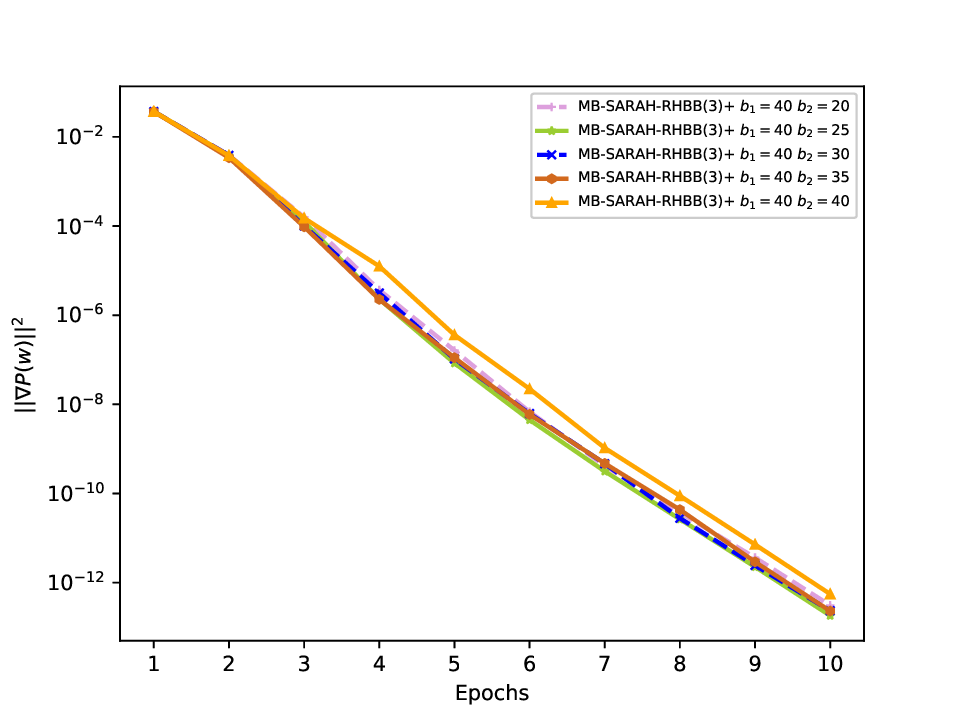}}
  		\caption{\footnotesize The performance of non-adaptive MB-SARAH-RHBB(3)+ and adaptive MB-SARAH-RHBB(3)+ ($\sigma_1=0.6, \sigma_2=0.2$), with different batch sizes on $a8a$. $Q$ is configured under \textbf{option II}.}
  		\label{fig32}
  	\end{figure*}

\begin{figure*}[htbp]
  		\centering
  		\subfigure[non-adaptive]
  		{\includegraphics[width=0.4\textwidth]{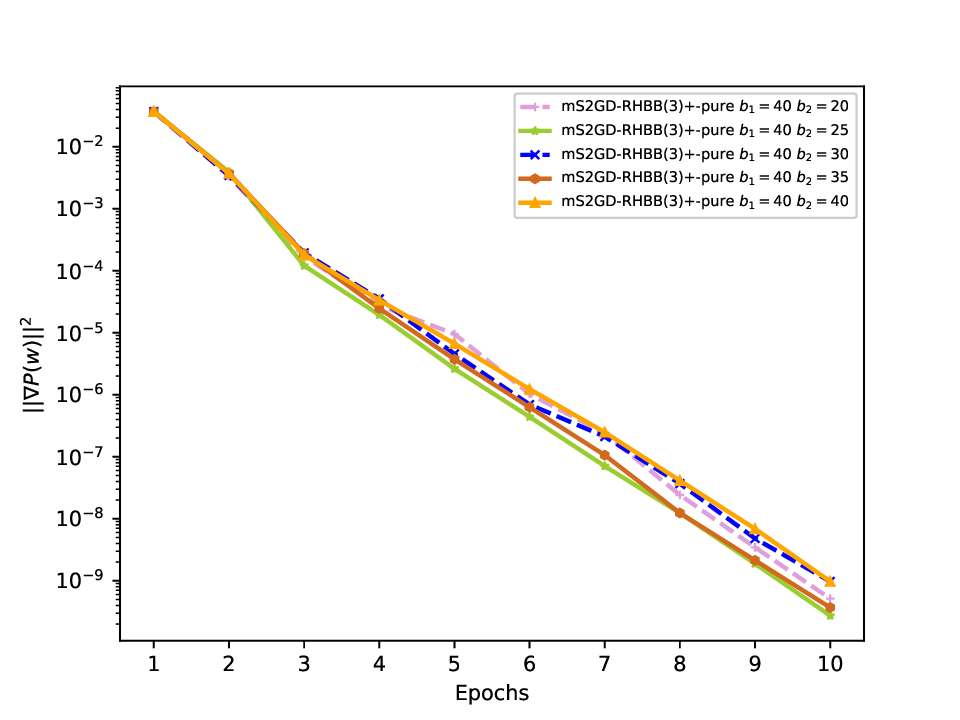}}
  		\subfigure[adaptive]
  		{\includegraphics[width=0.4\textwidth]{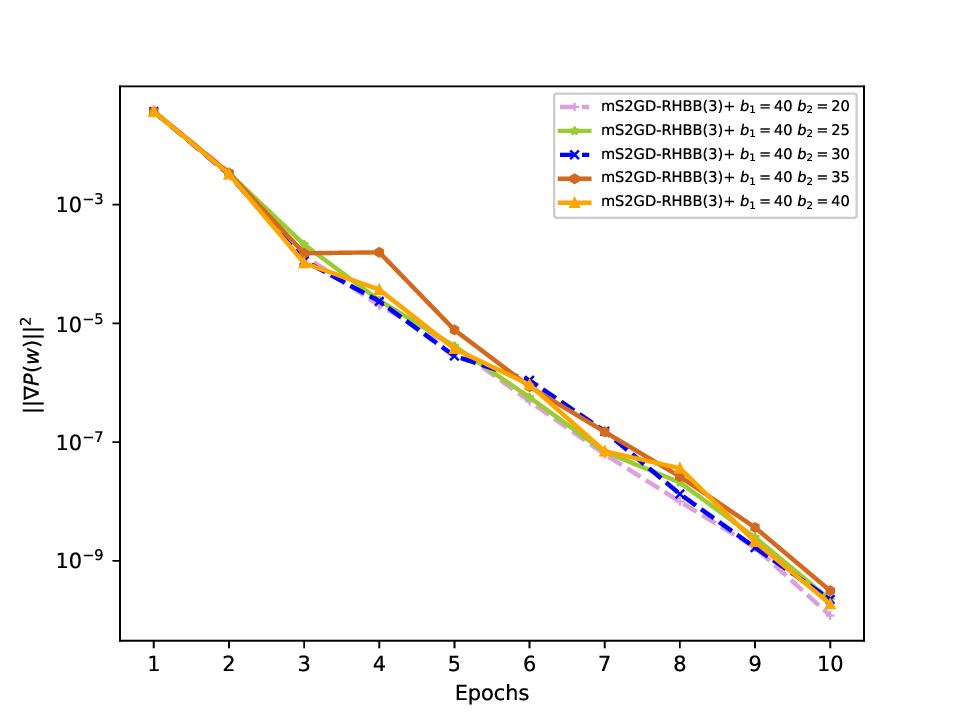}}
  		\caption{\footnotesize The performance of non-adaptive mS2GD-RHBB(3)+ and adaptive mS2GD-RHBB(3)+ ($\sigma_1=0.6, \sigma_2=0.2$), with different batch sizes on $a8a$. $Q$ is configured under \textbf{option II}.}
  		\label{fig33}
  	\end{figure*}

  	We technically supply this subsection to demonstrate that the performance of our algorithms is not sensitive to $b_1$ and $b_2$ under the same $b_H$. Here, we arrange $b_1=40$ and $b_2=20, 25, 30, 35, 40$ alternately for verification (all settings must ensure the batch correction $\overline{b}=\max\{b_1, b_2\}=40$ unvaried). We display the results from Fig. \ref{fig28} to Fig. \ref{fig33}. Hence, in all the previous experiments, setting the unified batch size of $b_1=b_2=b_H$ is a reasonable also economical choice.
  	
  	\section{Conclusion}
  	In this paper, we propose two novel and efficient rules for stochastic optimization, which are motivated by the random Barzilai-Borwein method, the important sampling technique and modern iterative adaptors. The idea of design is aggressive yet robust: by leveraging untapped curvature, we enlarge the random Barzilai-Borwein step sizes effectively, thereby accelerating stochastic algorithms with ease. 
  	
  	We take two prevalent stochastic frameworks, MB-SARAH and mS2GD, to verify their numerical efficiency. For MB-SARAH-RHBB/RHBB+ and mS2GD-RHBB/RHBB+, we rigorously analyze the adaptive acceleration mechanism and evaluate the corresponding complexity. Comprehensive tuning guidelines are provided for reference in practical implementations. We prove that they're both robust in ill-conditioned scenarios. Due to the flexibility, we can obtain different performance by trading-off related parameters.
  	
  	Numerical experiments have been conducted to present the properties of our four algorithms. Massive comparisons have been made in all-round aspects and demonstrate their superiority in modern stochastic optimization. Extensive explorations for the iterative adaptor show its promising scalability.

\paragraph{\bf Acknowledgments: } This work is supported partly by National Key R\&D Program of China 2021YFA1003600, 2021YFA1003601, and the State Key Program of National Natural Science of China (12331011) and National Natural Science Foundation of China (12071398). 
  	\bibliographystyle{abbrv}
  	\bibliography{ref}

\newpage
\appendix
\section*{\centerline{Appendix}}

\paragraph{A. Proof of Lemma \ref{lemma1}}

According to the strong convexity (\ref{C}) and the smoothness (\ref{L}) of $P(\cdot)$, we obtain an upper boundary for RHBB step size 
  	$$\begin{aligned}
  	(\eta_k^s)^{\mathrm{RHBB}}
  	& = \frac{\gamma}{\max\{|S_1|, |S_2|\}}\cdot\left( \frac{\alpha^{h(\sigma_1s+\sigma_2k)}\cdot\left\|w_k^s-w_{k-1}^s\right\|^2}{\left(\left(w_k^s-w_{k-1}^s\right)^T\left(\nabla P_{S_1}\left(w_k^s\right)-\nabla P_{S_1}\left(w_{k-1}^s\right)\right)\right)}\right.\\
  	&+ \left.\frac{\left(1-\alpha^{h(\sigma_1s+\sigma_2k)}\right)\cdot\left(\left(w_k^s-w_{k-1}^s\right)^T\left(\nabla P_{S_2}\left(w_k^s\right)-\nabla P_{S_2}\left(w_{k-1}^s\right)\right)\right)}{\left\|\nabla P_{S_2}\left(w_k^s\right)-\nabla P_{S_2}\left(w_{k-1}^s\right)\right\|^2}\right)\\
  	& \leq \frac{\gamma}{\overline{b}}\cdot \left( \frac{\hat{\alpha}\cdot\left\|w_k^s-w_{k-1}^s\right\|^2}{\mu\left\|w_k^s-w_{k-1}^s\right\|^2}+\frac{(1-\tilde{\alpha})\cdot\left\|\nabla P_{S_2}\left(w_k^s\right)-\nabla P_{S_2}\left(w_{k-1}^s\right)\right\|^2}{L\left\|\nabla P_{S_2}\left(w_k^s\right)-\nabla P_{S_2}\left(w_{k-1}^s\right)\right\|^2}\right)\\
  	& = \frac{\gamma}{\overline{b}}\cdot\frac{\hat{\alpha} L+(1-\tilde{\alpha})\mu}{\mu L}.\\
  	\end{aligned}$$
  	By the definition of $P^+(\cdot)$ (\ref{P+}) and $L_q$, $\mu_q$ (\ref{L_u_Q}), the individual $L$-smoothness of $f_i(\cdot)$ implies the uniform $L_q$-smoothness of $P^+(\cdot)$, we achieve
  	$$\begin{aligned}
  	(\eta_k^s)^{\mathrm{RHBB+}}& = \frac{\gamma}{\max\{|S_1|, |S_2|\}}\cdot\left( \frac{\alpha^{h(\sigma_1s+\sigma_2k)}\cdot\left\|w_k^s-w_{k-1}^s\right\|^2}{\left(\left(w_k^s-w_{k-1}^s\right)^T\left(\nabla P_{S_1}^+\left(w_k^s\right)-\nabla P_{S_1}^+\left(w_{k-1}^s\right)\right)\right)}\right.\\
  	&+ \left.\frac{\left(1-\alpha^{h(\sigma_1s+\sigma_2k)}\right)\cdot\left(\left(w_k^s-w_{k-1}^s\right)^T\left(\nabla P_{S_2}^+\left(w_k^s\right)-\nabla P_{S_2}^+\left(w_{k-1}^s\right)\right)\right)}{\left\|\nabla P_{S_2}^+\left(w_k^s\right)-\nabla P_{S_2}^+\left(w_{k-1}^s\right)\right\|^2}\right)\\
  	& \leq \frac{\gamma}{\overline{b}}\cdot\left(\frac{\hat{\alpha}\cdot\left\|w_k^s-w_{k-1}^s\right\|^2}{\frac{1}{|S_1|} \sum_{i \in S_1} \frac{\mu}{{n}q_i}}+\frac{\left(1-\tilde{\alpha}\right)\cdot\left\|\nabla P^+_{S_2}\left(w_k^s\right)-\nabla P^+_{S_2}\left(w_{k-1}^s\right)\right\|^2}{L_q\left\|\nabla P^+_{S_2}\left(w_k^s\right)-\nabla P^+_{S_2}\left(w_{k-1}^s\right)\right\|^2}\right)\\
  	& = \frac{\gamma}{\overline{b}}\cdot\frac{\hat{\alpha} L_q+(1-\tilde{\alpha})\mu_q}{\mu_q L_q}.\\
  	\end{aligned}$$
  	On the basis of (\ref{inequation1}), we have
  	$$\mathbb{E}\left[P\left(w_{k+1}^s\right)\right]\leq \mathbb{E}\left[P\left(w_k^s\right)\right]-\eta_k^s \mathbb{E}\left[\nabla P\left(w_k^s\right)^{\top} v_k^s\right]+\frac{L (\eta_k^s)^2}{2} \mathbb{E}\left[\left\|v_k^s\right\|^2\right].$$
  	Due to the fact $\theta_1^T\theta_2 = \frac{1}{2}\left[\|\theta_1\|^2+\|\theta_2\|^2-\|\theta_1-\theta_2\|^2\right]$, substituting the related boundary of RHBB step size, we obtain
  	$$
  	\begin{aligned}
  	& \mathbb{E}\left[P\left(w_{k+1}^s\right)\right]\\ 
  	& \leq \mathbb{E}\left[P\left(w_k^s\right)\right]-\frac{\gamma\hat{\alpha} L+ \gamma(1-\tilde{\alpha})\mu}{\overline{b}\mu L} \mathbb{E}\left[\nabla P\left(w_k^s\right)^{\top} v_k^s\right]+\frac{L \gamma^2}{2 \overline{b}^2} \cdot \left(\frac{\hat{\alpha} L+ (1-\tilde{\alpha})\mu}{\mu L}\right)^2\mathbb{E}\left[\left\|v_k^s\right\|^2\right] \\
  	& =\mathbb{E}\left[P\left(w_k^s\right)\right]-\frac{\gamma\hat{\alpha} L+ \gamma(1-\tilde{\alpha})\mu}{2 \overline{b}\mu L} \mathbb{E}\left[\left\|\nabla P\left(w_k^s\right)\right\|^2\right]+\frac{\gamma\hat{\alpha} L+ \gamma(1-\tilde{\alpha})\mu}{2 \overline{b}\mu L}\mathbb{E}\left[\| \nabla P\left(w_k^s\right)\right.\left.-v_k^s \|^2\right]\\
  	&-\frac{\hat{\alpha} L+ (1-\tilde{\alpha})\mu}{\mu L}\cdot\left(\frac{\gamma}{2 \overline{b}}-\frac{L \gamma^2}{2 \overline{b}^2}\cdot\frac{\hat{\alpha} L+ (1-\tilde{\alpha})\mu}{\mu L}\right) \mathbb{E}\left[\left\|v_k^s\right\|^2\right].
  	\end{aligned}
  	$$
  	Adding up $k$ from $0$ to $m$, we have
  	$$\begin{aligned} \mathbb{E}\left[P\left(w_{m+1}^s\right)\right]  & \leq  \mathbb{E}\left[P\left(w_0^s\right)\right]-\sum_{k=0}^m\frac{\gamma}{2\overline{b}}\cdot\frac{ \hat{\alpha} L +(1-\tilde{\alpha})\mu}{\mu L} \mathbb{E}\left[\left\|\nabla P\left(w_k^s\right)\right\|^2\right]\\ & +\sum_{k=0}^m\frac{\gamma}{2\overline{b}}\cdot\frac{\hat{\alpha} L+(1-\tilde{\alpha})\mu}{\mu L} \mathbb{E}\left[\| \nabla P\left(w_k^s\right)\right.  \left.-v_k^s \|^2\right]\\ &- \sum_{k=0}^m\frac{\gamma}{2\overline{b}}\frac{\hat{\alpha} L+(1-\tilde{\alpha})\mu}{\mu L}\left(1-\frac{L \gamma}{ \overline{b}}\cdot\frac{\hat{\alpha} L+(1-\tilde{\alpha})\mu}{\mu L}\right) \mathbb{E}\left[\left\|v_k^s\right\|^2\right].\end{aligned}$$
  	Since $w_*=\arg \min _w P(w)$, we ascertain that
  	$$\begin{aligned}\sum_{k=0}^m \mathbb{E}\left[\left\|\nabla P\left(w_k^s\right)\right\|^2\right] &\leq \frac{2 \overline{b}\mu L}{\gamma\hat{\alpha} L+\gamma(1-\tilde{\alpha})\mu} \mathbb{E}\left[P\left(w_0\right)-P\left(w_*\right)\right]+\sum_{k=0}^m \mathbb{E}\left[\left\|\nabla P\left(w_k^s\right)-v_k^s\right\|^2\right] \\ & -\left(1-\frac{L \gamma}{\overline{b}}\cdot\frac{\hat{\alpha} L+\left(1-\tilde{\alpha}\right)\mu}{\mu L}\right) \sum_{k=0}^m \mathbb{E}\left[\left\|v_k^s\right\|^2\right].\end{aligned}$$
  	
  	By using RHBB+ and the corresponding boundary, the remaining parts of Lemma \ref{lemma1} can be proven in a parallel manner. We will no longer expand in detail.

\vskip 0.5cm

\paragraph{B. Proof of Lemma \ref{lemma2}}
~\\

Based on Lemma 3 in \cite{nguyen2}, we readily obtain
  	$$ 
  	\mathbb{E}\left[\left\|\nabla P\left(w_k^s\right)-v_k^s\right\|^2\right] \leq \frac{1}{b}\left(\frac{n-b}{n-1}\right) L^2  \sum_{j=1}^k \left(\eta_{j-1}^s\right)^2 \mathbb{E}\left[\left\|v_{j-1}^s\right\|^2\right].
  	$$
  	By replacing $\eta_{j}^s$ (the step size) with related boundaries, we complete the proof.
\vskip 0.5cm

\paragraph{C. Proof of Theorem \ref{theorem1}}
~\\

Since $\left\|\nabla P\left(w_0^s\right)-v_0^s\right\|^2=0$, we apply (\ref{lemma21}) in Lemma \ref{lemma2} and sum over $k=0, ..., m$ to obtain
  	$$\begin{aligned}  \sum_{k=0}^m \mathbb{E}\left[\left\|\nabla P\left(w_k^s\right)-v_k^s\right\|^2\right] & \leq \frac{\left(\hat{\alpha}\gamma L^2+\left(1-\tilde{\alpha}\right)\gamma\mu L\right)^2}{b \overline{b}^2 \mu^2 L^2}\cdot \left(\frac{n-b}{n-1}\right)\\ & \cdot\left(m \mathbb{E}\left[\left\|v_0^s\right\|^2\right]+(m-1) \mathbb{E}\left[\left\|v_1^s\right\|^2\right]+\ldots+\mathbb{E}\left[\left\|v_{m-1}^s\right\|^2\right]\right).\end{aligned}$$
    Parameters $b$, $\gamma$ are chosen such that
    \begin{equation}
  		\label{theorem1i}
  		\frac{m\left(n-b\right)}{b\left(n-1\right)}\left(\frac{\hat{\alpha}\gamma L^2+\left(1-\tilde{\alpha}\right)\gamma\mu L}{\mu L\overline{b}}\right)^2+\frac{\hat{\alpha} \gamma L+(1-\tilde{\alpha}) \gamma \mu}{\mu \overline{b}} \leq 1.
  		\end{equation}
    Plugging (\ref{theorem1i}) in, we hence have
  	$$\begin{aligned} & \sum_{k=0}^m \mathbb{E}\left[\left\|\nabla P\left(w_k^s\right)-v_k^s\right\|^2\right]-\left(1-\frac{L \gamma}{\overline{b}}\frac{\hat{\alpha} L+(1-\tilde{\alpha})\mu}{\mu L}\right) \sum_{k=0}^m \mathbb{E}\left[\left\|v_k\right\|^2\right] \\ & \leq	\left(\left(\frac{\hat{\alpha}L+\left(1-\tilde{\alpha}\right)\mu}{\mu L}\right)^2\cdot\frac{L^2 \gamma^2}{b \overline{b}^2}\left(\frac{n-b}{n-1}\right) m-\left(1-\frac{\hat{\alpha} \gamma L+(1-\tilde{\alpha}) \gamma \mu)}{\mu \overline{b}}\right)\right)  \\ &
  	\times \left(\sum_{k=1}^m\mathbb{E}\left[\left\|v_{k-1}^s\right\|^2\right]\right) \leq 0. \end{aligned}$$
  	Using the Lemma \ref{lemma1}, we further derive
  	$$
  	\begin{aligned}
  	& \sum_{k=0}^m \mathbb{E}\left[\left\|\nabla P\left(w_k^s\right)\right\|^2\right] \leq \frac{2 \mu \overline{b} L}{\hat{\alpha} \gamma L+(1-\tilde{\alpha})\gamma\mu} \mathbb{E}\left[P\left(w_0^s\right)-P\left(w_*\right)\right] \\
  	& +\sum_{k=0}^m \mathbb{E}\left[\left\|\nabla P\left(w_k^s\right)-v_k^s\right\|^2\right]-\left(1-\frac{\hat{\alpha} \gamma L^2+(1-\tilde{\alpha}) \gamma L \mu}{\mu \overline{b} L}\right) \sum_{k=0}^m \mathbb{E}\left[\left\|v_k^s\right\|^2\right].
  	\end{aligned}
  	$$
  	By the definition of $\widetilde{w}_s$ and the outer update rule $\widetilde{w}_s=w_m^s$, we ascertain
  	$$\begin{aligned} \mathbb{E}\left[\left\|\nabla P\left(w_m^s\right)\right\|^2\right] & =\frac{1}{m+1} \sum_{k=0}^m \mathbb{E}\left[\left\|\nabla P\left(w_k^s\right)\right\|^2\right] \\ & \leq \frac{2 \overline{b} \mu L}{\gamma(m+1)(\hat{\alpha} L+(1-\tilde{\alpha})\mu)} \mathbb{E}\left[P\left(w_0^s\right)-P\left(w_*\right)\right].\end{aligned}$$
  	
  	By substituting with RHBB+ step sizes and the corresponding upper boundary, the remaining parts of Theorem \ref{theorem1} can be proven similarly in parallel. Notably here, parameters $b$, $\gamma$ are chosen such that
    \begin{equation}
  		\label{theorem1ii}
  		\frac{mL_r^2\left(n-b\right)}{b\left(n-1\right)}\left(\frac{\hat{\alpha}\gamma L_q+\left(1-\tilde{\alpha}\right)\gamma\mu_q}{\mu_q \overline{b}}\right)^2+L_r\frac{\hat{\alpha} \gamma  L_q+(1-\tilde{\alpha}) \gamma  \mu_q}{\mu_q \overline{b} } \leq 1,
  \end{equation}

\vskip 0.5cm

\paragraph{D. Proof of Theorem \ref{theorem2}}
~\\

Since $w_0^s=\widetilde{w}_{s-1}$ and $\widetilde{w}_s=w_m^s$, we apply the Theorem \ref{theorem1} and have
  	$$
  	\begin{aligned}
  	\mathbb{E}\left[\left\|\nabla P\left(\widetilde{w}_s\right) \mid \widetilde{w}_{s-1}\right\|^2\right] & =\mathbb{E}\left[\left\|\nabla P\left(\widetilde{w}_s\right) \mid w_0^s\right\|^2\right] \\
  	& \leq \frac{2 \overline{b}}{\gamma(m+1)}\cdot\frac{\mu L}{\hat{\alpha} L+(1-\tilde{\alpha})\mu} \mathbb{E}\left[P\left(w_0^s\right)-P\left(w_*\right)\right]. \\
  	\end{aligned}
  	$$
  	By taking expectation and using the convexity (\ref{convex1}), we obtain
  	$$
  	\begin{aligned}
  	\mathbb{E}\left[\left\|\nabla P\left(\widetilde{w}_s\right)\right\|^2\right] & \leq \frac{\overline{b} L}{\hat{\alpha} \gamma (m+1)L+(1-\tilde{\alpha})\gamma (m+1)\mu} \mathbb{E}\left[\left\|\nabla P\left(\widetilde{w}_{s-1}\right)\right\|^2\right] \\
  	& \leq\left[\frac{\overline{b} L}{\hat{\alpha} \gamma (m+1)L+(1-\tilde{\alpha})\gamma (m+1)\mu}\right]^s\left\|\nabla P\left(\widetilde{w}_0\right)\right\|^2.
  	\end{aligned}
  	$$

  	By substituting with RHBB+ step sizes and the corresponding upper boundary, the remaining parts of Theorem \ref{theorem2} follow a similar line of reasoning.

\vskip 0.5cm

\paragraph{E. Proof of Lemma \ref{lemma3}}
~\\

    Before the formal proof, let us define the $j$-th estimate at $w_k^s$ as $\tilde{v}_j=\nabla f_j\left(w_k^s\right)-\nabla f_j(\widetilde{w}_{s-1})+\nabla P(\widetilde{w}_{s-1})$, where $\nabla f_j$ represents the gradient of the $j$-th component function. According to $\tilde{v}_k^s=\frac{1}{b}\sum_{j \in S} \tilde{v}_j$, we obtain
  	$$
  	\begin{aligned}
  	\mathbb{E}\left[\left\|\tilde{v}_k^s\right\|^2\right]& =\frac{1}{b^2} \mathbb{E}\left[\|\sum_{j \in S} \tilde{v}_j\|^2\right] \\
  	& =\frac{1}{b^2} \mathbb{E}\left[\|\sum_{j \in S^{\prime}} \tilde{v}_j\|^2+2(\sum_{j \in S^{\prime}} \tilde{v}_j)^T (\tilde{v}_{j \in S-S^{\prime}})+\|\tilde{v}_{j \in S-S^{\prime}}\|^2\right] \\
  	& =\frac{1}{b^2}\left[\mathbb{E}[\|\sum_{j \in S^{\prime}} \tilde{v}_j\|^2]+2\|\nabla P\left(w_{k-1}^s\right)\|^2 +\mathbb{E}[\|\tilde{v}_{j \in S-S^{\prime}}\|^2]\right] \\
  	& =\cdots\\
  	& =\frac{1}{b^2}\left[\sum_{j \in S} \mathbb{E}\left[\|\tilde{v}_j\|^2\right]+2\left(b-1\right)\|\nabla P\left(w_{k-1}^s\right)\|^2\right] \\
  	& \leq \frac{1}{b^2}\left[\sum_{j \in S} \mathbb{E}\left[\|\tilde{v}_j\|^2\right]+2 b\|\nabla P\left(w_{k-1}^s\right)\|^2\right] \\
  	& \leq \frac{4 L}{b}\left[P\left(w_{k-1}^s\right)-P\left(w_*\right)+P(\widetilde{w}_{s-1})-P\left(w_*\right)\right]+\frac{2}{b}\|\nabla P\left(w_{k-1}^s\right)\|^2, \\
  	\end{aligned}
  	$$
  	where the subset $S^{\prime} \subset S$ with the
  	number of members of $|S-S^{\prime}|=1$. The last equality follows   Lemma 3 in \cite{yang2}.

\vskip 0.5cm

\paragraph{F. Proof of Theorem \ref{theorem3}}
~\\

By Lemma \ref{lemma3} and $\mathbb{E}\left[\tilde{v}_{k-1}^s\right]=\nabla P\left(w_{k-1}^s\right)$, we obtain 
  	$$\begin{aligned} &\mathbb{E}\left[\left\|w_k^s-w_*\right\|^2\right]\\
  	& = \mathbb{E}\left[\left\|w_{k-1}^s-\eta_{k-1}^s \tilde{v}_{k-1}^s-w_*\right\|_2^2\right]\\
  	& =  \left\|w_{k-1}^s-w_*\right\|_2^2-2 \tilde{\eta}_{k-1}^s \mathbb{E}\left[\left(w_{k-1}^s-w_*\right)^T \tilde{v}_{k-1}^s\right] +(\tilde{\eta}_{k-1}^s)^2 \mathbb{E}\left[\left\|\tilde{v}_{k-1}^s\right\|^2\right] \\ 
  	& \leq\left\|w_{k-1}^s-w_*\right\|^2-2 \tilde{\eta}_{k-1}^s\left(w_{k-1}^s-w_*\right)^T \nabla P\left(w_{k-1}^s\right) \\
  	&+\frac{4 L (\tilde{\eta}_{k-1}^s)^2}{b}\left[P\left(w_{k-1}^s\right)\right.   \left.-P\left(w_*\right)+P(\widetilde{w}_{s-1})-P\left(w_*\right)\right]+\frac{2 (\tilde{\eta}_{k-1}^s)^2}{b}\left\|\nabla P  \left(w_{k-1}^s\right)\right\|^2 \\  
  	& \leq \left\|w_{k-1}^s-w_*\right\|^2-2 \tilde{\eta}_{k-1}^s\left[P\left(w_{k-1}^s\right)-P\left(w_*\right)\right] \\ 
  	&+ \frac{4 L (\tilde{\eta}^s_{k-1})^2}{b}\left[P\left(w_{k-1}^s\right)-P\left(w_*\right)\right. \left.+P(\widetilde{w}_{s-1})-P\left(w_*\right)\right]+\frac{2 (\tilde{\eta}_{k-1}^s)^2}{b}\left\|\nabla P\left(w_{k-1}^s\right)\right\|^2 \\  
  	&\leq \left\|w_{k-1}^s-w_*\right\|^2-2 \tilde{\eta}_{k-1}^s\left(1-\frac{4 L \tilde{\eta}_{k-1}^s}{b}\right) \left[P\left(w_{k-1}^s\right)-P\left(w_*\right)\right]\\  &+\frac{4 L (\tilde{\eta}^s_{k-1})^2}{b} \cdot\left[P(\widetilde{w}_{s-1})-P\left(w_*\right)\right],
  	\end{aligned}$$
  	where we use the convexity of $P(\cdot)$ in the second inequality and (\ref{L2}) in the last.

  	We derive the upper boundary for RHBB step size in mS2GD as follows
  	$$\begin{aligned}
  	(\tilde{\eta}_k^s)^{\mathrm{RHBB}}& = \frac{\gamma_2}{\max\{|S_1|, |S_2|\}}\cdot\left( \frac{\alpha^{h(\sigma_1s+\sigma_2k)}\cdot\left\|w_k^s-w_{k-1}^s\right\|^2}{\left(\left(w_k^s-w_{k-1}^s\right)^T\left(\nabla P_{S_1}\left(w_k^s\right)-\nabla P_{S_1}\left(w_{k-1}^s\right)\right)\right)}\right.\\
  	&+ \left.\frac{\left(1-\alpha^{h(\sigma_1s+\sigma_2k)}\right)\cdot\left(\left(w_k^s-w_{k-1}^s\right)^T\left(\nabla P_{S_2}\left(w_k^s\right)-\nabla P_{S_2}\left(w_{k-1}^s\right)\right)\right)}{\left\|\nabla P_{S_2}\left(w_k^s\right)-\nabla P_{S_2}\left(w_{k-1}^s\right)\right\|^2}\right)\\
  	& \leq \frac{\gamma_2}{\overline{b}}\cdot \left( \frac{\hat{\alpha}\cdot\left\|w_k^s-w_{k-1}^s\right\|^2}{\mu\left\|w_k^s-w_{k-1}^s\right\|^2}+\frac{(1-\tilde{\alpha})\cdot\left\|\nabla P_{S_2}\left(w_k^s\right)-\nabla P_{S_2}\left(w_{k-1}^s\right)\right\|^2}{L\left\|\nabla P_{S_2}\left(w_k^s\right)-\nabla P_{S_2}\left(w_{k-1}^s\right)\right\|^2}\right)\\
  	& = \frac{\gamma_2}{\overline{b}}\cdot\frac{\hat{\alpha} L+(1-\tilde{\alpha})\mu}{\mu L}.
  	\end{aligned}$$
   Plugging it in, we have
  	$$\begin{aligned} \mathbb{E}\left\|w_k^s-w_*\right\|^2 & \leq  \left\|w_{k-1}^s-w_*\right\|^2 \\&-\frac{2\hat{\alpha}\gamma_2 L+2(1-\tilde{\alpha})\gamma_2\mu}{\mu L\overline{b}}\left(1-\frac{4 L(\hat{\alpha}\gamma_2 L+(1-\tilde{\alpha})\gamma_2\mu)}{b \overline{b}\mu L}\right)\left[P\left(w_{k-1}^s\right)-P\left(w_*\right)\right] \\ & +\frac{4 L}{b \overline{b}^2}\left(\frac{\hat{\alpha}\gamma_2 L+(1-\tilde{\alpha})\gamma_2\mu}{\mu L}\right)^2\left[P(\widetilde{w}_{s-1})-P\left(w_*\right)\right].\end{aligned}$$
  	By the definition of $\widetilde{w}_{s-1}$ in mS2GD-RHBB, we have (see in \cite{yang}\cite{konevcny})
  	$$
  	\mathbb{E}\left[P\left(\widetilde{w}_{s}\right)\right]=\frac{1}{m} \sum_{k=1}^m \mathbb{E}\left[P\left(w_k^s\right)\right].
  	$$
  	By summing over the previous inequality over $k$, we take expectation conditioned on history randomness. Since $\widetilde{w}_s=w_m^s$, $w_0^s=\widetilde{w}_{s-1}$, we obtain
  	$$\begin{aligned} & \mathbb{E}\left\|w_m^s-w_*\right\|^2+\frac{2m\hat{\alpha} \gamma_2L+2m(1-\tilde{\alpha})\gamma_2\mu}{\mu L\overline{b}}\left(1-\frac{4 L(\hat{\alpha}\gamma_2 L+(1-\tilde{\alpha})\gamma_2\mu)}{b \overline{b}\mu L}\right)\mathbb{E}\left[P\left(\widetilde{w}_s\right)-P\left(w_*\right)\right] \\ 
  	& \leq\mathbb{E}\left\|w_0^s-w_*\right\|^2+\frac{4m L}{b \overline{b}^2}\left(\frac{\hat{\alpha}\gamma_2 L+(1-\tilde{\alpha})\gamma_2\mu}{\mu L}\right)^2 \mathbb{E}\left[P(\widetilde{w}_{s-1})-P\left(w_*\right)\right]\\
  	& =\mathbb{E}\left\|\widetilde{w}_{s-1}-w_*\right\|_2^2+\frac{4m L}{b \overline{b}^2}(\frac{\hat{\alpha}\gamma_2 L+(1-\tilde{\alpha})\gamma_2\mu}{\mu L})^2 \mathbb{E}\left[P(\widetilde{w}_{s-1})-P\left(w_*\right)\right].\end{aligned}$$
  	Employing the strong convexity (\ref{convex1}), we further attain
  	$$ \begin{aligned} &\mathbb{E}\left\|\widetilde{w}_s-w_*\right\|_2^2+\frac{4m L}{b \overline{b}^2}\left(\frac{\hat{\alpha}\gamma_2 L+(1-\tilde{\alpha})\gamma_2\mu}{\mu L}\right)^2 \mathbb{E}\left[P(\widetilde{w}_{s})-P\left(w_*\right)\right] \\ & \leq \frac{2}{\mu} \mathbb{E}\left[P(\widetilde{w}_{s-1})-P\left(w_*\right)\right]+\frac{4m L}{b \overline{b}^2}\left(\frac{\hat{\alpha}\gamma_2 L+(1-\tilde{\alpha})\gamma_2\mu}{\mu L}\right)^2 \mathbb{E}\left[P(\widetilde{w}_{s-1})-P\left(w_*\right)\right] \\ & =
  	\left(\frac{2}{\mu}+\frac{4m L}{b \overline{b}^2}(\frac{\hat{\alpha}\gamma_2 L+(1-\tilde{\alpha})\gamma_2\mu}{\mu L})^2\right)\mathbb{E}\left[P(\widetilde{w}_{s-1})-P\left(w_*\right)\right].
  	\end{aligned}$$
  	By the definition of $\kappa_r$, we at last simplify it into
  	$$\mathbb{E}\left[P\left(\widetilde{w}_s\right)-P\left(w_*\right)\right] \leq \left(\frac{\kappa}{\gamma_2\kappa_r}\cdot\frac{b \overline{b}^2}{m(b \overline{b}-4\gamma_2\kappa_r)}+\frac{2\gamma_2 \kappa_r}{b \overline{b}-4\gamma_2\kappa_r}\right)\mathbb{E}\left[P\left(\widetilde{w}_{s-1}\right)-P\left(w_*\right)\right].$$
  	By recursively applying the previous procedures, we derive
  	$$\mathbb{E}\left[P\left(\widetilde{w}_s\right)-P\left(w_*\right)\right] \leq \left(\frac{\kappa}{\gamma_2\kappa_r}\cdot\frac{b \overline{b}^2}{m(b \overline{b}-4\gamma_2\kappa_r)}+\frac{2\gamma_2 \kappa_r}{b \overline{b}-4\gamma_2\kappa_r}\right)^s\mathbb{E}\left[P\left(\widetilde{w}_{0}\right)-P\left(w_*\right)\right].$$
  	
  	\vspace{8pt}
  	By substituting with RHBB+ step sizes and the corresponding boundary, the remaining parts in Theorem \ref{theorem3} can be proven similarly in parallel. We supply the relational boundary of RHBB+ in mS2GD as follows
  	$$\begin{aligned}
  	(\tilde{\eta}_k^s)^{\mathrm{RHBB+}}& = \frac{\gamma_2}{\max\{|S_1|, |S_2|\}}\cdot\left( \frac{\alpha^{h(\sigma_1s+\sigma_2k)}\cdot\left\|w_k^s-w_{k-1}^s\right\|^2}{\left(\left(w_k^s-w_{k-1}^s\right)^T\left(\nabla P_{S_1}^+\left(w_k^s\right)-\nabla P_{S_1}^+\left(w_{k-1}^s\right)\right)\right)}\right.\\
  	&+ \left.\frac{\left(1-\alpha^{h(\sigma_1s+\sigma_2k)}\right)\cdot\left(\left(w_k^s-w_{k-1}^s\right)^T\left(\nabla P_{S_2}^+\left(w_k^s\right)-\nabla P_{S_2}^+\left(w_{k-1}^s\right)\right)\right)}{\left\|\nabla P_{S_2}^+\left(w_k^s\right)-\nabla P_{S_2}^+\left(w_{k-1}^s\right)\right\|^2}\right)\\
  	& \leq \frac{\gamma_2}{\overline{b}}\cdot\left(\frac{\hat{\alpha}\cdot\left\|w_k^s-w_{k-1}^s\right\|^2}{\frac{1}{|S_1|} \sum_{i \in S_1} \frac{\mu}{{n}q_i}}+\frac{\left(1-\tilde{\alpha}\right)\cdot\left\|\nabla P^+_{S_2}\left(w_k^s\right)-\nabla P^+_{S_2}\left(w_{k-1}^s\right)\right\|^2}{L_q\left\|\nabla P^+_{S_2}\left(w_k^s\right)-\nabla P^+_{S_2}\left(w_{k-1}^s\right)\right\|^2}\right)\\
  	& = \frac{\gamma_2}{\overline{b}}\cdot\frac{\hat{\alpha} L_q+(1-\tilde{\alpha})\mu_q}{\mu_q L_q}.\\
  	\end{aligned}$$

  \end{document}